\documentclass{article}
\usepackage[textheight=9in,
textwidth=6.5in,
top=1in,
headheight=14pt,
headsep=25pt,
footskip=30pt]{geometry}
\usepackage{hyperref}
\usepackage{amssymb}
\usepackage{amsfonts}
\usepackage{amsmath}
\usepackage{eurosym}
\usepackage{amssymb}
\usepackage{amsfonts}
\usepackage{amsmath}
\usepackage{graphicx}
\usepackage[ruled]{algorithm2e}
\usepackage{subcaption}
\usepackage{xcolor}
\usepackage{titlesec}
\usepackage[titletoc,toc,title]{appendix}

\hypersetup{
colorlinks=true,
linkcolor=black,
anchorcolor=black,
citecolor=black,
urlcolor=black,
pdftitle={Multiple Bayesian Filtering as Message Passing},
pdfauthor={Giorgio M. Vitetta, Pasquale Di Viesti, Emilio Sirignano and Francesco Montorsi},
pdfkeywords={Hidden Markov Model, Factor Graph, Particle Filter, Sum-Product Algorithm, Marginalized Particle Filter, Kalman Filter, Multiple Particle Filtering}}
\graphicspath{{Figures/}}

\begin{document}

\title{Multiple Bayesian Filtering as Message Passing}
\author{}
\maketitle

\begin{abstract}
In this manuscript, a general method for deriving filtering algorithms that
involve a network of interconnected Bayesian filters is proposed. This
method is based on the idea that the processing accomplished inside each of
the Bayesian filters and the interactions between them can be represented as
message passing algorithms over a proper graphical model. The usefulness of
our method is exemplified by developing new filtering techniques, based on
the interconnection of a particle filter and an extended Kalman filter, for
conditionally linear Gaussian systems. Numerical results for two specific
dynamic systems evidence that the devised algorithms can achieve a better
complexity-accuracy tradeoff than marginalized particle filtering and
multiple particle filtering.
\end{abstract}

\begin{center}
\begin{tabular}{cc}
Giorgio M. Vitetta$^\dagger$ & Pasquale Di Viesti$^\dagger$ \\ 
\vspace{.5cm} \href{mailto:giorgio.vitetta@unimore.it}{%
giorgio.vitetta@unimore.it} & \href{mailto:pasquale.diviesti@unimore.it}{%
pasquale.diviesti@unimore.it} \\ 
Emilio Sirignano$^\dagger$ & Francesco Montorsi \\ 
\href{mailto:emilio.sirignano@unimore.it}{emilio.sirignano@unimore.it} & 
\href{mailto:francesco.montorsi@gmail.com}{francesco.montorsi@gmail.com} \\ 
& 
\end{tabular}

\bigskip $^\dagger${Dept. of Engineering "Enzo Ferrari", University of
Modena and Reggio Emilia.}
\end{center}

\bigskip


\bigskip

\section{Introduction\label{sec:intro}}

It is well known that \emph{Bayesian filtering} represents a general
recursive solution to the \emph{nonlinear filtering problem }(e.g., see \cite%
[Sect. II, eqs. (3)-(5)]{Arulampalam_2002}), i.e. to the problem of
inferring the posterior distribution of the hidden state of a nonlinear 
\emph{state-space model} (SSM). Unluckily, this solution can be put in
closed form in few cases \cite{Anderson_1979}. For this reason, various
filtering methods generating a functional approximation of the desired
posterior pdf have been developed; these can be divided into \emph{local}
and \emph{global} methods on the basis of the way the posterior \emph{%
probability density function} (pdf) is approximated \cite{Mazuelas_2013}, 
\cite{Smidl_2008}. On the one hand, local techniques, like \emph{extended
Kalman filtering} (EKF) \cite{Anderson_1979}, are computationally efficient,
but may suffer from error accumulation over time; on the other hand, global
techniques, like \emph{particle filtering} (PF) \cite{Doucet_2001}--\cite%
{Andrieu_2002}, may achieve high accuracy at the price, however, of
unacceptable complexity and numerical problems when the dimension of the
state space becomes large \cite{Daum_2003}--\cite{Djuric_2013}. These
considerations have motivated the investigation of \emph{various methods
able to achieve high accuracy under given computational constraints}. Some
of such solutions are based on the idea of \emph{combining} \emph{local and
global methods}; relevant examples of this approach are represented by: 1) 
\emph{Rao-Blackwellized particle filtering} (RBPF; also known as \emph{%
marginalized particle filtering}) \cite{Schon_2005} and other techniques
related to it (e.g., see \cite{Smidl_2008}); 2) \emph{cascaded architectures}
based on the joint use of EKF and PF (e.g., see \cite{Krach_2008_bis}). Note
that, in the first case, the state vector is split into two disjoint
components, namely, a \emph{linear state component} and a \emph{nonlinear
state component}; moreover, these are estimated by a bank of Kalman filters
and by a particle filter, respectively. In the second case, instead, an
extended Kalman filter and a particle filter are run over partially
overlapped state vectors. In both cases, however, two heterogeneous
filtering methods are combined in a way that the resulting overall algorithm
is forward only and, within each of its recursions, both methods are
executed only once. Another class of solutions, known as \emph{multiple
particle filtering} (MPF), is based on the idea of partitioning the state
vector into multiple substates and running multiple particle filters in
parallel, one on each subspace \cite{Djuric_2013}, \cite{Closas_2012}-\cite%
{Hoteit_2016}. The resulting network of particle filters requires the mutual
exchange of statistical information (in the form of estimates/predictions of
the tracked substates or parametric distributions), so that, within each
filter, the unknown portion of the state vector can be integrated out in
both weight computation and particle propagation. In principle, MPF can be
employed only when the selected substates are separable in the state
equation, even if approximate solutions can be devised to circumvent this
problem \cite{Hoteit_2016}. Moreover, the technical literature about MPF has
raised three interesting technical issues that have received limited
attention until now. The first issue refers to the possibility of coupling
an extended Kalman filter with each particle filter of the network; the
former filter should provide the latter one with the statistical information
required for integrating out the unknown portion of the state vector (see 
\cite[Par. 3.2]{Djuric_2007}). The second one concerns the use of filters
having \emph{partially overlapped} substates (see \cite[Sec.1]{Bugallo_2014}%
). The third (and final) issue, instead, concerns the iterative exchange of
statistical information among the interconnected filters of the network.
Some work related to the first issue can be found in \cite{Chavali_2012},
where the application of MPF to target tracking in a cognitive radar network
has been investigated. In this case, however, the proposed solution is based
on Rao-Blackwellisation; for this reason, each particle filter of the
network is not coupled with a single extended Kalman filter, but with a bank
of Kalman filters. The second issue has not been investigated at all,
whereas limited attention has been paid to the third one; in fact, the last
problem has been investigated only in \cite{Closas_2012}, where a specific
iterative method based on game theory has been developed. The need of
employing iterative methods in MPF has been also explicitly recognised in 
\cite{Hoteit_2016}, but no solution has been developed to meet it.

In this manuscript, we first focus on the general problem of developing
filtering algorithms that involve multiple interconnected Bayesian filters;
these filters are run over distinct (but not necessarily disjoint) subspaces
and can exploit iterative methods in their exchange of statistical
information. The solution devised for this problem (and called \emph{%
multiple Bayesian filtering}, MBF, since it represents a generalisation of
the MPF approach) is based on previous work on the application of \emph{%
factor graph theory} to the filtering and smoothing problems \cite%
{Loeliger_2007}--\cite{Vitetta_2019}. More specifically, we show that: a) a 
\emph{graphical model} can be developed for a network of Bayesian filters by
combining multiple factor graphs, each referring to one of the involved
filters; b) the pdfs computed by all these filters can be represented as
messages passed on such a graphical model. This approach offers various
important advantages. In fact, all the expressions of the passed messages
can be derived by applying the same rule, namely the so called \emph{%
sum-product algorithm} (SPA) \cite{Loeliger_2007}, \cite{Kschischang_2001},
to the graphical model devised for the whole network. Moreover, iterative
algorithms can be developed in a natural fashion once the cycles contained
in this graphical model have been identified and the order according to
which messages are passed on them (i.e., the \emph{message scheduling}) has
been established. The usefulness of our approach is exemplified by mainly
illustrating its application to a network made of two Bayesian filters. More
specifically, we investigate the interconnection of an extended Kalman
filter with a particle filter, and develop two new filtering algorithms
under the assumption that the considered SSM is \emph{conditionally linear
Gaussian} (CLG). Simulation results for two specific SSMs evidence that the
devised algorithms perform similarly or better than RBPF and MPF, but
require a smaller computational effort.

The remaining parts of this manuscript are organized as follows. In Section %
\ref{sec:Factorgraphs}, after introducing factor graph theory and the SPA,
the filtering problem is analysed from a factor graph perspective for a
network of multiple interconnected Bayesian filters. In Section \ref%
{sec:Message-Passing}, the tools illustrated in the previous section are
applied to a network consisting of an extended Kalman filter interconnected
with a particle filter, two new MBF algorithms are derived and their
computational complexity is analysed in detail. The developed MBF algorithms
are compared with EKF and RBPF, in terms of accuracy and execution time, in
Section \ref{num_results}. Finally, some conclusions are offered in Section %
\ref{sec:conc}.

\section{Graphical Modelling for Multiple Bayesian Filtering\label%
{sec:Factorgraphs}}

In this paragraph, we illustrate some basic concepts about factor graphs and
the computation of the messages passed over them. Then, we derive a
graphical model for representing the overall processing accomplished by
multiple interconnected Bayesian filters as a message passing on it.

\subsection{Factor Graphs and the Sum-Product Algorithm\label{GF_and_SPA}}

\begin{figure}[tbp]
\centering
\includegraphics[width=0.5\textwidth]{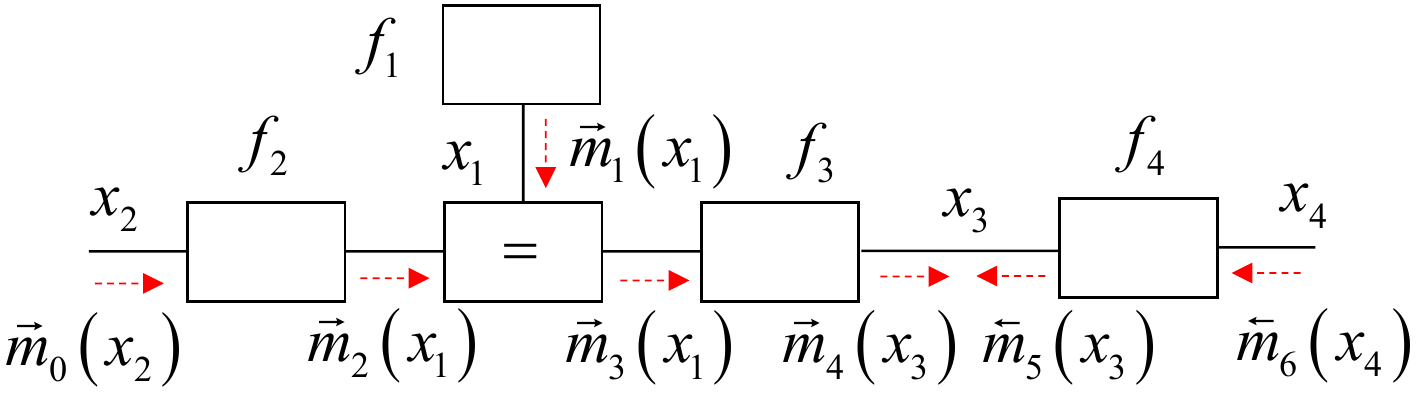}
\caption{Factor graph representing the structure of the function $%
f(x_{1},x_{2},x_{3},x_{4})$ (\protect\ref{factorised_func}) and message
passing on it for the evaluation of the marginal $f(x_{3})$.}
\label{Fig_1}
\end{figure}

A \emph{factor graph} is a graphical model representing the factorization of
any function $f(\cdot )$ expressible as a product of factors $\left\{
f_{i}(\cdot )\right\} $, each depending on a set of variables $\left\{
x_{l}\right\} .$ In the following, Forney-style factor graphs are considered 
\cite{Loeliger_2007}. This means that the factor graph associated with the
function $f(\cdot )$ consists of \emph{nodes}, \emph{edges} (connecting
distinct nodes) and \emph{half-edges} (connected to a single node only).
Moreover, the following rules are employed for its construction: a) every 
\emph{factor} is represented by a single node (a \emph{rectangle} in our
pictures); b) every \emph{variable} is represented by a unique \emph{edge}
or \emph{half edge}; c) the node representing a factor $f_{i}(\cdot )$ is
connected with the edge (or half-edge) representing the variable $x_{l}$ if
and only if such a factor depends on $x_{l}$; d) an \emph{equality
constraint node }(represented by a rectangle labelled by \textquotedblleft
=\textquotedblright ) is used as a branching point when more than two
factors are required to share the same variable. For instance, the
factorisation of the function 
\begin{equation}
f\left( x_{1},x_{2},x_{3},x_{4}\right) =f_{1}\left( x_{1}\right)
\,f_{2}\left( x_{1},x_{2}\right) \,f_{3}\left( x_{1},x_{3}\right)
\,f_{4}\left( x_{3},x_{4}\right)  \label{factorised_func}
\end{equation}%
can be represented through the factor graph shown in Fig. \ref{Fig_1}.

\begin{figure}[tbp]
\centering
\includegraphics[width=0.5\textwidth]{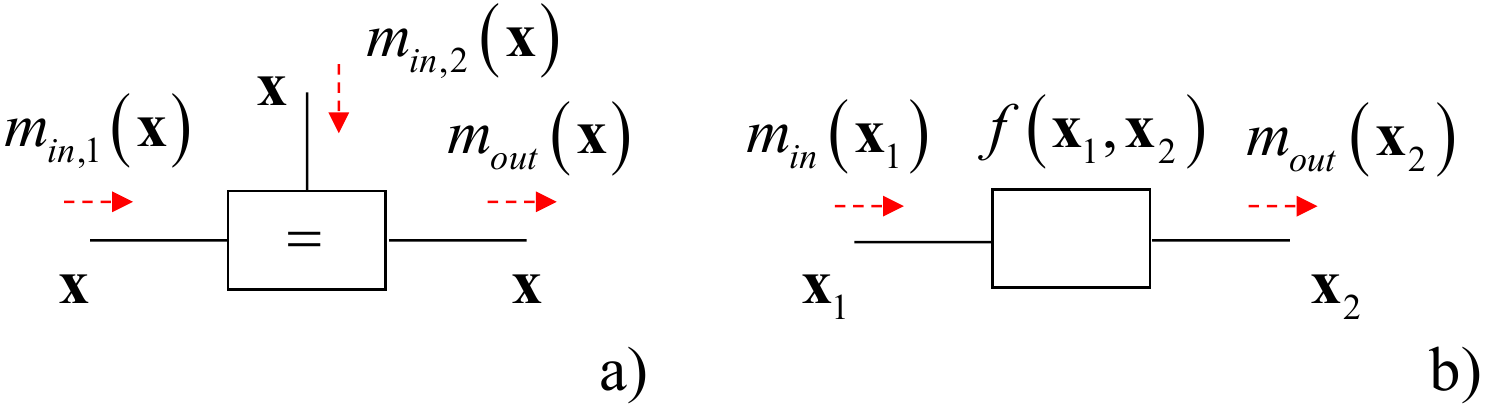}
\caption{Representation of the graphical models which Eqs. (\protect\ref%
{CR_1}) (diagram a)) and (\protect\ref{CR_2}) (diagram b)) refer to.}
\label{Fig_2}
\end{figure}

In this manuscript, factorisable functions represent joint pdfs. It is well
known that the \emph{marginalization} of $f(\cdot )$ with respect to one or
more of its variables can be usually split into a sequence of simpler
marginalizations; our interest in the graph representing $f(\cdot )$ is
motivated by the fact that the function resulting from each of these
marginalizations can be represented as a \emph{message} (conveying a joint
pdf of the variables it depends on) passed along an edge of the graph
itself. In this work, the computation of all the messages is based on the
SPA (also known as \emph{belief propagation}). This algorithm can be
formulated as follows (e.g., see \cite[Sec. IV]{Loeliger_2007}): the message
emerging from a node, representing a factor $f_{i}(\cdot )$, along the edge
associated with a variable $x_{l}$ is expressed by the product of $%
f_{i}(\cdot )$ and the messages along all the incoming edges (except that
associated with $x_{l}$), integrated over all the involved variables except $%
x_{l}$. Two simple applications of the SPA are illustrated in Fig. \ref%
{Fig_2}-a) and in Fig. \ref{Fig_2}-b), that refer to an \emph{equality} 
\emph{constraint node} and to a \emph{function node}, respectively (note
that, generally speaking, these nodes are connected to edges representing 
\emph{vectors of variables}). On the one hand, the message $\vec{m}_{out}(%
\mathbf{x})$ emerging from the \emph{equality node} shown in Fig. \ref{Fig_2}%
-a) is evaluated as 
\begin{equation}
\vec{m}_{out}\left( \mathbf{x}\right) =\vec{m}_{in,1}\left( \mathbf{x}%
\right) \,\vec{m}_{in,2}\left( \mathbf{x}\right) ,  \label{CR_1}
\end{equation}%
where $\vec{m}_{in,1}\left( \mathbf{x}\right) $ and $\vec{m}_{in,2}\left( 
\mathbf{x}\right) $ are the two messages entering the node itself (if a
single message $\vec{m}\left( \mathbf{x}\right) $ enters an equality node,
the two messages emerging from are simply copies of it) and $\mathbf{x}$ is
the vector of variables all these message refer to. On the other hand, the
message $\vec{m}_{out}\left( \mathbf{x}_{2}\right) $ emerging from the \emph{%
function node} shown Fig. \ref{Fig_2}-b), that refers to the function $%
f\left( \mathbf{x}_{1},\mathbf{x}_{2}\right) $ depending on the vectors of
variables $\mathbf{x}_{1}$ and $\mathbf{x}_{2}$, is given by 
\begin{equation}
\vec{m}_{out}\left( \mathbf{x}_{2}\right) =\int \vec{m}_{in}\left( \mathbf{x}%
_{1}\right) \,f\left( \mathbf{x}_{1},\mathbf{x}_{2}\right) \,d\mathbf{x}_{1},
\label{CR_2}
\end{equation}%
where $\vec{m}_{in}\left( \mathbf{x}_{1}\right) $ denotes the message
entering it.

In applying the SPA, it is important to keep in mind that: a) the marginal
pdf $f\left( x_{l}\right) $, referring to the variable $x_{l}$ only, is
expressed by the product of two messages associated with the edge $x_{l}$,
but coming from opposite directions; b) the half-edge associated with a
variable $x_{l}$ may be thought as carrying a constant message of unit value
as incoming message; c) if a marginal pdf is required to be known up to a
scale factor, the involved messages can be freely scaled in their
computation. The use of the last rules and of those expressed by Eqs. (\ref%
{CR_1}) and (\ref{CR_2}) can be exemplified by taking into consideration
again the function $f\left( x_{1},x_{2},x_{3},x_{4}\right) $ (\ref%
{factorised_func}) (which is assumed now to represent the joint pdf of four
continuous random variables) and showing how, thanks to these rules, the
marginal pdf $f\left( x_{3}\right) $ can be evaluated in a step-by-step
fashion. If the messages $\vec{m}_{1}\left( x_{1}\right) =f_{1}\left(
x_{1}\right) $, $\vec{m}_{0}\left( x_{2}\right) =1$ and $\overset{\leftarrow 
}{m}_{6}\left( x_{4}\right) =1$ are defined, applying Eqs. (\ref{CR_1})--(%
\ref{CR_2}) to the factor graph shown in Fig. \ref{Fig_1} leads to the
ordered computation of the messages 
\begin{equation}
\,\vec{m}_{2}\left( x_{1}\right) =\int f_{2}\left( x_{1},x_{2}\right) \,\vec{%
m}_{0}\left( x_{2}\right) \,dx_{2},  \label{marginal_5}
\end{equation}%
\begin{equation}
\,\vec{m}_{3}\left( x_{1}\right) =\vec{m}_{1}\left( x_{1}\right) \,\vec{m}%
_{2}\left( x_{1}\right) ,  \label{marginal_4}
\end{equation}%
\begin{equation}
\vec{m}_{4}\left( x_{3}\right) =\int f_{3}\left( x_{1},x_{3}\right) \,\vec{m}%
_{3}\left( x_{1}\right) \,dx_{1},\,  \label{marginal_3}
\end{equation}%
and%
\begin{equation}
\reflectbox{\ensuremath{\vec{\reflectbox{\ensuremath{m}}}}}_{5}\left(
x_{3}\right) =\int f_{4}\left( x_{3},x_{4}\right) \,\reflectbox{\ensuremath{%
\vec{\reflectbox{\ensuremath{m}}}}}_{6}\left( x_{4}\right) \,dx_{4}.
\label{marginal_2}
\end{equation}%
Then, given the messages $\vec{m}_{4}(x_{3})$ (\ref{marginal_3}) and $%
\reflectbox{\ensuremath{\vec{\reflectbox{\ensuremath{m}}}}}_{5}\left(
x_{3}\right) $ (\ref{marginal_2}), referring to the same edge, but
originating from opposite directions, the required marginal is evaluated as 
\begin{equation}
f\left( x_{3}\right) =\vec{m}_{4}\left( x_{3}\right) \,\,\reflectbox{%
\ensuremath{\vec{\reflectbox{\ensuremath{m}}}}}_{5}\left( x_{3}\right) .
\label{marginal_1}
\end{equation}%
This result is \emph{exact} since the graph representing the joint pdf $%
f\left( x_{1},x_{2},x_{3},x_{4}\right) $ (\ref{factorised_func}) is \emph{%
cycle free}, i.e. it does not contain closed paths. When the considered
graph does not have this property, the SPA can still be employed (e.g., see 
\cite[Par. III.A]{Loeliger_2007} and \cite[Sec. V]{Kschischang_2001}), but
its application leads to \emph{iterative} message passing algorithms, that,
in general, produce approximate results. Moreover, the order according to
which messages are passed on a cycle (i.e., the \emph{message scheduling})
has to be properly selected. Despite this, it is widely accepted that the
most important applications of the SPA refer to cyclic graphs \cite%
{Kschischang_2001}.

The last important issue related to the application of the SPA\ is the
availability of closed form expressions for the passed messages when, like
in the filtering problem investigated in this manuscript, the involved
variables are continuous. In the following, the pdfs of all the considered
random vectors are Gaussian or are approximated through a set of $N_{p}$
weighted particles. In the first case, the pdf of a random vector $\mathbf{x}
$ is conveyed by the message%
\begin{equation}
\vec{m}_{G}\left( \mathbf{x}\right) =\mathcal{N(}\mathbf{x};\mathbf{\eta },%
\mathbf{C}),  \label{Gaussian_message}
\end{equation}%
where $\mathbf{\eta }$ and $\mathbf{C}$ denote the mean and the covariance
of $\mathbf{x}$, respectively. In the second case, instead, its pdf is
conveyed by the message%
\begin{equation}
\vec{m}_{P}\left( \mathbf{x}\right) =\sum_{j=1}^{N_{p}}\vec{m}_{P,j}\left( 
\mathbf{x}\right) ,  \label{particle message}
\end{equation}%
where%
\begin{equation}
\vec{m}_{P,j}\left( \mathbf{x}\right) \triangleq w_{j}\,\delta \left( 
\mathbf{x}-\mathbf{x}_{j}\right)  \label{particle message_j}
\end{equation}%
represents the $j-$th \emph{component} of the message $\vec{m}_{P}(\mathbf{x}%
)$ (\ref{particle message}), i.e. the contribution of the $j-$th particle $%
\mathbf{x}_{j}$ and its weight $w_{j}$ to such a message. Luckily, various
closed form results are available for these two types of messages; the few
mathematical rules required in the computation of all the messages appearing
in our filtering algorithms can be found in Tables I--III of \cite[App. A,
p. 1534]{Vitetta_2019}.

\subsection{Graphical Modelling for a Network of Bayesian Filters and
Message Passing on it}

In this manuscript, we consider a discrete-time SSM whose $D-$dimensional 
\emph{hidden state} in the $k-$th interval is denoted $\mathbf{x}%
_{k}\triangleq \lbrack x_{0,k},x_{1,k},...,$ $x_{D-1,k}]^{T}$, and whose 
\emph{state update} and \emph{measurement models} are expressed by%
\begin{equation}
\mathbf{x}_{k+1}=\mathbf{f}_{k}\left( \mathbf{x}_{k}\right) +\mathbf{w}_{k}
\label{eq:X_update}
\end{equation}%
and%
\begin{eqnarray}
\mathbf{y}_{k} &\triangleq& [y_{0,k},y_{1,k},...,y_{P-1,k}]^{T}  \notag \\
&=&\mathbf{h}_{k}\left( \mathbf{x}_{k}\right) +\mathbf{e}_{k},
\label{meas_mod}
\end{eqnarray}%
respectively. Here, $\mathbf{f}_{k}\left( \mathbf{x}_{k}\right) $ ($\mathbf{h%
}_{k}\left( \mathbf{x}_{k}\right) $) is a time-varying $D-$dimensional ($P-$%
dimensional) real function and $\mathbf{w}_{k}$ ($\mathbf{e}_{k}$) the $k-$%
th element of the process (measurement) noise sequence $\left\{ \mathbf{w}%
_{k}\right\} $ ($\left\{ \mathbf{e}_{k}\right\} $); this sequence consists
of $D-$dimensional ($P-$dimensional) \emph{independent and identically
distributed} (iid) Gaussian noise vectors, each characterized by a zero mean
and a covariance matrix $\mathbf{C}_{w}$ ($\mathbf{C}_{e}$). Moreover,
statistical independence between $\left\{ \mathbf{e}_{k}\right\} $ and $\{%
\mathbf{w}_{k}\}$ is assumed for simplicity. Note that, from a statistical
viewpoint, the SSM described by Eqs. (\ref{eq:X_update})--(\ref{meas_mod})
is characterized by the \emph{Markov model} $f(\mathbf{x}_{k+1}|\mathbf{x}%
_{k})$ and the \emph{observation model} $f(\mathbf{y}_{k}|\mathbf{x}_{k})$
for any $k$.

\begin{figure}[tbp]
\centering
\includegraphics[width=0.4\textwidth]{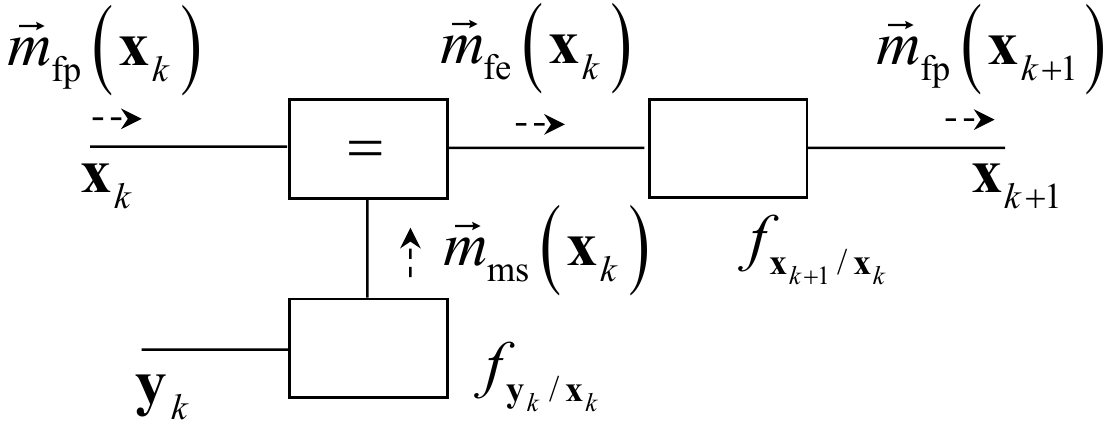}
\caption{Message passing over the factor graph representing the $k-$th
recursion of Bayesian filtering. A SSM characterized by the Markov model $f(%
\mathbf{x}_{k+1}|\mathbf{x}_{k})$ and the observation model $f(\mathbf{y}%
_{k}|\mathbf{x}_{k})$ is considered.}
\label{Fig_3}
\end{figure}

In the following sections, we focus on the so-called \emph{filtering problem}%
, which concerns the evaluation of the posterior pdf $f(\mathbf{x}_{k}|%
\mathbf{y}_{1:t})$ at an instant $t\geq 1$, given a) the initial pdf $f(%
\mathbf{x}_{1})$ and b) the $t\cdot P$-dimensional measurement vector $%
\mathbf{y}_{1:t}=\left[ \mathbf{y}_{1}^{T},\mathbf{y}_{2}^{T},...,\mathbf{y}%
_{t}^{T}\right] ^{T}$. It is well known that, if the pdf $f(\mathbf{x}_{1})$
referring to the first observation interval is known, the computation of the
posterior (i.e., filtered) pdf $f(\mathbf{x}_{t}|\mathbf{y}_{1:t})$ for $%
t\geq 1$ can be accomplished by means of an exact \emph{Bayesian recursive
procedure}, consisting of a \emph{measurement update} step followed by a 
\emph{time update} step. In \cite[Sec. III]{Vitetta_2019}, it is shown that,
if this procedure is formulated with reference to the joint pdf $f(\mathbf{x}%
_{t},\mathbf{y}_{1:t})$ (in place of the associated a posteriori pdf $f(%
\mathbf{x}_{t}|\mathbf{y}_{1:t})$), its $k-$th recursion (with $k=1,2,...,t$%
) can be represented as a \emph{forward only} message passing algorithm over
the cycle free factor graph shown in Fig. \ref{Fig_3}. In the measurement
update, the message $\vec{m}_{\mathrm{fe}}(\mathbf{x}_{k})$ going out of the
equality node is computed as\footnote{%
In the following the acronyms \textrm{fp}, \textrm{fe}, \textrm{ms} and 
\textrm{pm} are employed in the subscripts of various messages, so that
readers can easily understand their meaning; in fact, the messages these
acronyms refer to convey a \emph{forward prediction}, a \emph{forward
estimate}, \emph{measurement }information and \emph{pseudo-measurement}
information, respectively.} (see Eq. (\ref{CR_1})) 
\begin{eqnarray}
\vec{m}_{\mathrm{fe}}\left( \mathbf{x}_{k}\right) &=&\vec{m}_{\mathrm{fp}%
}\left( \mathbf{x}_{k}\right) \,\vec{m}_{\mathrm{ms}}\left( \mathbf{x}%
_{k}\right)  \notag \\
&=&f(\mathbf{x}_{k},\mathbf{y}_{1:k}),  \label{eq_forward_est}
\end{eqnarray}%
where%
\begin{equation}
\vec{m}_{\mathrm{fp}}(\mathbf{x}_{k})\triangleq f(\mathbf{x}_{k},\mathbf{y}%
_{1:k-1})  \label{m_fp_k}
\end{equation}%
is the message feeding the considered graph. Note that the messages $\vec{m}%
_{\mathrm{fp}}(\mathbf{x}_{k})$ (\ref{m_fp_k}) and $\vec{m}_{\mathrm{fe}}(%
\mathbf{x}_{k})$ convey the predicted pdf (i.e., the \emph{forward prediction%
}) of $\mathbf{x}_{k}$ computed in the previous (i.e., in the $(k-1)-$th)
recursion and the filtered pdf (i.e., the \emph{forward estimate}) of $%
\mathbf{x}_{k}$ computed in the considered recursion, respectively, whereas
the message $\vec{m}_{\mathrm{ms}}(\mathbf{x}_{k})\triangleq f\left( \mathbf{%
y}_{k}\left\vert \mathbf{x}_{k}\right. \right) $ conveys the statistical
information provided by the measurement $\mathbf{y}_{k}$ (\ref{meas_mod}).

In the time update, the message that emerges from the function node
referring to the pdf $f(\mathbf{x}_{k+1}|\mathbf{x}_{k})$ is evaluated as
(see Eq. (\ref{CR_2})) 
\begin{equation}
\int f\left( \mathbf{x}_{k+1}\left\vert \mathbf{x}_{k}\right. \right) \,\vec{%
m}_{\mathrm{fe}}\left( \mathbf{x}_{k}\right) d\mathbf{x}_{k}=f(\mathbf{x}%
_{k+1},\mathbf{y}_{1:k});  \label{eq:message_fp_new}
\end{equation}%
such a message is equal to $\vec{m}_{\mathrm{fp}}(\mathbf{x}_{k+1})$ (see
Eq. (\ref{m_fp_k}))

Let us take into consideration now a \emph{network of }$N_{F}$ \emph{%
interconnected\ Bayesian filters}. In the following, we assume that:

a) All the filters of the network are fed by the same measurement vector
(namely, $\mathbf{y}_{k}$ (\ref{meas_mod})), work in parallel and cooperate
in order to estimate the state vector $\mathbf{x}_{k}$; in doing so, they
can fully share their statistical information.

b) The $i-$th filter of the network (with $i=1$, $2$, $...$, $N_{F}$),
denoted F$_{i}$, works on a lower dimensional space and, in particular,
estimates the portion $\mathbf{x}_{k}^{(i)}$ (having size $D_{i}$, with $%
D_{i}\leq D$) of the state vector $\mathbf{x}_{k}$; therefore, the substate $%
\mathbf{\bar{x}}_{k}^{(i)}$, representing the portion of $\mathbf{x}_{k}$
not included in $\mathbf{x}_{k}^{(i)}$, can be considered as a \emph{nuisance%
} vector for F$_{i}$.

c) The set $\{\mathbf{x}_{k}^{(i)}\}$, collecting the substates estimated by
all the filters of the network, covers $\mathbf{x}_{k}$, but does not
necessarily represent a partition of it. In other words, unlike MPF, some
overlapping between the substates estimated by different filters is
admitted. This means that the filtering algorithm running on the whole
network may contain a form of \emph{redundancy}, since one or more elements
of the state vector can be independently estimated by different Bayesian
filters.

We are interested in developing recursive filtering algorithms for the whole
network of Bayesian filters. The approach we propose to solve this problem
consists of the following three steps: S1) building a factor graph that
allows us to represent the measurement and time updates accomplished by each
filter of the network and its interactions with the other filters as message
passing algorithms on it; S2) developing a graphical model for the whole
network on the basis of the factor graphs devised in the first step; S3)
deriving new filtering methods as message passing algorithms over the whole
graphical model obtained in the second step.

Let us focus, now, on step S1. In developing a graphical model for filter F$%
_{i}$, the following considerations must be kept into account:

1) Since the portion $\mathbf{\bar{x}}_{k}^{(i)}$ of $\mathbf{x}_{k}$ is
unknown to F$_{i}$ (and, consequently, represents a \emph{nuisance state}),
an \emph{estimate} of its pdf $f_{k}(\mathbf{\bar{x}}_{k}^{(i)})$ must be
provided by the other filters of the network; this allows F$_{i}$ to
integrate out the dependence of its Markov model $f(\mathbf{x}_{k+1}^{(i)}|%
\mathbf{x}_{k}^{(i)},\mathbf{\bar{x}}_{k}^{(i)})$ and of its observation
model $f(\mathbf{y}_{k}|\mathbf{x}_{k}^{(i)},\mathbf{\bar{x}}_{k}^{(i)})$ on 
$\mathbf{\bar{x}}_{k}^{(i)}$.

2) Filter F$_{i}$ can benefit from the \emph{pseudo-measurements} computed
on the basis of the statistical information provided by the other filters of
the network.

As far as the last point is concerned, it is worth pointing out that, in
this manuscript, any pseudo-measurement represents a \emph{fictitious}
measurement computed on the basis of the statistical information provided by
a filtering algorithm different from the one benefiting from it; despite
this, it can be processed as if it was a \emph{real} measurement, provided
that its statistical model is known. In practice, a pseudo-measurement $%
\mathbf{z}_{k}^{(i)}$ made available to the filter F$_{i}$ is a $P_{i}-$%
dimensional random vector that, similarly as the real measurement $\mathbf{y}%
_{k}$ (\ref{meas_mod}), can be modelled as\footnote{%
The possible dependence of the pseudo-measurement $\mathbf{z}_{k}^{(i)}$ (%
\ref{pseudo_meas_model}) on the substate $\mathbf{\bar{x}}_{k}^{(i)}$ is
ignored here, for simplicity.}%
\begin{equation}
\mathbf{z}_{k}^{(i)}=\mathbf{\tilde{h}}_{k}\left( \mathbf{x}%
_{k}^{(i)}\right) +\mathbf{\tilde{e}}_{k}^{(i)},  \label{pseudo_meas_model}
\end{equation}%
where $\mathbf{\tilde{h}}_{k}\left( \mathbf{x}_{k}\right) $ is a
time-varying $P_{i}-$dimensional function and $\mathbf{\tilde{e}}_{k}^{(i)}$
is a zero mean $P_{i}-$dimensional noise vector. The evaluation of these
fictitious measurements is often based on the \emph{mathematical constraints}
established by the Markov model of the considered SSM, as shown in the
following section, where a specific network of filters is considered.

Based on the considerations illustrated above, the equations describing the
measurement/time updates accomplished by F$_{i}$ in the $k-$th recursion of
the network can be formulated as follows. At the beginning of this
recursion, F$_{i}$ is fed by the forward prediction%
\begin{equation}
\vec{m}_{\mathrm{fp}}\left( \mathbf{x}_{k}^{(i)}\right) =f(\mathbf{x}%
_{k}^{(i)},\mathbf{y}_{1:k-1}),  \label{eq_fp_Fi}
\end{equation}%
originating from the previous recursion. In its first step (i.e., in its
measurement update), it computes two filtered pdfs (i.e., two forward
estimates), the first one based on the measurement $\mathbf{y}_{k}$ (\ref%
{meas_mod}), the second one on the pseudo-measurement $\mathbf{z}_{k}^{(i)}$
(\ref{pseudo_meas_model}). The first filtered pdf is evaluated as (see Eq. (%
\ref{eq_forward_est}))%
\begin{equation}
\vec{m}_{\mathrm{fe}1}\left( \mathbf{x}_{k}^{(i)}\right) =\vec{m}_{\mathrm{fp%
}}\left( \mathbf{x}_{k}^{(i)}\right) \vec{m}_{\mathrm{ms}}\left( \mathbf{x}%
_{k}^{(i)}\right) ,  \label{eq_fe1_Fi}
\end{equation}%
where 
\begin{equation}
m_{\mathrm{ms}}\left( \mathbf{x}_{k}^{(i)}\right) \triangleq \int f\left( 
\mathbf{y}_{k}\left\vert \mathbf{x}_{k}^{(i)},\mathbf{\bar{x}}%
_{k}^{(i)}\right. \right) \,m_{\mathrm{mg}1}\left( \mathbf{\bar{x}}%
_{k}^{(i)}\right) \,d\mathbf{\bar{x}}_{k}^{(i)}  \label{eq_ms_Fi}
\end{equation}%
and $m_{\mathrm{mg}1}(\mathbf{\bar{x}}_{k}^{(i)})$ are the messages
conveying \emph{measurement} information and a filtered (or predicted) pdf
of $\mathbf{\bar{x}}_{k}^{(i)}$ provided by the other filters, respectively.
Similarly, the second filtered pdf is evaluated as (see Eq. (\ref%
{eq_forward_est})) 
\begin{equation}
\vec{m}_{\mathrm{fe}2}\left( \mathbf{x}_{k}^{(i)}\right) =\vec{m}_{\mathrm{fe%
}1}\left( \mathbf{x}_{k}^{(i)}\right) \vec{m}_{\mathrm{pm}}\left( \mathbf{x}%
_{k}^{(i)}\right) ,  \label{eq_fe2_Fi}
\end{equation}%
where\footnote{%
If the pseudo-measurement $\mathbf{z}_{k}^{(i)}$ (\ref{pseudo_meas_model})
depends also on $\mathbf{\bar{x}}_{k}^{(i)}$, marginalization with respect
to this substate is required in the computation of the following message.}%
\begin{equation}
m_{\mathrm{pm}}\left( \mathbf{x}_{k}^{(i)}\right) \triangleq f\left( \mathbf{%
z}_{k}^{(i)}\left\vert \mathbf{x}_{k}^{(i)}\right. \right)  \label{eq_pm_Fi}
\end{equation}%
is the message conveying \emph{pseudo-measurement} information. Then, in its
second step (i.e., in its time update), F$_{i}$ computes the new forward
prediction (see Eq. (\ref{eq:message_fp_new})) 
\begin{equation}
\vec{m}_{\mathrm{fp}}\left( \mathbf{x}_{k+1}^{(i)}\right) =\int \int f\left( 
\mathbf{x}_{k+1}^{(i)}\left\vert \mathbf{x}_{k}^{(i)},\mathbf{\bar{x}}%
_{k}^{(i)}\right. \right) \vec{m}_{\mathrm{fe}2}\left( \mathbf{x}%
_{k}^{(i)}\right) \cdot m_{\mathrm{mg}2}\left( \mathbf{\bar{x}}%
_{k}^{(i)}\right) d\mathbf{x}_{k}\,d\mathbf{\bar{x}}_{k}^{(i)},
\label{eq_fpnew_Fi}
\end{equation}%
where $m_{\mathrm{mg}2}(\mathbf{\bar{x}}_{k}^{(i)})$ has the same meaning as 
$m_{\mathrm{mg}1}(\mathbf{\bar{x}}_{k}^{(i)})$ (see Eq. (\ref{eq_ms_Fi})),
but is not necessarily equal to it (since more refined information about $%
\mathbf{\bar{x}}_{k}^{(i)}$ could be made available by the other filters of
the network after that the message $m_{\mathrm{ms}}(\mathbf{x}_{k}^{(i)})$ (%
\ref{eq_ms_Fi}) has been computed).

\begin{figure}[tbp]
\centering
\includegraphics[width=0.6\textwidth]{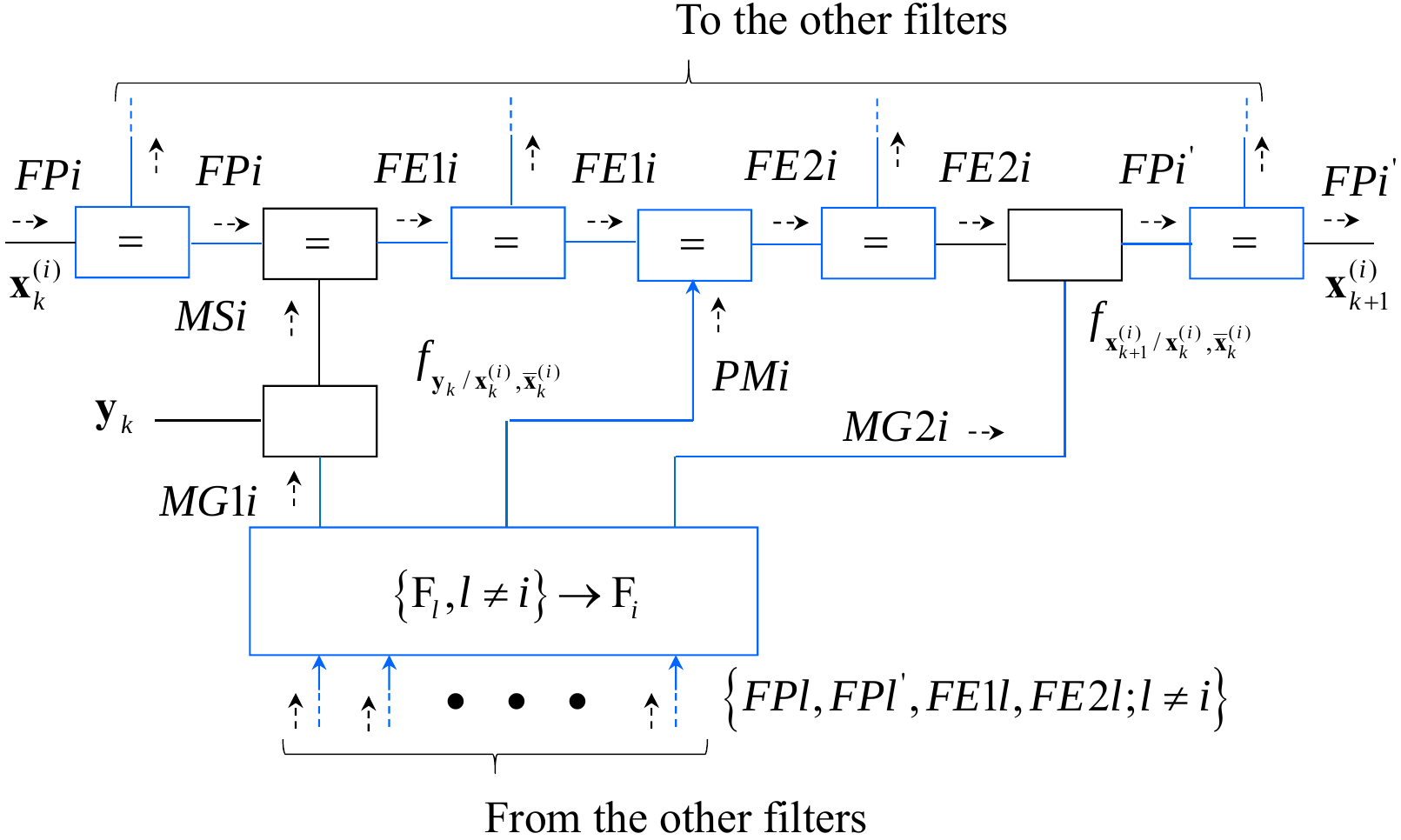}
\caption{Message passing over the factor graph representing the couple of
measurement updates and the time update accomplished by the $i-$th Bayesian
filter in the $k-$th recursion of the network it belongs to. The messages $%
\vec{m}_{\mathrm{fp}}(\mathbf{x}_{k}^{(i)})$, $\vec{m}_{\mathrm{fp}}(\mathbf{%
\ x}_{k+1}^{(i)})$, $\vec{m}_{\mathrm{ms}}(\mathbf{x}_{k}^{(i)})$, $\vec{m}%
_{ \mathrm{mg}l}(\mathbf{\bar{x}}_{k}^{(i)})$, $\vec{m}_{\mathrm{pm}}(%
\mathbf{x} _{k}^{(i)})$ and $\vec{m}_{\mathrm{fe}l}(\mathbf{x}_{k}^{(i)})$
are denoted $FPi$, $FPi^{^{\prime }}$, $MSi$, $MGli$, $PMi$\ and $FEli$,
respectively, to ease reading.}
\label{Fig_4}
\end{figure}

Formulas (\ref{eq_fe1_Fi})-(\ref{eq_fe2_Fi}) and (\ref{eq_fpnew_Fi}) involve
only products of pdfs and integrations of products; for this reason, their
evaluation can be represented as a \emph{forward only} message passing over
the \emph{cycle free} factor graph shown in Fig. \ref{Fig_4}. Note that, if
this graph is compared with the one shown in Fig. \ref{Fig_3}, the following
additional elements (identified by blue lines) are found:

1) Five equality nodes - Four of them allow to generate copies of the
messages $\vec{m}_{\mathrm{fp}}(\mathbf{x}_{k}^{(i)})$, $\vec{m}_{\mathrm{fe}%
1}(\mathbf{x}_{k}^{(i)})$, $\vec{m}_{\mathrm{fe}2}(\mathbf{x}_{k}^{(i)})$
and $\vec{m}_{\mathrm{fp}}(\mathbf{x}_{k+1}^{(i)})$, to be shared with the
other filters of the network, whereas the remaining one is involved in the
second measurement update of F$_{i}$.

2) A block in which the predicted/filtered pdfs $\{\vec{m}_{\mathrm{fp}}(%
\mathbf{x}_{k}^{(l)}),\vec{m}_{\mathrm{fp}}(\mathbf{x}_{k+1}^{(l)})$, $\vec{m%
}_{\mathrm{fe}q}(\mathbf{x}_{k}^{(l)})$; $q=1$, $2$ and $l\neq i\}$ provided
by the other filters of the network are processed - In this block, the
messages $m_{\mathrm{mg}q}(\mathbf{\bar{x}}_{k}^{(i)})$ (with $q=1$ and $2$)
and $m_{\mathrm{pm}}(\mathbf{x}_{k}^{(i)})$ are computed (see Eqs. (\ref%
{eq_ms_Fi}), (\ref{eq_pm_Fi}) and (\ref{eq_fpnew_Fi})); this block is
connected to \emph{oriented} edges only, i.e. to edges on which the flow of
messages is unidirectional.

\begin{figure}[tbp]
\centering
\includegraphics[width=0.6\textwidth]{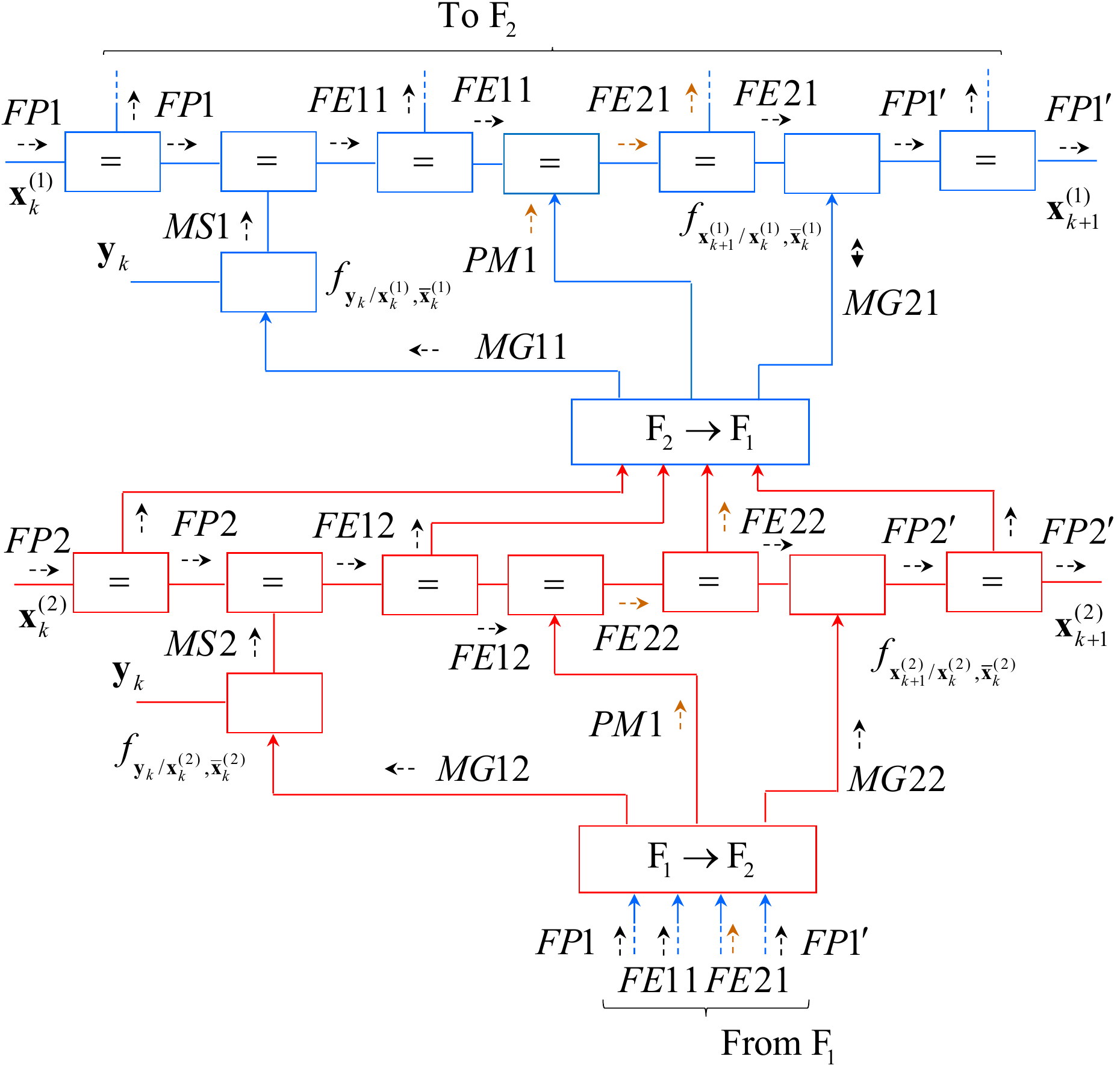}
\caption{Graphical model based on the factor graph shown in Fig. \protect\ref%
{Fig_4} and referring to the interconnection of two Bayesian filters; the
presence of a closed path (cycle) on which messages can be passed multiple
times is highlighted by brown arrows.}
\label{Fig_5}
\end{figure}

Given the graphical model represented in Fig. \ref{Fig_4}, step S2 can be
accomplished by adopting the same conceptual approach as \cite[Sec. III]%
{Vitetta_2019}, where the factor graph on which RBPF and dual RBPF are based
is devised by merging two sub-graphs, that refer to distinct substates. For
this reason, a graphical model for the whole network of $N_{F}$ Bayesian
filters can be developed by interconnecting $N_{F}$ distinct factor graphs,
each structured like the one shown in that Figure. For instance, if $N_{F}=2$
is assumed for simplicity, this procedure results in the graphical model
shown in Fig. \ref{Fig_5}. It is important to note that, in this case, if
the substates $\mathbf{x}_{k}^{(1)}$ and $\mathbf{x}_{k}^{(2)}$ estimated by
F$_{1}$ and F$_{2}$, respectively, do not form a partition of state vector $%
\mathbf{x}_{k}$, they share a portion of it; this consists of $%
N_{d}\triangleq D_{1}+D_{2}-D$ state variables, that are separately
estimated by the two Bayesian filters. The parameter $N_{d}$ can be
considered as the \emph{degree of redundancy} characterizing the considered
network of filters. The presence of redundancy in a filtering algorithm may
result in an improvement of estimation accuracy and/or tracking capability;
however, this is obtained at the price of an increased complexity with
respect to the case in which F$_{1}$ and F$_{2}$ are run on disjoint
substates.

Once the graphical model for the whole network has been developed, step S3
can be easily accomplished. In fact, recursive filtering algorithms for the
considered network can be derived by systematically applying the SPA to its
graphical model after that a proper \emph{scheduling} has been established
for the exchange of messages among its $N_{F}$ Bayesian filters. Moreover,
in developing a specific filtering algorithm to be run on a network of
Bayesian filters, we must always keep in mind that:

1) Its $k-$th recursion is fed by the set of forward predictions $\{\vec{m}_{%
\mathrm{fp}}(\mathbf{x}_{k}^{(i)})$, $i=1$, $2$, $...$, $N_{F}\}$, and
generates $N_{F}$ couples of filtered densities $\{(\vec{m}_{\mathrm{fe}1}(%
\mathbf{x}_{k}^{(i)}),\vec{m}_{\mathrm{fe}2}(\mathbf{x}_{k}^{(i)}))$, $i=1$, 
$2$, $...$, $N_{F}\}$ and $N_{F}$ new forward predictions $\{\vec{m}_{%
\mathrm{fp}}(\mathbf{x}_{k+1}^{(i)})$, $i=1$, $2$, $...$, $N_{F}\}$.
Moreover, similarly as MPF, a \emph{joint} filtered density for the whole
state $\mathbf{x}_{k}$ is \emph{unavailable} (unless the substate of one or
more of the employed Bayesian filters coincides with $\mathbf{x}_{k}$) and
multiple filtered/predicted pdfs are available for any substate shared by
distinct filters.

2) Specific algorithms are needed to compute the pseudo-measurement and the
nuisance substate pdfs in the $\{$F$_{l}$, $l\neq i\}{\rightarrow} $F$_{i}$
block appearing in Fig. \ref{Fig_5}. These algorithms depend on the
considered SSM and on the selected message scheduling; for this reason, a
general description of their structure cannot be provided.

3) The graphical model shown in Fig. \ref{Fig_5}, unlike the one illustrated
in Fig. \ref{Fig_3}, is \emph{not cycle free}; the presence of cycles is
highlighted in the considered figure by showing the flow of messages along
one of them. The presence of cycles raises the problems of a) identifying
all the messages that can be iteratively refined and b) establishing the
order according to which they are computed. Generally speaking, iterative
message passing on the graphical model referring to a network of filters
involves both the couple of measurement updates and the time update
accomplished by all the interconnected filters. In fact, this should allow
each Bayesian filter to a) progressively refine the nuisance substate
density employed in its measurement/time updates, and b) improve the quality
of the pseudo-measurements exploited in its second measurement update. For
this reason, if $n_{i}$ iterations are run, the overall computational
complexity of each recursion is multiplied by $n_{i}$.

In the following section, a specific application of the general principles
illustrated in this paragraph is analysed.

\section{Filtering Algorithms Based on the Interconnection of an Extended
Kalman Filter with a Particle Filter\label{sec:Message-Passing}}

In this section we focus on the development of two new filtering algorithms
based on the interconnection of an extended Kalman filter with a particle
filter. We first describe the graphical models on which these algorithms are
based. Then, we provide a detailed description of the computed messages and
their scheduling in a specific case. Finally, we provide a detailed analysis
of the computational complexity of the devised algorithms.

\subsection{Graphical Modelling}

In this section, we develop new filtering algorithms for the class of \emph{%
conditionally linear Gaussian} SSMs \cite{Schon_2005}, \cite%
{Vitetta_2018_bis}, \cite{Vitetta_2019}; this allows us to partition the
state vector in the $k-$th interval as $\mathbf{x}_{k}=[(\mathbf{x}%
_{k}^{(L)})^{T},(\mathbf{x}_{k}^{(N)})^{T}]^{T}$, where $\mathbf{x}%
_{k}^{(L)}\triangleq \lbrack x_{0,k}^{(L)}$, $%
x_{1,k}^{(L)},...,x_{D_{L}-1,k}^{(L)}]^{T}$ ($\mathbf{x}_{k}^{(N)}\triangleq
\lbrack x_{0,k}^{(N)},x_{1,k}^{(N)},...,x_{D_{N}-1,k}^{(N)}]^{T}$) is its 
\emph{linear }(\emph{nonlinear}) \emph{component} (with $D_{N}+D_{L}=D$).
The devised algorithms rely on the following assumptions:

1) They involve two interconnected Bayesian filters, denoted F$_{1}$ and F$%
_{2}$.

2) Filter F$_{2}$ is a \emph{particle filter}\footnote{%
In particular, a \emph{sequential importance resampling} filter is employed 
\cite{Arulampalam_2002}.} and estimates the nonlinear state component only
(so that $\mathbf{x}_{k}^{(2)}=\mathbf{x}_{k}^{(N)}$ and $\mathbf{\bar{x}}%
_{k}^{(2)}=\mathbf{x}_{k}^{(L)}$).

3) Filter F$_{1}$ is an\emph{\ extended Kalman filter} and works on the 
\emph{whole system state} or on the \emph{linear state component} only.
Consequently, in the first case (denoted \textbf{C.1} in the following), $%
\mathbf{x}_{k}^{(1)}=\mathbf{x}_{k}$ and $\mathbf{\bar{x}}_{k}^{(1)}$ is
empty, and both the interconnected filters estimate the nonlinear state
component (for this reason, the corresponding degree of redundancy is $%
N_{d}=D_{N}$). In the second case (denoted \textbf{C.2} in the following),
instead, $\mathbf{x}_{k}^{(1)}=\mathbf{x}_{k}^{(L)}$ and $\mathbf{\bar{x}}%
_{k}^{(1)}=\mathbf{x}_{k}^{(N)}$, and the two filters estimate disjoint
substates (consequently, $N_{d}=0$).

This network configuration has been mainly inspired by RBPF. In fact,
similarly as RBPF, the filtering techniques we develop are based on the idea
of concatenating a local filtering method (EKF) with a global method (PF).
However, unlike RBPF, a \emph{single\ }extended Kalman filter is employed in
place of a bank of Kalman filters. It is also worth remembering that, on the
one hand, the use of a particle filter interconnected with an extended
Kalman filter for tracking \emph{disjoint} substates has been suggested in 
\cite[Par. 3.2]{Djuric_2007}, where, however, no filtering algorithm based
on this idea has been derived. On the other hand, a filtering scheme based
on the interconnection of the same filters, but working on \emph{partially
overlapped substates}, has been derived in \cite{Montorsi_2013}, where it
has also been successfully applied to inertial navigation.

\begin{figure}[tbp]
\centering
\includegraphics[width=0.6\textwidth]{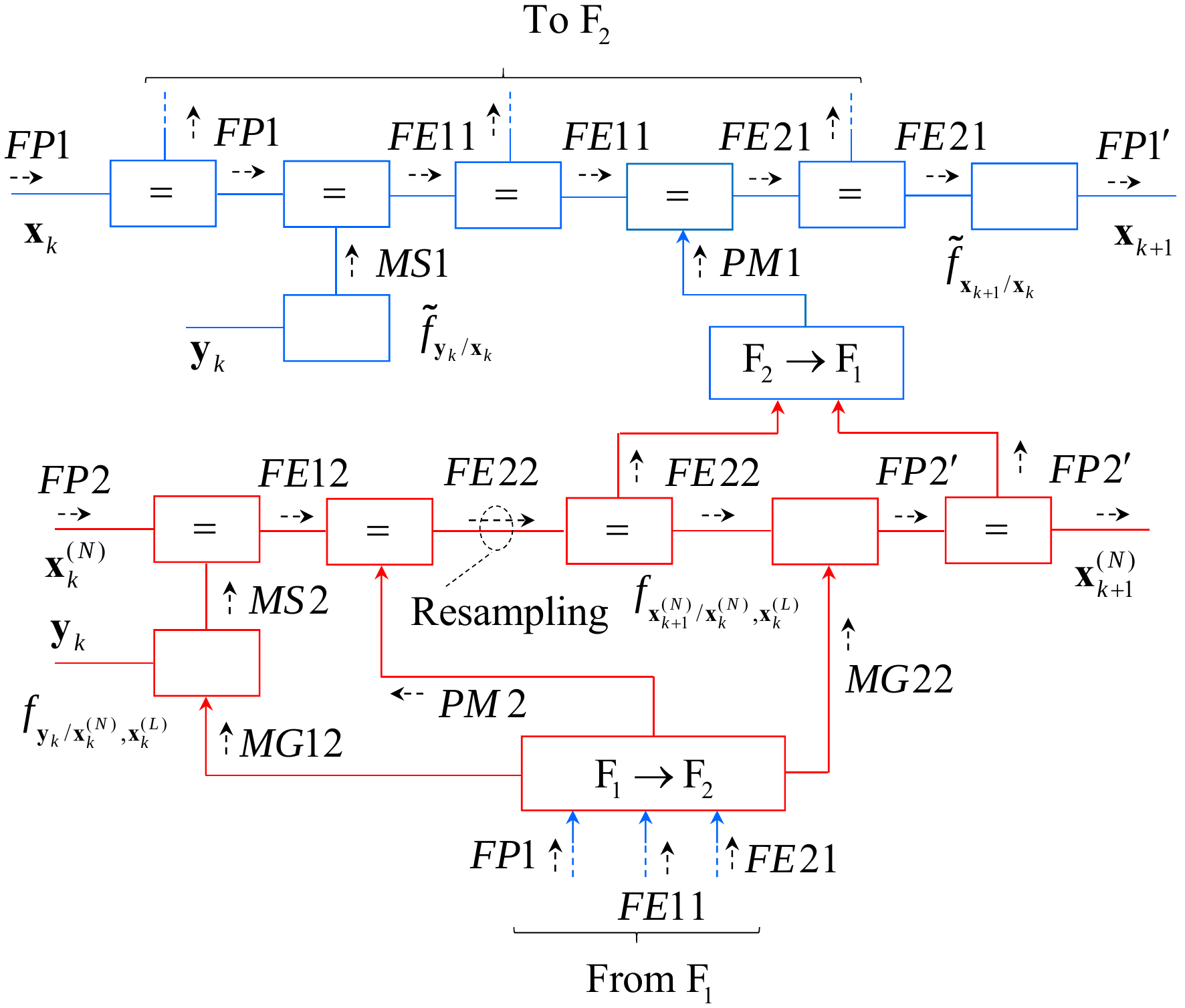}
\caption{Graphical model based on the factor graph shown in Fig. \protect\ref%
{Fig_5} and referring to the interconnection of an extended Kalman filter (F$%
_{1}$) with a particle filter (F$_{2}$).}
\label{Fig_6}
\end{figure}

Based on the graphical model shown in Fig. \ref{Fig_5}, the factor graph
illustrated in Fig. \ref{Fig_6} can be drawn for\ case \textbf{C.1}. It is
important to point out that:

1) Filter F$_{1}$ is based on \emph{linearised} (and, consequently, \emph{%
approximate}) Markov/measurement models of the considered SSM, whereas
filter F$_{2}$ relies on \emph{exact} models, as explained in more detail
below.

2) Since the nuisance substate $\mathbf{\bar{x}}_{k}^{(1)}$ is empty, no
marginalization is required in F$_{1}$; for this reason, the messages $\{%
\vec{m}_{\mathrm{mg}q}(\mathbf{\bar{x}}_{k}^{(1)})$; $q=1,2\}$ (i.e., $MG11$
and $MG21$) visible in Fig. \ref{Fig_5} do not appear in Fig. \ref{Fig_6}.

3) The new predicted pdf $\vec{m}_{\mathrm{fp}}(\mathbf{x}_{k+1}^{(2)})=\vec{%
m}_{\mathrm{fp}}(\mathbf{x}_{k+1}^{(N)})$ and the second filtered pdf $\vec{m%
}_{\mathrm{fe}2}(\mathbf{x}_{k}^{(2)})=\vec{m}_{\mathrm{fe}2}(\mathbf{x}%
_{k}^{(N)})$ computed by F$_{2}$ (i.e., the messages $FP2^{^{\prime }}$ and $%
FE22$, respectively) feed the F$_{2}{\rightarrow }$F$_{1}$ block, where they
are jointly processed to generate the pseudo-measurement message $\vec{m}_{%
\mathrm{pm}}(\mathbf{x}_{k}^{(1)})=\vec{m}_{\mathrm{pm}}(\mathbf{x}_{k})$ ($%
PM1$) feeding F$_{1}$. Similarly, as shown below, the computation of the
pseudo-measurement message exploited by F$_{2}$ (i.e., of the message $\vec{m%
}_{\mathrm{pm}}(\mathbf{x}_{k}^{(2)})=\vec{m}_{\mathrm{pm}}(\mathbf{x}%
_{k}^{(N)})$, $PM2$) requires the knowledge of a new predicted pdf that
refers, however, to the \emph{linear state component only}. In our graphical
model, the computation of this prediction is accomplished by the F$_{1}{%
\rightarrow }$F$_{2}$ block; this explains why the new predicted pdf $\vec{m}%
_{\mathrm{fp}}(\mathbf{x}_{k+1}^{(1)})=\vec{m}_{\mathrm{fp}}(\mathbf{x}%
_{k+1})$ ($FP1^{^{\prime }}$) evaluated by F$_{1}$ and referring to the 
\emph{whole state} of the considered SSM, does not feed the F$_{1}{%
\rightarrow }$F$_{2}$ block.

4) Particle \emph{resampling} with replacement has been included in the
portion of the graphical model referring to filter F$_{2}$. This important
task, accomplished after the second measurement update of this filter, does
not emerge from the application of the SPA to our graphical model and
ensures that the particles emerging from it are all equally likely. Note
also that, because of the presence of particle resampling, two versions of
the second filtered pdf $\vec{m}_{\mathrm{fe}2}(\mathbf{x}_{k}^{(2)})=\vec{m}%
_{\mathrm{fe}2}(\mathbf{x}_{k}^{(N)})$ ($FE22$) become available, one before
resampling, the other one after it. As shown in the next paragraph, the
second version of this message is exploited in the computation of the
pseudo-measurement message $\vec{m}_{\mathrm{pm}}(\mathbf{x}_{k}^{(1)})=\vec{%
m}_{\mathrm{pm}}(\mathbf{x}_{k})$ ($PM1$).

In the remaining part of this paragraph, we first provide various details
about the filters F$_{1}$ and F$_{2}$, and the way pseudo-measurements are
computed for each of them; then, we comment on how the factor graph shown in
Fig. \ref{Fig_6} should be modified if case \textbf{C.2} is considered.

\emph{Filter} F$_{1}$ - Filter F$_{1}$ is based on the \emph{linearized}
versions of Eqs. (\ref{eq:X_update}) and (\ref{meas_mod}), i.e. on the
models (e.g., see \cite[pp. 194-195]{Anderson_1979})%
\begin{equation}
\mathbf{x}_{k+1}=\mathbf{F}_{k}\,\mathbf{x}_{k}+\mathbf{u}_{k}+\mathbf{w}_{k}
\label{state_up_approx}
\end{equation}%
and 
\begin{equation}
\mathbf{y}_{k}=\mathbf{H}_{k}^{T}\,\mathbf{x}_{k}+\mathbf{v}_{k}+\mathbf{e}%
_{k},  \label{meas_mod_approx}
\end{equation}%
respectively; here, $\mathbf{F}_{k}\triangleq \lbrack \partial \mathbf{f}%
_{k}\left( \mathbf{x}\right) /\partial \mathbf{x}]_{\mathbf{x=x}_{\mathrm{fe}%
,k}}$, $\mathbf{u}_{k}\triangleq \mathbf{f}_{k}\left( \mathbf{x}_{\mathrm{fe}%
,k}\right) -\mathbf{F}_{k}\,\mathbf{x}_{\mathrm{fe},k}$, $\mathbf{H}%
_{k}^{T}\triangleq \lbrack \partial \mathbf{h}_{k}\left( \mathbf{x}\right)
/\partial \mathbf{x}]_{\mathbf{x=x}_{\mathrm{fp},k}}$, $\mathbf{v}%
_{k}\triangleq \mathbf{h}_{k}\left( \mathbf{x}_{\mathrm{fp},k}\right) -%
\mathbf{H}_{k}^{T}\mathbf{x}_{\mathrm{fp},k}$ and $\mathbf{x}_{\mathrm{fp}%
,k} $ ($\mathbf{x}_{\mathrm{fe},k}$) is the \emph{forward prediction }(\emph{%
forward estimate}) of $\mathbf{x}_{k}$ computed by F$_{1}$ in its $(k-1)-$th
($k-$th) recursion. Consequently, the approximate models%
\begin{equation}
\tilde{f}\left( \mathbf{x}_{k+1}\left\vert \mathbf{x}_{k}\right. \right) =%
\mathcal{N}\left( \mathbf{x}_{k};\mathbf{F}_{k}\,\mathbf{x}_{k}+\mathbf{u}%
_{k},\mathbf{C}_{w}\right)  \label{Markov_ekf_approx}
\end{equation}%
and 
\begin{equation}
\tilde{f}\left( \mathbf{y}_{k}\left\vert \mathbf{x}_{k}\right. \right) =%
\mathcal{N}\left( \mathbf{x}_{k};\mathbf{H}_{k}^{T}\,\mathbf{x}_{k}+\mathbf{v%
}_{k},\mathbf{C}_{e}\right)  \label{Measurement_ekf_approx}
\end{equation}%
appear in the graphical model shown in Fig. \ref{Fig_6}.

\emph{Filter} F$_{2}$ - In developing filter F$_{2}$, we assume that the
portion of Eq. (\ref{eq:X_update}) referring to the nonlinear state
component (i.e., the last $D_{N}$ lines of the considered Markov model) and
that the observation model (\ref{meas_mod}) can be put in the form (e.g.,
see \cite[eqs. (3)-(4)]{Vitetta_2019}) 
\begin{equation}
\mathbf{x}_{k+1}^{(N)}=\mathbf{A}_{k}^{(N)}\left( \mathbf{x}%
_{k}^{(N)}\right) \mathbf{x}_{k}^{(L)}+\mathbf{f}_{k}^{(N)}\left( \mathbf{x}%
_{k}^{(N)}\right) +\mathbf{w}_{k}^{(N)}  \label{eq:XN_update}
\end{equation}%
and%
\begin{equation}
\mathbf{y}_{k}=\mathbf{g}_{k}\left( \mathbf{x}_{k}^{(N)}\right) +\mathbf{B}%
_{k}\left( \mathbf{x}_{k}^{(N)}\right) \mathbf{x}_{k}^{(L)}+\mathbf{e}_{k},
\label{eq:y_t}
\end{equation}%
respectively. In Eq. (\ref{eq:XN_update}), $\mathbf{f}_{k}^{(N)}(\mathbf{x}%
_{k}^{(N)})$ ($\mathbf{A}_{k}^{(N)}(\mathbf{x}_{k}^{(N)})$) is a
time-varying $D_{N}-$dimensional real function ($D_{N}\times D_{L}$ real
matrix) and $\mathbf{w}_{k}^{(N)}$ consists of the last $D_{N}$ elements of
the noise term $\mathbf{w}_{k}$ appearing in Eq. (\ref{eq:X_update}) (the
covariance matrix of $\mathbf{w}_{k}^{(N)}$ is denoted $\mathbf{C}_{w}^{(N)}$%
); moreover, in Eq. (\ref{eq:y_t}), $\mathbf{g}_{k}(\mathbf{x}_{k}^{(N)})$ ($%
\mathbf{B}_{k}(\mathbf{x}_{k}^{(N)})$) is a time-varying $P-$dimensional
real function ($P\times D_{L}$ real matrix). This explains why filter F$_{2}$
is based on the exact pdfs 
\begin{equation}
f\left( \mathbf{x}_{k+1}^{(N)}\left\vert \mathbf{x}_{k}^{(N)},\mathbf{x}%
_{k}^{(L)}\right. \right) =\mathcal{N}\left( \mathbf{x}_{k}^{(N)};\mathbf{A}%
_{k}^{(N)}\left( \mathbf{x}_{k}^{(N)}\right) \mathbf{x}_{k}^{(L)}+\mathbf{f}%
_{k}^{(N)}\left( \mathbf{x}_{k}^{(N)}\right) ,\mathbf{C}_{w}^{(N)}\right)
\label{f_x_N}
\end{equation}%
and 
\begin{equation}
f\left( \mathbf{y}_{k}\left\vert \mathbf{x}_{k}^{(N)},\mathbf{x}%
_{k}^{(L)}\right. \right) = \mathcal{N}\left( \mathbf{x}_{k};\mathbf{g}%
_{k}\left( \mathbf{x}_{k}^{(N)}\right) +\mathbf{B}_{k}\left( \mathbf{x}%
_{k}^{(N)}\right) \mathbf{x}_{k}^{(L)},\mathbf{C}_{e}\right) ,  \label{f_y_N}
\end{equation}%
that appear in the graphical model shown in Fig. \ref{Fig_6}.

\emph{Computation of the pseudo-measurements for filter} F$_{1}$ - Filter F$%
_{1}$ is fed by pseudo-measurement information about the \emph{whole state} $%
\mathbf{x}_{k}$, i.e. about both the substates $\mathbf{x}_{k}^{(L)}$ and $%
\mathbf{x}_{k}^{(N)}$. On the one hand, $N_{p}$ pseudo-measurements about
the nonlinear state component are provided by the $N_{p}$ particles
contributing to the filtered pdf $\vec{m}_{\mathrm{fe}2}(\mathbf{x}%
_{k}^{(N)})$ ($FE22$) available \emph{after particle resampling}. On the
other hand, $N_{p}$ pseudo-measurements about the linear state component are
evaluated by means of the same method employed by RBPF for this task. This
method is based on the idea that the random vector (see \cite[Par. II.D, p.
2283, eq. (24a)]{Schon_2005} and \cite[Sec. III, p. 1524, eq. (9)]%
{Vitetta_2019})%
\begin{equation}
\mathbf{z}_{k}^{(L)}\triangleq \mathbf{x}_{k+1}^{(N)}-\mathbf{f}%
_{k}^{(N)}\left( \mathbf{x}_{k}^{(N)}\right) ,  \label{eq:z_L_l}
\end{equation}%
depending on the \emph{nonlinear state component} \emph{only}, must equal
the sum (see Eq. (\ref{eq:XN_update}))%
\begin{equation}
\mathbf{A}_{k}^{(N)}\left( \mathbf{x}_{k}^{(N)}\right) \mathbf{x}_{k}^{(L)}+%
\mathbf{w}_{k}^{(N)},  \label{eq:z_L_l_bis}
\end{equation}%
that depends on the \emph{linear state component}. For this reason, $N_{p}$
realizations of $\mathbf{z}_{k}^{(L)}$ (\ref{eq:z_L_l}) are computed in the F%
$_{2}{\rightarrow }$F$_{1}$ block on the basis of the messages $\vec{m}_{%
\mathrm{fe}2}(\mathbf{x}_{k}^{(N)})$ ($FE22$) and $\vec{m}_{\mathrm{fp}}(%
\mathbf{x}_{k+1}^{(N)})$ ($FP2^{^{\prime }}$) and are treated as
measurements about $\mathbf{x}_{k}^{(L)}$.

\emph{Computation of the pseudo-measurements for filter} F$_{2}$ - The
messages feeding F$_{1}{\rightarrow }$F$_{2}$ block are employed for: a)
generating a pdf of $\mathbf{x}_{k}^{(L)}$, so that the dependence of the
state update and measurement models (i.e., of the densities $f(\mathbf{x}%
_{k+1}^{(N)}|\mathbf{x}_{k}^{(N)}$, $\mathbf{x}_{k}^{(L)})$ (\ref{f_x_N})
and $f(\mathbf{y}_{k}|\mathbf{x}_{k}^{(N)},\mathbf{x}_{k}^{(L)})$ (\ref%
{f_y_N}), respectively) on this substate can be integrated out; b) computing
pseudo-measurement information about $\mathbf{x}_{k}^{(N)}$. As far as the
last point is concerned, the approach we adopt is the same as that developed
for \emph{dual} RBPF in \cite[Sec. V, pp. 1528-1529]{Vitetta_2019}. Such an
approach relies on the Markov model 
\begin{equation}
\mathbf{x}_{k+1}^{(L)}=\mathbf{A}_{k}^{(L)}\left( \mathbf{x}%
_{k}^{(N)}\right) \mathbf{x}_{k}^{(L)}+\mathbf{f}_{k}^{(L)}\left( \mathbf{x}%
_{k}^{(N)}\right) +\mathbf{w}_{k}^{(L)},  \label{eq:XL1_update}
\end{equation}%
referring to the \emph{linear} state component \cite{Vitetta_2018_bis}, \cite%
{Vitetta_2019}; in the last expression, $\mathbf{f}_{k}^{(L)}(\mathbf{x}%
_{k}^{(N)})$ ($\mathbf{A}_{k}^{(L)}(\mathbf{x}_{k}^{(N)})$) is a
time-varying $D_{L}-$dimensional real function ($D_{L}\times D_{L}$ real
matrix), and $\mathbf{w}_{k}^{(N)}$ consists of the first $D_{L}$ elements
of the noise term $\mathbf{w}_{k}$ appearing in Eq. (\ref{eq:X_update}) (the
covariance matrix of $\mathbf{w}_{k}^{(L)}$ is denoted $\mathbf{C}_{w}^{(L)}$%
, and independence between $\{\mathbf{w}_{k}^{(L)}\}$ and $\{\mathbf{w}%
_{k}^{(N)}\}$ is assumed for simplicity). From Eq. (\ref{eq:XL1_update}) it
is easily inferred that the random vector 
\begin{equation}
\mathbf{z}_{k}^{(N)}\triangleq \mathbf{x}_{k+1}^{(L)}-\mathbf{A}%
_{k}^{(L)}\left( \mathbf{x}_{k}^{(N)}\right) \,\mathbf{x}_{k}^{(L)}\text{,}
\label{z_N_l}
\end{equation}%
equals the sum 
\begin{equation}
\mathbf{f}_{k}^{(L)}\left( \mathbf{x}_{k}^{(N)}\right) +\mathbf{w}_{k}^{(L)},
\label{z_N_l_bis}
\end{equation}%
that depends on $\mathbf{x}_{k}^{(N)}$ \emph{only}; for this reason, $%
\mathbf{z}_{k}^{(N)}$ (\ref{z_N_l}) can be interpreted as a
pseudo-measurement about $\mathbf{x}_{k}^{(N)}$. In this case, the
generation of pseudo-measurement information can be summarised as follows.
First, $N_{p}$ pdfs, one for each of the particles conveyed by the message $%
\vec{m}_{\mathrm{fe}2}(\mathbf{x}_{k}^{(N)})$ ($FE22$), are computed for the
random vector $\mathbf{z}_{k}^{(N)}$ (\ref{z_N_l}) by exploiting the
statistical information about the linear state component made available by F$%
_{1}$. Then, each of these pdfs is \emph{correlated} with the pdf obtained
for $\mathbf{z}_{k}^{(N)}$ under the assumption that this vector is
expressed by Eq. (\ref{z_N_l_bis}); this procedure results in a set of $%
N_{p} $ particle weights, different from those computed on the basis of $%
\mathbf{y}_{k}$ (\ref{eq:y_t}) in the first measurement update of F$_{2}$.

A graphical model similar to the one shown in Fig. \ref{Fig_6} can be easily
derived from the general model appearing in Fig. \ref{Fig_5} for\ case 
\textbf{C.2} too. The relevant differences with respect to case \textbf{C.1}
can be summarized as follows:

1) Filters F$_{1}$ and F$_{2}$ estimate $\mathbf{x}_{k}^{(1)}=\mathbf{x}%
_{k}^{(L)}$ and $\mathbf{x}_{k}^{(2)}=\mathbf{x}_{k}^{(N)}$, respectively;
consequently, their nuisance substates are $\mathbf{\bar{x}}_{k}^{(1)}=%
\mathbf{x}_{k}^{(N)}$ and $\mathbf{\bar{x}}_{k}^{(2)}=\mathbf{x}_{k}^{(L)}$,
respectively.

2) The F$_{2}{\rightarrow }$F$_{1}$ block is fed by the predicted/filtered
pdfs computed by F$_{2}$; such pdfs are employed for: a) for providing F$%
_{1} $ with a pdf for $\mathbf{x}_{k}^{(N)}$, so that dependence of the
Markov model (see Eq. (\ref{eq:XL1_update}))%
\begin{equation}
f\left( \mathbf{x}_{k+1}^{(L)}\left\vert \mathbf{x}_{k}^{(N)},\mathbf{x}%
_{k}^{(L)}\right. \right) =\mathcal{N}\left( \mathbf{x}_{k}^{(L)};\mathbf{A}%
_{k}^{(L)}\left( \mathbf{x}_{k}^{(N)}\right) \mathbf{x}_{k}^{(L)}+\mathbf{f}%
_{k}^{(L)}\left( \mathbf{x}_{k}^{(N)}\right) ,\mathbf{C}_{w}^{(L)}\right)
\label{f_x_L}
\end{equation}%
and of the measurement model $f(\mathbf{y}_{k}|\mathbf{x}_{k}^{(N)},\mathbf{x%
}_{k}^{(L)})$ (\ref{f_y_N}) on this substate can be integrated out;\ b)
generating pseudo-measurement information about the substate $\mathbf{x}%
_{k}^{(L)}$ \emph{only}. As far as point a) is concerned, it is also
important to point out that the approximate model $\tilde{f}(\mathbf{y}_{k}|%
\mathbf{x}_{k}^{(L)})$ ($\tilde{f}(\mathbf{x}_{k+1}^{(L)}|\mathbf{x}%
_{k}^{(L)})$) on which F$_{1}$ is based can be derived from Eq. (\ref{f_y_N}%
) (Eq. (\ref{f_x_L})) after setting $\mathbf{x}_{k}^{(N)}=\mathbf{x}_{%
\mathrm{fp},k}^{(N)}$ ($\mathbf{x}_{k}^{(N)}=\mathbf{x}_{\mathrm{fe}%
,k}^{(N)} $); here, $\mathbf{x}_{\mathrm{fp},k}^{(N)}$ ($\mathbf{x}_{\mathrm{%
fe},k}^{(N)}$) denote the prediction (the estimate) of $\mathbf{x}_{k}^{(N)}$
evaluated on the basis of the message $\vec{m}_{\mathrm{fp}}(\mathbf{x}%
_{k}^{(N)})$ ($\vec{m}_{\mathrm{fe}2}(\mathbf{x}_{k}^{(N)})$) computed by F$%
_{2}$. Moreover, since Eqs. (\ref{eq:y_t}) and (\ref{eq:XL1_update}) exhibit
a \emph{linear} dependence on $\mathbf{x}_{k}^{(L)}$, F$_{1}$ becomes a
standard Kalman filter.

The derivation of a specific filtering algorithm based on the graphical
models described in this paragraph requires defining the scheduling of the
messages passed on them and deriving mathematical expressions for such
messages. These issues are investigated in detail in the following paragraph.

\subsection{Message Scheduling and Computation}

\begin{figure}[tbp]
\centering
\includegraphics[width=0.6\textwidth]{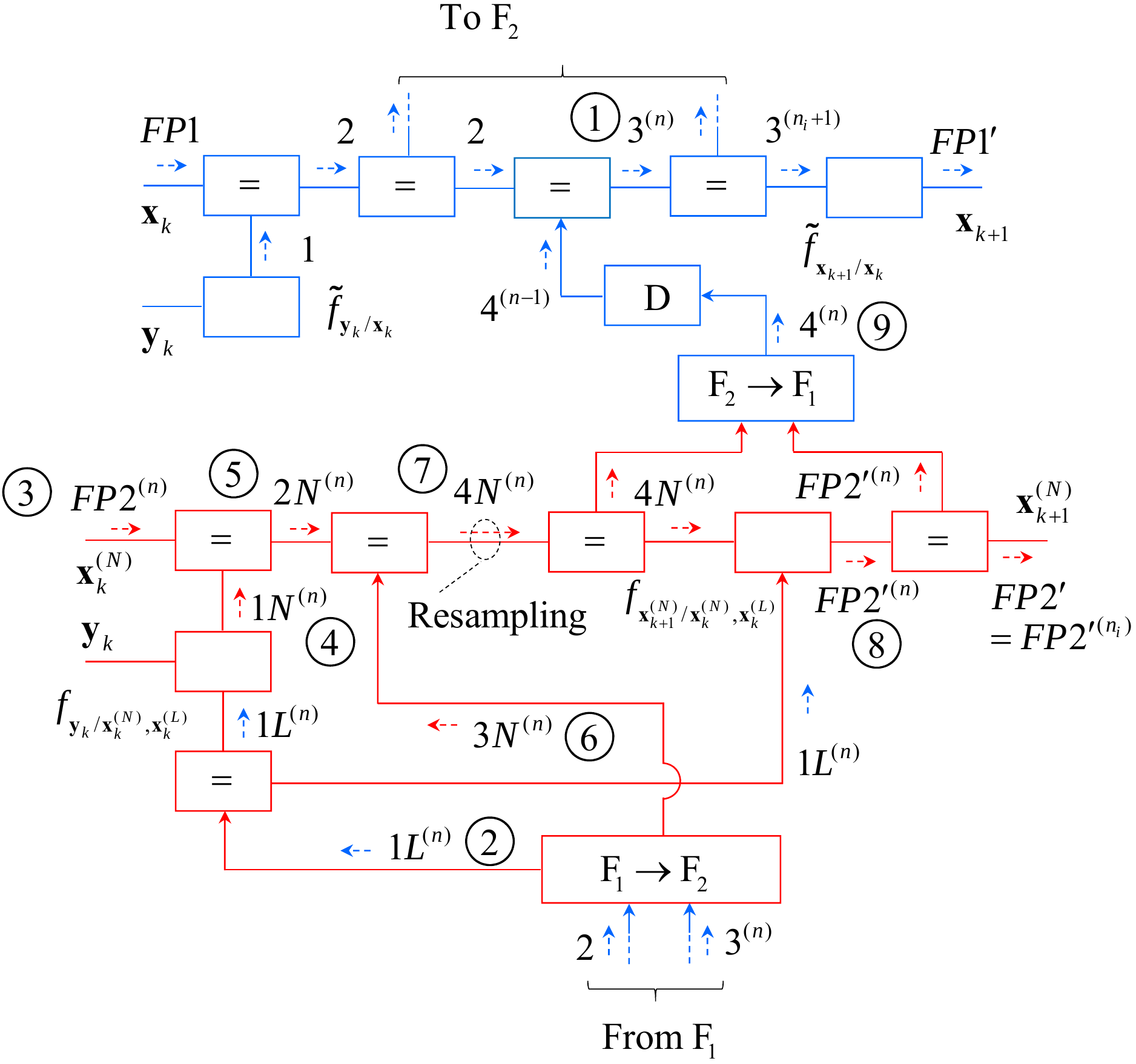}
\caption{Message scheduling adopted in the $k$-th recursion of the DBF
technique. The circled integers $1-9$ specify the order according to which
nine distinct messages are computed in the $n$-th iteration of phase II.}
\label{Fig_7}
\end{figure}

In this paragraph, a \emph{recursive filtering technique}, called \emph{dual
Bayesian filtering }(DBF) and based on the graphical model illustrated in
Fig. \ref{Fig_6}, is developed. In each recursion of the DBF technique, F$%
_{1}$ is run before F$_{2}$; moreover, the presence of cycles in the graph
on which it is based is accounted for by including a procedure for the
iterative computation of the messages passed on them. Our description of the
selected scheduling relies on Fig. \ref{Fig_7}, that refers to the $k-$th
recursion and to the $n-$th iteration accomplished within this recursion
(with $n=1$, $2$, $...$, $n_{i}$, where $n_{i}$ represents the overall
number of iterations). It is important to point out that the following
changes have been made in Fig. \ref{Fig_7} with respect to Fig. \ref{Fig_6}:

1) A simpler notation has been adopted for the messages to ease reading. In
particular, the symbols $FP2^{(n)}$, $FP2^{^{\prime }(n)}$, $q$ ($q^{(n)}$), 
$qL$ ($qL^{(n)}$) and $qN$ ($qN^{(n)}$) represent the messages $\vec{m}_{%
\mathrm{fp}}^{(n)}(\mathbf{x}_{k}^{(N)})$, $\vec{m}_{\mathrm{fp}}^{(n)}(%
\mathbf{x}_{k+1}^{(N)})$, $\vec{m}_{q}(\mathbf{x}_{k})$ ($\vec{m}_{q}^{(n)}(%
\mathbf{x}_{k})$), $\vec{m}_{q}(\mathbf{x}_{k}^{(L)})$ ($\vec{m}_{q}^{(n)}(%
\mathbf{x}_{k}^{(L)})$) and $\vec{m}_{q}(\mathbf{x}_{k}^{(N)})$ ($\vec{m}%
_{q}^{(n)}(\mathbf{x}_{k}^{(N)})$), respectively; moreover, the integer
parameter $n$ appearing in the superscript of some of them represents the
iteration index.

2) Blue (red) arrows have been employed to identify Gaussian messages
(messages in other forms).

3) The F$_{1}{\rightarrow }$F$_{2}$ block is fed by the two filtered pdfs of 
$\mathbf{x}_{k}$ computed by F$_{1}$ (i.e., by the messages $\vec{m}_{2}(%
\mathbf{x}_{k})$ and $\vec{m}_{3}^{(n)}(\mathbf{x}_{k})$), but not by the
predicted pdf $\vec{m}_{\mathrm{fp}}(\mathbf{x}_{k})$, since the last
message is useless.

3) The forward prediction $\vec{m}_{\mathrm{fp}}^{(n)}(\mathbf{x}_{k}^{(N)})$
feeding F$_{2}$ is involved in the proposed iterative procedure and may
change from iteration to iteration because of resampling (in fact, this may
lead to discarding a portion of the particles conveyed by this message); for
this reason, its dependence on the iteration index $n$ has been explicitly
indicated.

4) The same message (namely, $\vec{m}_{1}^{(n)}(\mathbf{x}_{k}^{(L)})$) is
employed in F$_{2}$ for integrating out the dependence of the Markov model $%
f(\mathbf{x}_{k+1}^{(N)}|\mathbf{x}_{k}^{(N)},\mathbf{x}_{k}^{(L)})$ (\ref%
{f_x_N}) and of the measurement model $f(\mathbf{y}_{k}|\mathbf{x}_{k}^{(N)},%
\mathbf{x}_{k}^{(L)})$ (\ref{f_y_N}) on the linear component $\mathbf{x}%
_{k}^{(L)}$.

5) A memory cell (identified by the label `D') has been added to store the
last message evaluated in each iteration (i.e., the pseudo-measurement
message $m_{4}^{(n)}(\mathbf{x}_{k})$), so that it can be made available to F%
$_{1}$ at the beginning of the next iteration.

The DBF technique, at the beginning of its $k-$th recursion, is fed by the
message 
\begin{equation}
\vec{m}_{\mathrm{fp}}\left( \mathbf{x}_{k}\right) =\mathcal{N}\left( \mathbf{%
x}_{k};\mathbf{\eta }_{\mathrm{fp},k},\mathbf{C}_{\mathrm{fp},k}\right)
\label{eq:eq:message_fp_l_EKF}
\end{equation}%
and 
\begin{equation}
\vec{m}_{\mathrm{fp}}\left( \mathbf{x}_{k}^{(N)}\right) =\sum_{j=1}^{N_{p}}%
\vec{m}_{\mathrm{fp},j}\left( \mathbf{x}_{k}^{(N)}\right) ,
\label{eq:eq:message_fp_PF}
\end{equation}%
that corresponds to $FP2^{(1)}$ in Fig. \ref{Fig_7}; here,%
\begin{equation}
\vec{m}_{\mathrm{fp},j}\left( \mathbf{x}_{k}^{(N)}\right) =w_{p}\,\mathbf{%
\delta }\left( \mathbf{x}_{k}^{(N)}-\mathbf{x}_{k,j}^{(N)}\right)
\label{eq:message_N_pred-1}
\end{equation}%
is the $j-$th component of $\vec{m}_{\mathrm{fp}}(\mathbf{x}_{k}^{(N)})$, $%
\mathbf{x}_{k,j}^{(N)}$ is the $j-$th particle \emph{predicted} in the
previous (i.e., in the $(k-1)-$th) recursion and $w_{p}\triangleq 1/N_{p}$
is its weight. The DBF processes the messages $\vec{m}_{\mathrm{fp}}(\mathbf{%
x}_{k})$ (\ref{eq:eq:message_fp_l_EKF}) and $\vec{m}_{\mathrm{fp}}(\mathbf{x}%
_{k}^{(N)})$ (\ref{eq:eq:message_fp_PF}), and the new measurement $\mathbf{y}%
_{k}$ (\ref{eq:y_t}), and generates: a) a couple of filtered densities for
both $\mathbf{x}_{k}$ and $\mathbf{x}_{k}^{(N)}$; b) the output messages $%
\vec{m}_{\mathrm{fp}}(\mathbf{x}_{k+1})$ and $\vec{m}_{\mathrm{fp}}(\mathbf{x%
}_{k+1}^{(N)})$, having the \emph{same functional form} as $\vec{m}_{\mathrm{%
fp}}(\mathbf{x}_{k})$ (\ref{eq:eq:message_fp_l_EKF}) and $\vec{m}_{\mathrm{fp%
}}(\mathbf{x}_{k}^{(N)})$ (\ref{eq:eq:message_fp_PF}), respectively. The
message passing accomplished to achieve these results can be divided in the
three consecutive phases listed below.

\textbf{I} - In the first phase, filter F$_{1}$ accomplishes its first
measurement update on the basis of the forward prediction $\vec{m}_{\mathrm{%
fp}}\left( \mathbf{x}_{k}\right) $ and of the new measurement $\mathbf{y}%
_{k} $. This leads to the ordered computation of the messages $\vec{m}_{1}(%
\mathbf{x}_{k})$ and $\vec{m}_{2}(\mathbf{x}_{k})$.

\textbf{II} - In the second phase, an iterative procedure involving the
first measurement update and the time update of F$_{2}$, and the computation
of pseudo-measurements and their exploitation in the second measurement
update of each filter is carried out. The $n-$th iteration of this procedure
can be divided into six consecutive steps and leads to the ordered
computation of the following messages: 1) $\vec{m}_{3}^{(n)}\left( \mathbf{x}%
_{k}\right) $, $\vec{m}_{1}^{(n)}(\mathbf{x}_{k}^{(L)})$; 2) $\vec{m}_{%
\mathrm{fp}}^{(n)}(\mathbf{x}_{k}^{(N)})$, $\vec{m}_{1}^{(n)}(\mathbf{x}%
_{k}^{(N)})$, $\vec{m}_{2}^{(n)}(\mathbf{x}_{k}^{(N)})$; 3) $\vec{m}%
_{3}^{(n)}(\mathbf{x}_{k}^{(N)})$; 4) $\vec{m}_{4}^{(n)}(\mathbf{x}%
_{k}^{(N)})$; 5) $\vec{m}_{\mathrm{fp}}^{(n)}(\mathbf{x}_{k+1}^{(N)})$; 6) $%
\vec{m}_{4}^{(n)}(\mathbf{x}_{k})$.

\textbf{III} - In the third phase, the new predictions $\vec{m}_{\mathrm{fp}%
}\left( \mathbf{x}_{k+1}\right) $ and $\vec{m}_{\mathrm{fp}}(\mathbf{x}%
_{k+1}^{(N)})$ are generated by F$_{1}$ and F$_{2}$, respectively. This
involves the ordered computation of the following messages: $\vec{m}_{%
\mathrm{fp}}(\mathbf{x}_{k+1}^{(N)})$, $\vec{m}_{3}^{(n_{i}+1)}\left( 
\mathbf{x}_{k}\right) $ and $\vec{m}_{\mathrm{fp}}(\mathbf{x}_{k+1})$.

In the remaining part of this paragraph, the expressions of all the messages
computed in each of the three phases described above are provided; the
derivation of these expressions is sketched in Appendix \ref{app:A}.

\textbf{Phase I} - In this phase, the forward prediction $\vec{m}_{\mathrm{fp%
}}(\mathbf{x}_{k})$ (\ref{eq:eq:message_fp_l_EKF}) feeding filter F$_{1}$\
is merged with the message%
\begin{equation}
\vec{m}_{1}\left( \mathbf{x}_{k}\right) =\mathcal{N}\left( \mathbf{x}_{k};%
\mathbf{\eta }_{1,k},\mathbf{C}_{1,k}\right) ,  \label{m_ms_EKF_bis}
\end{equation}%
conveying measurement information; the covariance matrix $\mathbf{C}_{1,k}$
and the mean vector $\mathbf{\eta }_{1,k}$ of the last message are evaluated
on the basis of the associated precision matrix 
\begin{equation}
\mathbf{W}_{1,k}\triangleq \left( \mathbf{C}_{1,k}\right) ^{-1}=\mathbf{H}%
_{k}\mathbf{W}_{e}\mathbf{H}_{k}^{T},  \label{w_ms_EKF}
\end{equation}%
and of the associated transformed mean vector 
\begin{equation}
\mathbf{w}_{1,k}\triangleq \mathbf{W}_{1,k}\,\mathbf{\eta }_{1,k}=\mathbf{H}%
_{k}\mathbf{W}_{e}\left( \mathbf{y}_{k}-\mathbf{v}_{k}\right) ,
\label{W_ms_EKF}
\end{equation}%
respectively, with $\mathbf{W}_{e}\triangleq \mathbf{C}_{e}^{-1}$. This
results in the \emph{first filtered pdf} (see Fig. \ref{Fig_7})%
\begin{eqnarray}
\vec{m}_{2}\left( \mathbf{x}_{k}\right) &=&\vec{m}_{\mathrm{fp}}(\mathbf{x}%
_{k})\,\vec{m}_{1}\left( \mathbf{x}_{k}\right)  \label{m_fe1_EKF_bis} \\
&=&\mathcal{N}\left( \mathbf{x}_{k};\mathbf{\eta }_{2,k},\mathbf{C}%
_{2,k}\right)  \label{m_fe1_EKF_tris}
\end{eqnarray}%
computed by filter F$_{1}$; here, the covariance matrix $\mathbf{C}_{2,k}$
and the mean vector $\mathbf{\eta }_{2,k}$ are evaluated on the basis of the
associated precision matrix 
\begin{equation}
\mathbf{W}_{2,k}\triangleq \left( \mathbf{C}_{2,k}\right) ^{-1}=\mathbf{W}_{%
\mathrm{fp},k}+\mathbf{W}_{1,k},  \label{W_fe1_EKFa}
\end{equation}%
and of the associated transformed mean vector 
\begin{equation}
\mathbf{w}_{2,k}\triangleq \mathbf{W}_{2,k}\,\mathbf{\eta }_{2,k}=\mathbf{w}%
_{\mathrm{fp},k}+\mathbf{w}_{1,k},  \label{w_fe1_EKF}
\end{equation}%
respectively, $\mathbf{W}_{\mathrm{fp},k}\triangleq (\mathbf{C}_{\mathrm{fp}%
,k})^{-1}$ and $\mathbf{w}_{\mathrm{fp},k}\triangleq \mathbf{W}_{\mathrm{fp}%
,k}\,\mathbf{\eta }_{\mathrm{fp},k}$.

\textbf{Phase II} - A short description of the six steps accomplished in the 
$n-$th iteration of this phase is provided in the following. As shown below,
the elements of the particle set processed by F$_{2}$ can change from
iteration to iteration, even if its cardinality remains the same. In the
following, the particle set available\emph{\ at the beginning} of the $n-$th
iteration is denoted $S_{k}[n]=\{\mathbf{x}_{k,j}^{(N)}[n]$; $j=1$, $2$, $%
... $, $N_{p}\}$; note that the initial particle set is $S_{k}[1]\triangleq
\{\mathbf{x}_{k,j}^{(N)},\,j=1$, $2$, $...$, $N_{p}\}$ (i.e., $\mathbf{x}%
_{k,j}^{(N)}[1]=\mathbf{x}_{k,j}^{(N)}$ for any $j$) and collects the $N_{p}$
\emph{predicted particles} conveyed by the message $\vec{m}_{\mathrm{fp}}(%
\mathbf{x}_{k}^{(N)})$ (\ref{eq:eq:message_fp_PF}).

1) \emph{Second measurement update\ in }F$_{1}$ - The \emph{second filtered
pdf} (see Fig. \ref{Fig_7})%
\begin{eqnarray}
\vec{m}_{3}^{(n)}\left( \mathbf{x}_{k}\right) &=&\vec{m}_{2}\left( \mathbf{x}%
_{k}\right) \,\vec{m}_{4}^{(n-1)}\left( \mathbf{x}_{k}\right) \,
\label{m_fe2_EKF} \\
&=&\mathcal{N}\left( \mathbf{x}_{k};\mathbf{\eta }_{3,k}^{(n)},\mathbf{C}%
_{3,k}^{(n)}\right)  \label{m_fe2_EKF_bis}
\end{eqnarray}%
is computed by F$_{1}$ in order to exploit the pseudo-measurement message $%
\vec{m}_{4}^{(n-1)}(\mathbf{x}_{k})$ (evaluated in the previous iteration);
since $\vec{m}_{4}^{(n-1)}(\mathbf{x}_{k})=1$ for $n=1$ (note that F$_{1}$
cannot benefit from pseudo-measurement information at the beginning of the
first iteration) and $\vec{m}_{4}^{(n-1)}(\mathbf{x}_{k})=\mathcal{N}(%
\mathbf{x}_{k};\mathbf{\eta }_{4,k}^{(n-1)},\mathbf{C}_{4,k}^{(n-1)})$ for $%
n>1$ (see Eq. (\ref{message_pm_x_l_k})), it easy to show that 
\begin{equation}
\mathbf{C}_{3,k}^{(n)}=\mathbf{C}_{2,k}  \label{C_3_1}
\end{equation}%
and%
\begin{equation}
\mathbf{\eta }_{3,k}^{(n)}=\mathbf{\eta }_{2,k}  \label{eta_3_1}
\end{equation}%
for $n=1$, whereas 
\begin{equation}
\mathbf{C}_{3,k}^{(n)}=\mathbf{W}_{k}^{(n-1)}\mathbf{C}_{4,k}^{(n-1)},
\label{W_fe2_EKF_k}
\end{equation}%
and%
\begin{equation}
\mathbf{\eta }_{3,k}^{(n)}=\mathbf{W}_{k}^{(n-1)}\left[ \mathbf{C}%
_{4,k}^{(n-1)}\mathbf{w}_{2,k}+\mathbf{\eta }_{4,k}^{(n-1)}\right]
\label{w_fe2_EKF_k}
\end{equation}%
for $n>1$; here, $\mathbf{W}_{k}^{(n-1)}\triangleq \lbrack \mathbf{C}%
_{4,k}^{(n-1)}\mathbf{W}_{2,k}+\mathbf{I}_{D}]^{-1}$. Then, the message $%
\vec{m}_{3}^{(n)}\left( \mathbf{x}_{k}\right) $ (\ref{m_fe2_EKF_bis}) is
marginalized with respect to $\mathbf{x}_{k}^{(N)}$ in the F$_{1}{%
\rightarrow }$F$_{2}$ block; this results in the message 
\begin{eqnarray}
\vec{m}_{1}^{(n)}\left( \mathbf{x}_{k}^{(L)}\right) &\triangleq &\int \vec{m}%
_{3}^{(n)}\left( \mathbf{x}_{k}\right) \,d\mathbf{x}_{k}^{(N)}  \notag \\
&=&\mathcal{N(}\mathbf{x}_{k}^{(L)};\mathbf{\tilde{\eta}}_{1,k}^{(n)},%
\mathbf{\tilde{C}}_{1,k}^{(n)}),  \label{m_fe_L_EKF_2}
\end{eqnarray}%
where $\mathbf{\tilde{C}}_{1,k}^{(n)}$ and $\mathbf{\tilde{\eta}}%
_{1,k}^{(n)} $ are easily extracted from $\mathbf{C}_{3,k}^{(n)}$ (\ref%
{W_fe2_EKF_k}) and $\mathbf{\eta }_{3,k}^{(n)}$ (\ref{w_fe2_EKF_k}) for $n>1$
($\mathbf{C}_{3,k}^{(n)}$ (\ref{C_3_1}) and $\mathbf{\eta }_{3,k}^{(n)}$ (%
\ref{eta_3_1}) for $n=1$), respectively, since $\mathbf{x}_{k}^{(L)}$
consists of the first $D_{L}$ elements of $\mathbf{x}_{k}$.

2) \emph{First measurement update\ in} F$_{2}$ - This step concerns the
computation of the message (see Fig. \ref{Fig_7})%
\begin{equation}
\vec{m}_{2}^{(n)}\left( \mathbf{x}_{k}^{(N)}\right) =\vec{m}_{\mathrm{fp}%
}^{(n)}\left( \mathbf{x}_{k}^{(N)}\right) \,\vec{m}_{1}^{(n)}\left( \mathbf{x%
}_{k}^{(N)}\right) \,,  \label{m_fe1_k_PF}
\end{equation}%
that represents the \emph{first filtered pdf} computed by F$_{2}$. The
message $\vec{m}_{\mathrm{fp}}^{(n)}(\mathbf{x}_{k}^{(N)})$ conveys a set of
predicted particles; its $j-$th component is given by 
\begin{equation}
\vec{m}_{\mathrm{fp},j}^{(n)}\left( \mathbf{x}_{k}^{(N)}\right) =w_{p}\,%
\mathbf{\delta }\left( \mathbf{x}_{k}^{(N)}-\mathbf{x}_{k,j}^{(N)}\left[ n%
\right] \right)  \label{eq:message_fp_l_j_PF_bis}
\end{equation}%
and, consequently, coincides with $\vec{m}_{\mathrm{fp},j}(\mathbf{x}%
_{k}^{(N)})$ (\ref{eq:message_N_pred-1}) for $n=1$ only; note also that the
same weight is assigned to all the messages $\{\vec{m}_{\mathrm{fp},j}^{(n)}(%
\mathbf{x}_{k}^{(N)})\}$ for any $n$, since particle resampling is employed
in each iteration of this phase (see step 4)). The message (see Fig. \ref%
{Fig_7})%
\begin{equation}
\vec{m}_{1}^{(n)}(\mathbf{x}_{k}^{(N)})=\int f(\mathbf{y}_{k}|\mathbf{x}%
_{k}^{(N)},\,\mathbf{x}_{k}^{(L)})\,\vec{m}_{1}^{(n)}(\mathbf{x}%
_{k}^{(L)})\,d\mathbf{x}_{k}^{(L)},  \label{m_1_x_N}
\end{equation}%
instead, conveys \emph{measurement information}, that is the information
about $\mathbf{x}_{k}^{(N)}$ provided by $\mathbf{y}_{k}$ (\ref{eq:y_t}). In
particular, the value%
\begin{equation}
w_{1,k,j}^{(n)}=\mathcal{N}\left( \mathbf{y}_{k};\mathbf{\tilde{\eta}}%
_{1,k,j}^{(n)},\mathbf{\tilde{C}}_{1,k,j}^{(n)}\right)  \label{m_ms_j_k}
\end{equation}%
taken on by the message $\vec{m}_{1}^{(n)}(\,\mathbf{x}_{k}^{(N)})$ (\ref%
{m_1_x_N}) for $\mathbf{x}_{k}^{(N)}=\mathbf{x}_{k,j}^{(N)}[n]$ represents
the measurement-based \emph{weight} assigned to the $j-$th particle $\mathbf{%
x}_{k,j}^{(N)}[n]$;\ here,

\begin{equation}
\mathbf{\tilde{\eta}}_{1,k,j}^{(n)}=\mathbf{B}_{k,j}[n]\,\mathbf{\tilde{\eta}%
}_{1,k}^{(n)}+\mathbf{g}_{k,j}[n],  \label{eq:eta_fe_PF}
\end{equation}%
\begin{equation}
\mathbf{\tilde{C}}_{1,k,j}^{(n)}=\mathbf{B}_{k,j}[n]\,\mathbf{\tilde{C}}%
_{1,k}^{(n)}\,\left( \mathbf{B}_{k,j}[n]\right) ^{T}+\mathbf{C}_{e},
\label{eq:C_fe_PF}
\end{equation}%
$\mathbf{g}_{k,j}[n]\triangleq \mathbf{g}_{k}(\mathbf{x}_{k,j}^{(N)}[n])$
and $\mathbf{B}_{k,j}[n]\triangleq \mathbf{B}_{k}(\mathbf{x}_{k,j}^{(N)}[n])$%
. From (\ref{m_fe1_k_PF}), (\ref{eq:message_fp_l_j_PF_bis}) and (\ref%
{m_ms_j_k}) it is easily inferred that $\vec{m}_{2}^{(n)}(\mathbf{x}%
_{k}^{(N)})$ (\ref{m_fe1_k_PF}) conveys the same set of particles as $\vec{m}%
_{\mathrm{fp}}^{(n)}(\mathbf{x}_{k}^{(N)})$ and that its $j-$th component is%
\begin{equation}
\vec{m}_{2,j}^{(n)}\left( \mathbf{x}_{k}^{(N)}\right)
=w_{p}\,w_{1,k,j}^{(n)}\,\mathbf{\delta }\left( \mathbf{x}_{k}^{(N)}-\mathbf{%
x}_{k,j}^{(N)}\left[ n\right] \right) .  \label{eq:message_fe_l_j_PF_bis}
\end{equation}

3) \emph{Computation of the pseudo-measurements} \emph{for} F$_{2}$ - This
step is accomplished in the F$_{1}{\rightarrow }$F$_{2}$ block and aims at
computing the message $\vec{m}_{3}^{(n)}(\mathbf{x}_{k}^{(N)})$; this
conveys the statistical information about $\mathbf{x}_{k}^{(N)}$ that
originates from the pseudo-measurement $\mathbf{z}_{k}^{(N)}$ (\ref{z_N_l})
(further details about this message and its meaning are provided in Appendix %
\ref{app:A}). Actually, what is really required in the next step is the
value taken on by this message for $\mathbf{x}_{k}^{(N)}=\mathbf{x}%
_{k,j}^{(N)}\left[ n\right] $ (with $j=1,2,...,N_{p}$), because of the Dirac
delta function conveyed by the message $\vec{m}_{2,j}^{(n)}(\mathbf{x}%
_{k}^{(N)})$ (\ref{eq:message_fe_l_j_PF_bis}) and appearing in the \emph{%
right-hand side} (RHS) of Eq. (\ref{m_fe2_k_PF}); such a value is 
\begin{eqnarray}
w_{3,k,j}^{(n)} &=&\breve{D}_{k,j}^{(n)} \cdot \exp \left[ \frac{1}{2}\left(
\left( \mathbf{\check{\eta}}_{3,k,j}^{(n)}\right) ^{T}\mathbf{\check{W}}%
_{3,k,j}^{(n)}\,\mathbf{\check{\eta}}_{3,k,j}^{(n)}-\left( \mathbf{\check{%
\eta}}_{z,k,j}^{(n)}\right) ^{T}\right. \right.  \notag \\
&&\left. \left. \cdot \mathbf{\check{W}}_{z,k,j}^{(n)}\,\mathbf{\check{\eta}}%
_{z,k,j}^{(n)}-\left( \mathbf{f}_{k,j}^{(L)}[n]\right) ^{T}\mathbf{W}%
_{w}^{(L)}\mathbf{f}_{k,j}^{(L)}[n]\right) \right] ,  \notag \\
&&  \label{m_pm_x_N_l_ja}
\end{eqnarray}%
and represents a new weight to be assigned to $\mathbf{x}_{k,j}^{(N)}\left[ n%
\right] $, i.e. to the $j-$th particle of the set $S_{k}[n]$; here,%
\begin{equation}
\mathbf{\check{W}}_{3,k,j}^{(n)}\triangleq \left( \mathbf{\check{C}}%
_{3,k,j}^{(n)}\right) ^{-1}=\mathbf{\check{W}}_{z,k,j}^{(n)}+\mathbf{W}%
_{w}^{(L)},  \label{W_pm_x_N_l_j}
\end{equation}%
\begin{equation}
\mathbf{\check{w}}_{3,k,j}^{(n)}\triangleq \mathbf{\check{W}}_{3,k,j}^{(n)}\,%
\mathbf{\check{\eta}}_{3,k,j}^{(n)}=\mathbf{\check{w}}_{z,k,j}^{(n)}+\mathbf{%
W}_{w}^{(L)}\mathbf{f}_{k,j}^{(L)}[n],  \label{w_pm_x_N_l_j}
\end{equation}%
$\mathbf{W}_{w}^{(L)}\triangleq \lbrack \mathbf{C}_{w}^{(L)}]^{-1}$, $%
\mathbf{f}_{k,j}^{(L)}[n]\triangleq \mathbf{f}_{k}^{(L)}(\mathbf{x}%
_{k,j}^{(N)}[n])$, $\mathbf{\check{W}}_{z,k,j}^{(n)}\triangleq (\mathbf{%
\check{C}}_{z,k,j}^{(n)})^{-1}$, $\mathbf{\check{w}}_{z,k,j}^{(n)}\triangleq 
\mathbf{\check{W}}_{z,k,j}^{(n)}\mathbf{\check{\eta}}_{z,k,j}^{(n)}$, 
\begin{equation}
\mathbf{\check{C}}_{z,k,j}^{(n)}=\mathbf{C}_{w}^{(L)}+\mathbf{A}%
_{k,j}^{(L)}[n]\left[ \mathbf{\tilde{C}}_{3,k}^{(n)}-\mathbf{\tilde{C}}_{2,k}%
\right] \left( \mathbf{A}_{k,j}^{(L)}[n]\right) ^{T},  \label{Cov_pm_x_N_l_j}
\end{equation}%
\begin{equation}
\mathbf{\check{\eta}}_{z,k,j}^{(n)}=\mathbf{A}_{k,j}^{(L)}[n]\left[ \mathbf{%
\tilde{\eta}}_{3,k}^{(n)}-\mathbf{\tilde{\eta}}_{2,k}\right] +\mathbf{f}%
_{k,j}^{(L)}[n],  \label{eta_pm_x_N_l_j}
\end{equation}%
$\mathbf{A}_{k,j}^{(L)}[n]\triangleq \mathbf{A}_{k}^{(L)}(\mathbf{x}%
_{k,j}^{(N)}[n])$, $\breve{D}_{k,j}^{(n)}\triangleq \lbrack \det (\mathbf{%
\breve{C}}_{k,j}^{(n)})]^{-1/2}$, $\mathbf{\breve{C}}_{k,j}^{(n)}\triangleq 
\mathbf{\check{C}}_{z,k,j}^{(n)}+\mathbf{C}_{w}^{(L)}$, and $\mathbf{\tilde{%
\eta}}_{2,k}$ and $\mathbf{\tilde{C}}_{2,k}$ ($\mathbf{\tilde{\eta}}%
_{3,k}^{(n)}$ and $\mathbf{\tilde{C}}_{3,k}^{(n)}$)\ are extracted from $%
\mathbf{\eta }_{2,k}$ and $\mathbf{C}_{2,k}$ ($\mathbf{\eta }_{3,k}^{(n)}$
and $\mathbf{C}_{3,k}^{(n)}$), respectively (see Eqs. (\ref{m_fe1_EKF_tris})
and (\ref{m_fe2_EKF_bis})), since they refer to the first $D_{L}$ elements
of $\mathbf{x}_{k}$.

4) \emph{Second measurement update\ in }F$_{2}$ - In this step, the weights
of the particles forming the set $S_{k}[n]$\ are updated on the basis of the
weights $\{w_{3,k,j}^{(n)}\}$ computed in the previous step (see Eq. (\ref%
{m_pm_x_N_l_ja})). The new weight for the $j-$th particle $\mathbf{x}%
_{k,j}^{(N)}[n]$ is computed as 
\begin{equation}
w_{4,k,j}^{(n)}\triangleq w_{p}\cdot w_{1,k,j}^{(n)}\cdot w_{3,k,j}^{(n)}
\label{w_fe_2_x_N_l}
\end{equation}%
and combines the initial weight $w_{p}$ (originating from $\vec{m}_{\mathrm{%
fp},j}^{(n)}(\mathbf{x}_{k}^{(N)})$ (\ref{eq:message_fp_l_j_PF_bis})) with
the weights $w_{1,k,j}^{(n)}$ (\ref{m_ms_j_k}) and $w_{3,k,j}^{(n)}$ (\ref%
{m_pm_x_N_l_ja}) related to the measurement $\mathbf{y}_{k}$ (\ref{eq:y_t})
and the pseudo-measurement $\mathbf{z}_{k}^{(N)}$ (\ref{z_N_l}),
respectively. Note also that the weight $w_{4,k,j}^{(n)}$ (\ref{w_fe_2_x_N_l}%
) is conveyed by the message (see Fig. \ref{Fig_7})%
\begin{align}
\vec{m}_{4,j}^{(n)}\left( \mathbf{x}_{k}^{(N)}\right) & =\vec{m}%
_{2,j}^{(n)}\left( \mathbf{x}_{k}^{(N)}\right) \,\vec{m}_{3}^{(n)}\left( 
\mathbf{x}_{k}^{(N)}\right)  \label{m_fe2_k_PF} \\
& =w_{4,k,j}^{(n)}\,\delta \left( \mathbf{x}_{k}^{(N)}-\mathbf{x}%
_{k,j}^{(N)}[n]\right) ,  \label{m_fe_2_x_N_l_bis}
\end{align}%
that represents the $j-$th component of the message $\vec{m}_{4}^{(n)}(%
\mathbf{x}_{k}^{(N)})$ (with $j=1,2,...,N_{p}$).

Once all the weights $\{w_{4,k,j}^{(n)}\}$ are available, their
normalization is accomplished; this produces the normalised weights%
\begin{equation}
W_{4,k,j}^{(n)}\triangleq C_{k}^{(n)}\,w_{4,k,j}^{(n)},  \label{W_fe_2_x_N_l}
\end{equation}%
where $C_{k}^{(n)}\triangleq 1/\sum\limits_{j=1}^{N_{p}}w_{4,k,j}^{(n)}$.
The particles $\{\mathbf{x}_{k,j}^{(N)}[n]\}$ and their weights $%
\{W_{4,k,j}^{(n)}\}$ represent the\emph{\ second} \emph{filtered pdf} of $%
\mathbf{x}_{k}^{(N)}$ computed by F$_{2}$\ in the $n-$th iteration of the
considered recursion; consequently, the \emph{final} filtered pdf evaluated
by F$_{2}$ is represented by the particles $\{\mathbf{x}_{k,j}^{(N)}[n_{i}]%
\} $ ad their weights $\{W_{4,k,j}^{(n_{i})}\}$ computed in the last
iteration.

Resampling with replacement is now accomplished for the particle set $%
S_{k}[n]$ on the basis of the new weights $\{W_{4,k,j}^{(n)}\}$ (see Eq. (%
\ref{W_fe_2_x_N_l})). This entails that the $N_{p}$ particles $\{\mathbf{x}%
_{k,j}^{(N)}[n]\}$ and their associated weights $\{W_{4,k,j}^{(n)}\}$ are
replaced by the new particles $\{\mathbf{x}_{k,j}^{(N)}[n+1]\}$, forming
the\ set $S_{k}[n+1]$ and having identical weights (all equal to $%
w_{p}\triangleq 1/N_{p}$). Consequently, the effect of resampling can be
represented as turning the message $\vec{m}_{4,j}^{(n)}(\mathbf{x}%
_{k}^{(N)}) $ (\ref{m_fe_2_x_N_l_bis}) into the message 
\begin{equation}
\vec{m}_{4,j}^{(n)}\left( \mathbf{x}_{k}^{(N)}\right) =w_{p}\,\delta \left( 
\mathbf{x}_{k}^{(N)}-\mathbf{x}_{k,j}^{(N)}[n+1]\right) ,
\label{m_fe_2_x_N_l_tris}
\end{equation}%
with $j=1,2,..,N_{p}$.

5) \emph{Time update in} F$_{2}$ - In this step, the message $\vec{m}_{%
\mathrm{fp}}^{(n)}(\mathbf{x}_{k+1}^{(N)})$, conveying the predicted pdf of $%
\mathbf{x}_{k+1}^{(N)}$, is computed using the same method as RBPF (e.g.,
see \cite[Sec. IV, p. 1526]{Vitetta_2019}). For this reason, for any $j$,
the pdf (see Fig. \ref{Fig_7}) 
\begin{eqnarray}
&&\int \int f\left( \mathbf{x}_{k+1}^{(N)}\left\vert \mathbf{x}_{k}^{(L)},%
\mathbf{x}_{k}^{(N)}\right. \right) \cdot \vec{m}_{4,j}^{(n)}\left( \mathbf{x%
}_{k}^{(N)}\right) \vec{m}_{1}^{(n)}\left( \mathbf{x}_{k}^{(L)}\right) d%
\mathbf{x}_{k}^{(L)}d\mathbf{x}_{k}^{(N)}  \label{eq:message_5_Na} \\
&=&\mathcal{N}\left( \mathbf{x}_{k+1}^{(N)};\mathbf{\eta }_{3,k,j}^{(N)},%
\mathbf{C}_{3,k,j}^{(N)}\right) ,  \label{eq:message_5_Nb}
\end{eqnarray}%
representing a prediction of $\mathbf{x}_{k+1}^{(N)}$ conditioned on $%
\mathbf{x}_{k}^{(N)}=\mathbf{x}_{k,j}^{(N)}[n+1]$ is computed first; here, 
\begin{equation}
\mathbf{\eta }_{3,k,j}^{(N)}\triangleq \mathbf{A}_{k,j}^{(N)}\,\left[ n+1%
\right] \,\,\mathbf{\tilde{\eta}}_{1,k}^{(n)}+\mathbf{f}_{k,j}^{(N)}\left[
n+1\right] ,  \label{eq:eta_5_N}
\end{equation}%
\begin{equation}
\mathbf{C}_{3,k,j}^{(N)}\triangleq \mathbf{A}_{k,j}^{(N)}\,\left[ n+1\right]
\,\mathbf{\tilde{C}}_{1,k}^{(n)}\left( \mathbf{A}_{k,j}^{(N)}\left[ n+1%
\right] \right) ^{T}+\mathbf{C}_{w}^{(N)},  \label{eq:C_5_N}
\end{equation}%
$\mathbf{A}_{k,j}^{(N)}[n+1]\triangleq \mathbf{A}_{k}^{(N)}(\mathbf{x}%
_{k,j}^{(N)}[n+1])$ and $\mathbf{f}_{k,j}^{(N)}[n+1]\triangleq \mathbf{f}%
_{k}^{(N)}(\mathbf{x}_{k,j}^{(N)}[n+1])$. Then, the sample $\mathbf{\bar{x}}%
_{k+1,j}^{(N)}\left[ n+1\right] $ is drawn from the Gaussian function (\ref%
{eq:message_5_Nb}) and the weight $w_{p}$ is assigned to it; these
information are conveyed by the $j-$th component 
\begin{equation}
\vec{m}_{\mathrm{fp},j}^{(n)}\left( \mathbf{x}_{k+1}^{(N)}\right) =w_{p}\,%
\mathbf{\delta }\left( \mathbf{x}_{k+1}^{(N)}-\mathbf{\bar{x}}_{k+1,j}^{(N)}%
\left[ n+1\right] \right) ,  \label{eq:message_5_N_l+1}
\end{equation}%
of the message $\vec{m}_{\mathrm{fp}}^{(n)}(\mathbf{x}_{k+1}^{(N)})$.

6) \emph{Computation of the pseudo-measurements\ for} F$_{1}$ - This step is
accomplished in the F$_{2}{\rightarrow }$F$_{1}$ block and aims at
generating the message (see Fig. \ref{Fig_7}) 
\begin{equation}
\vec{m}_{4}^{(n)}\left( \mathbf{x}_{k}\right) =\mathcal{\mathcal{N}}\left( 
\mathbf{x}_{k};\mathbf{\eta }_{4,k}^{(n)},\mathbf{C}_{4,k}^{(n)}\right) ,
\label{message_pm_x_l_k}
\end{equation}%
that conveys the pseudo-measurement information exploited by F$_{1}$\ in its 
\emph{second measurement update} of the next iteration. The mean vector $%
\mathbf{\eta }_{4,k}^{(n)}$ is evaluated as 
\begin{equation}
\mathbf{\eta }_{4,k}^{(n)}=\left[ \left( \mathbf{\tilde{\eta}}%
_{4,k}^{(n)}\right) ^{T},\left( \mathbf{\check{\eta}}_{4,k}^{(n)}\right) ^{T}%
\right] ^{T},  \label{eta_pm_l_k}
\end{equation}%
where 
\begin{equation}
\mathbf{\tilde{\eta}}_{4,k}^{(n)}\triangleq \frac{1}{N_{p}}\sum_{j=1}^{N_{p}}%
\mathbf{\tilde{\eta}}_{4,k,j}^{(n)}  \label{eta_pm_l_L_k}
\end{equation}%
and 
\begin{equation}
\mathbf{\check{\eta}}_{4,k}^{(n)}\triangleq \frac{1}{N_{p}}\sum_{j=1}^{N_{p}}%
\mathbf{x}_{k,j}^{(N)}[n]  \label{eta_pm_l_N_k}
\end{equation}%
are a $D_{L}-$dimensional mean vector and a $D_{N}-$dimensional mean vector,
respectively. The covariance matrix $\mathbf{C}_{4,k}^{(n)}$, instead, is
computed as%
\begin{equation}
\mathbf{C}_{4,k}^{(n)}=\left[ 
\begin{array}{cc}
\mathbf{\tilde{C}}_{4,k}^{(n)} & \mathbf{\dot{C}}_{4,k}^{(n)} \\ 
\left( \mathbf{\dot{C}}_{4,k}^{(n)}\right) ^{T} & \mathbf{\check{C}}%
_{4,k}^{(n)}%
\end{array}%
\right] ,  \label{C_pm_l_k}
\end{equation}%
where 
\begin{equation}
\mathbf{\tilde{C}}_{4,k}^{(n)}\triangleq \frac{1}{N_{p}}\sum_{j=1}^{N_{p}}%
\mathbf{\tilde{r}}_{4,k,j}^{(n)}-\mathbf{\tilde{\eta}}_{4,k}^{(n)}\left( 
\mathbf{\tilde{\eta}}_{4,k}^{(n)}\right) ^{T},  \label{C_pm_l_L_k_bis}
\end{equation}%
is a $D_{L}\times D_{L}$ covariance matrix,%
\begin{equation}
\mathbf{\check{C}}_{4,k}^{(n)}\triangleq \frac{1}{N_{p}}\sum_{j=1}^{N_{p}}%
\mathbf{\check{r}}_{4,k,j}^{(n)}[n]-\mathbf{\check{\eta}}_{4,k}^{(n)}\left( 
\mathbf{\check{\eta}}_{4,k}^{(n)}\right) ^{T},  \label{C_pm_l_N_k}
\end{equation}%
is a $D_{N}\times D_{N}$\ covariance matrix and%
\begin{equation}
\mathbf{\dot{C}}_{4,k}^{(n)}\triangleq \frac{1}{N_{p}}\sum_{j=1}^{N_{p}}%
\mathbf{\dot{r}}_{4,k,j}^{(n)}-\mathbf{\tilde{\eta}}_{4,k}^{(n)}\left( 
\mathbf{\check{\eta}}_{4,k}^{(n)}\right) ^{T},  \label{C_pm_l_LN_k}
\end{equation}%
is $D_{L}\times D_{N}$ cross-covariance matrix. Moreover, $\mathbf{\tilde{r}}%
_{4,k,j}^{(n)}\triangleq \mathbf{\tilde{C}}_{4,k,j}^{(n)}+\mathbf{\tilde{\eta%
}}_{4,k,j}^{(n)}(\mathbf{\tilde{\eta}}_{4,k,j}^{(n)})^{T}$, $\mathbf{\check{r%
}}_{4,k,j}^{(N)}[n]\triangleq \mathbf{x}_{k,j}^{(N)}[n+1](\mathbf{x}%
_{k,j}^{(N)}[n+1])^{T}$ $\mathbf{\dot{r}}_{4,k,j}^{(n)}\triangleq \mathbf{%
\tilde{\eta}}_{4,k,j}^{(n)}(\mathbf{x}_{k,j}^{(N)}[n+1])^{T}$, the
covariance matrix $\mathbf{\tilde{C}}_{4,k,j}^{(n)}$ and the mean vector $%
\mathbf{\tilde{\eta}}_{4,k,j}^{(n)}$ are computed on the basis of the
associated precision matrix%
\begin{equation}
\mathbf{\tilde{W}}_{4,k,j}^{(n)}\triangleq \left( \mathbf{\tilde{C}}%
_{4,k,j}^{(n)}\right) ^{-1}=\left( \mathbf{A}_{k,j}^{(N)}\left[ n+1\right]
\right) ^{T}\mathbf{W}_{w}^{(N)}\mathbf{A}_{k,j}^{(N)}\left[ n+1\right]
\label{eq:W_pm_L_j1}
\end{equation}%
and of the associated transformed mean vector 
\begin{equation}
\mathbf{\tilde{w}}_{4,k,j}^{(n)}\triangleq \mathbf{\tilde{W}}_{4,k,j}^{(n)}%
\mathbf{\tilde{\eta}}_{4,k,j}^{(n)}=\left( \mathbf{A}_{k,j}^{(N)}\left[ n+1%
\right] \right) ^{T}\mathbf{W}_{w}^{(N)}\mathbf{z}_{k,j}^{(L)}\left[ n+1%
\right] ,  \label{eq:w_pm_L_j}
\end{equation}%
respectively, and%
\begin{equation}
\mathbf{z}_{k,j}^{(L)}\left[ n+1\right] \triangleq \mathbf{\bar{x}}%
_{k+1,j}^{(N)}[n+1]-\mathbf{f}_{k,j}^{(N)}[n+1].  \label{eq:z_L_evaluateda}
\end{equation}

The computation of $\vec{m}_{4}^{(n)}(\mathbf{x}_{k})$ (\ref%
{message_pm_x_l_k}) concludes step 6) and, consequently, the $n-$th
iteration of phase II. Then, if the iteration index $n$ is less than $n_{i}$%
, it is increased by one, so that a new iteration can be started by going
back to step 1); otherwise, phase III is accomplished.

\textbf{Phase III }- In this phase, the message (see Fig. \ref{Fig_7}) 
\begin{equation}
\vec{m}_{3}^{(n_{i}+1)}\left( \mathbf{x}_{k}\right) =\mathcal{\mathcal{N}}%
\left( \mathbf{x}_{k};\mathbf{\eta }_{3,k}^{(n_{i}+1)},\mathbf{C}%
_{3,k}^{(n_{i}+1)}\right) ,  \label{m_fe2_EKF_final}
\end{equation}%
conveying the final filtered pdf provided by F$_{1}$, is computed on the
basis of Eqs. (\ref{m_fe2_EKF})--(\ref{w_fe2_EKF_k}) as if a new iteration
(corresponding to $n=n_{i}+1$) was started. Then, if $k<t$, the output
messages $\vec{m}_{\mathrm{fp}}(\mathbf{x}_{k+1}^{(N)})$ and $\vec{m}_{%
\mathrm{fp}}\left( \mathbf{x}_{k+1}\right) $ (i.e., the new predicted
densities) are computed; otherwise, DBF processing is over, since the final
measurement has been processed. In the first case, the $j-$th component of $%
\vec{m}_{\mathrm{fp}}(\mathbf{x}_{k+1}^{(N)})$ is generated by F$_{1}$ as
(see Fig. \ref{Fig_5}) 
\begin{equation}
\vec{m}_{\mathrm{fp},j}\left( \mathbf{x}_{k+1}^{(N)}\right) =\vec{m}_{%
\mathrm{fp},j}^{(n_{i})}\left( \mathbf{x}_{k+1}^{(N)}\right)
\label{m_fp_N_l+1_final}
\end{equation}%
for $j=1$, $...$, $N_{p}$ (see Eq. (\ref{eq:message_fp_l_j_PF_bis})); this
means that the particle set $S_{k+1}[1]$ available at the beginning of the
next recursion consists of the particles $\{\mathbf{x}_{k+1,j}^{(N)}=\mathbf{%
\bar{x}}_{k+1,j}^{(N)}\left[ n_{i}+1\right] $; $j=1$, $2$, $...$, $N_{p}\}$.
Then, the predicted pdf $\vec{m}_{\mathrm{fp}}\left( \mathbf{x}_{k+1}\right) 
$ is computed by F$_{1}$ as (see Fig. \ref{Fig_7}) 
\begin{eqnarray}
\vec{m}_{\mathrm{fp}}\left( \mathbf{x}_{k+1}\right) &=&\int \,\tilde{f}%
\left( \mathbf{x}_{k+1}\left\vert \mathbf{x}_{k}\right. \right) \,\vec{m}%
_{3}^{(n_{i}+1)}\left( \mathbf{x}_{k}\right) \,d\mathbf{x}_{k}
\label{m_fp_l+1_whole} \\
&=&\mathcal{N}\left( \mathbf{x}_{k+1};\mathbf{\eta }_{\mathrm{fp},k+1},%
\mathbf{C}_{\mathrm{fp},k+1}\right) ,  \label{m_fp_l+1_whole_tris}
\end{eqnarray}%
where%
\begin{equation}
\mathbf{\eta }_{\mathrm{fp},k+1}\triangleq \mathbf{F}_{k}\,\mathbf{\eta }%
_{3,k}^{(n_{i}+1)}+\mathbf{u}_{k},  \label{eta_fp_l+1_whole1}
\end{equation}%
and%
\begin{equation}
\mathbf{C}_{\mathrm{fp},k+1}\triangleq \mathbf{C}_{w}+\mathbf{F}_{k}\,%
\mathbf{C}_{3,k}^{(n_{i}+1)}\mathbf{F}_{k}^{T}.  \label{C_fp_l+1_whole}
\end{equation}%
This concludes the $k-$th recursion of the DBF technique.

The algorithm described above needs a proper initialization. In our work, a
(known) Gaussian pdf $f(\mathbf{x}_{1})=\mathcal{\mathcal{N(}}\mathbf{x}_{1};%
\mathbf{\eta }_{1},\mathbf{C}_{1})$ is assumed for the initial $\mathbf{x}%
_{1}$; for this reason, DBF\ is initialised by setting $\vec{m}_{\mathrm{fp}%
}\left( \mathbf{x}_{1}\right) =f(\mathbf{x}_{1})$ for F$_{1}$ and by
sampling the pdf $f(\mathbf{x}_{1}^{(N)})$ (that results from the
marginalization of $f(\mathbf{x}_{1})$ with respect to $\mathbf{x}_{1}^{(L)}$%
) $N_{p}$ times in order to generate the\ initial particle set $S_{1}[1]=\{%
\mathbf{x}_{1,j}^{(N)},\,j=1,...,N_{p}\}$; then, the same weight ($%
w_{p}=1/N_{p}$) is assigned to each particle.

All the processing tasks accomplished by the DBF technique are summarized in
Algorithm 1. Note also that, at the end of the $k-$th recursion, estimates $%
\mathbf{\hat{x}}_{\mathrm{fe},k}^{(N)}$ and $\mathbf{\hat{x}}_{\mathrm{fe}%
,k}^{(L)}$ of $\mathbf{x}_{k}^{(N)}$ and $\mathbf{x}_{k}^{(L)}$,
respectively, can be evaluated as: a) $\mathbf{\hat{x}}_{\mathrm{fe}%
,k}^{(N)}=\sum_{j=1}^{N_{p}}W_{4,k,j}^{(n_{i})}\,\mathbf{x}%
_{k,j}^{(N)}[n_{i}]$ (see our comments following Eq. (\ref{W_fe_2_x_N_l}))
or $\mathbf{\hat{x}}_{\mathrm{fe},k}^{(N)}=\bar{\eta}_{3,k}^{(n_{i}+1)}$,
where $\bar{\eta}_{3,k}^{(n_{i}+1)}$ consists of the \emph{last} $D_{N}$
elements of $\eta _{3,k}^{(n_{i}+1)}$ (see Eq. (\ref{m_fe2_EKF_final})); b) $%
\mathbf{\hat{x}}_{\mathrm{fe},k}^{(L)}=\mathbf{\tilde{\eta}}%
_{3,k}^{(n_{i}+1)}$, where $\mathbf{\tilde{\eta}}_{3,k}^{(n_{i}+1)}$
consists of the \emph{first} $D_{L}$ elements of $\mathbf{\eta }%
_{3,k}^{(n_{i}+1)}$.

\begin{algorithm}
	\SetKw{a}{a-}
	\SetKw{b}{b-}
	\SetKw{c}{c-}
	\SetKw{d}{c1-}
	\SetKw{e}{c2-}
	\SetKw{f}{c3-}
	\SetKw{g}{d-}
	\SetKw{h}{e-}
	\SetKw{i}{f-}
	\SetKw{j}{g-}
	\SetKw{k}{h-}
		
	\nl\textbf{Initialisation:} 
	For $j=1$ to $N_{p}$: sample the pdf $f(\mathbf{x}_{1}^{(N)})$ to generate the particles $\mathbf{x}_{1,j}^{(N)}$ (forming the set $S_{1}[1]$), and
	assign the weight $w_{p}=1\/N_{p}$ to each of them. Set $\mathbf{W}_{\mathrm{fp},1}=\mathbf{W}_{1}=[\mathbf{C}_{1}]^{-1}$, $\mathbf{w}_{\mathrm{fp},1}=\mathbf{W}_{1}{\eta}_{1}$.
	
	\nl\textbf{Filtering:} For $k=1$ to $t$:
	
	\a \emph{First measurement update in} F$_1$: Compute $\mathbf{W}_{2,k}$ (\ref{W_fe1_EKFa}) and $\mathbf{w}_{2,k}$ (\ref{w_fe1_EKF}), $\mathbf{C}_{2,k}=[\mathbf{W}_{2,k}]^{-1}$ and $\mathbf{\eta }_{2,k}=\mathbf{C}_{2,k}\mathbf{w}_{2,k}$. Then, extract $\mathbf{\tilde{\eta}}_{2,k}$ and $\mathbf{\tilde{C}}_{2,k}$ from $\mathbf{\eta}_{2,k}$ and $\mathbf{C}_{2,k}$, respectively, and set $\mathbf{W}_{4,k}^{(0)}=\mathbf{0}_{D,D}$ and $\mathbf{w}_{4,k}^{(0)}=\mathbf{0}_{D}$.
	
	\For {$n=1$ to $n_{i}$}{
		
		\b  \emph{Second measurement update in} F$_1$: Compute $\mathbf{C}_{3,k}^{(n)}$  and $\mathbf{\eta}_{3,k}^{(n)}$ (see Eqs.  (\ref{C_3_1})--\ref{w_fe2_EKF_k}); then, extract $\mathbf{\tilde{\eta}}_{1,k}^{(n)}$ and $\mathbf{\tilde{C}}_{1,k}^{(n)}$ from $\mathbf{\eta}_{3,k}^{(n)}$ and $\mathbf{C}_{3,k}^{(n)}$, respectively.
		
		\c \emph{Measurement updates in} F$_2$:\\
		\For {$j=1$ to $N_{p}$} 
		{\d \emph{First measurement update}: compute $\mathbf{\tilde{\eta}}_{1,k,j}^{(n)}$ (\ref{eq:eta_fe_PF}), $\mathbf{\tilde{C}}_{1,k,j}^{(n)}$ (\ref{eq:C_fe_PF}) and $w_{1,k,j}^{(n)}$ (\ref{m_ms_j_k}).

		\e  \emph{Computation of the pseudo-measurements for} F$_2$: compute $\mathbf{\check{C}}_{z,k,j}^{(n)}$ (\ref{Cov_pm_x_N_l_j}), $\mathbf{\check{\eta}}_{z,k,j}^{(n)}$ (\ref{eta_pm_x_N_l_j}), $\mathbf{\check{W}}_{z,k,j}^{(n)}= [\mathbf{\check{C}}%
		_{z,k,j}^{(n)}]^{-1}$ and $\mathbf{\check{w}}_{z,k,j}^{(n)}= \mathbf{\check{W}}_{z,k,j}^{(n)}\mathbf{\check{\eta}}_{z,k,j}^{(n)}$. Then, compute $\mathbf{\check{W}}_{3,k,j}^{(n)}$ (\ref{W_pm_x_N_l_j}), $\mathbf{\check{w}}_{3,k,j}^{(n)}$ (\ref{w_pm_x_N_l_j}), $\mathbf{\check{C}}_{3,k,j}^{(n)} = [\mathbf{\check{W}}_{3,k,j}^{(n)}]^{-1}$ and ${\check{\eta}}_{3,k,j}^{(n)}=\mathbf{\check{C}}_{3,k,j}^{(n)} \mathbf{\check{w}}_{3,k,j}^{(n)}$. Finally, compute $ w_{3,k,j}^{(n)}$ (\ref{m_pm_x_N_l_ja}).
			
		\f \emph{Second measurement update}: compute $w_{4,k,j}^{(n)}$ (\ref{w_fe_2_x_N_l}). 	
	}
		\g \emph{Normalization of particle weights}: compute the normalized weights $\{W_{4,k,j}^{(n)}\}$ according to Eq. (\ref{W_fe_2_x_N_l}).
			
		\h \emph{Resampling with replacement}: generate the new particle set $S_k[n+1] = \{\mathbf{x}_{k,j}^{(N)}[n+1]\}$ by resampling $S_k[n] $ on the basis of the weights $\{W_{4,k,j}^{(n)}\}$.
			
		\i  \emph{Time update in} F$_2$: For $j=1$ to $N_{p}$: Compute $\mathbf{\eta}_{3,k,j}^{(N)}$ (\ref{eq:eta_5_N}) and $\mathbf{C}_{3,k,j}^{(N)}$ (\ref{eq:C_5_N}), and sample the pdf $\mathcal{N}( 
		\mathbf{x}_{k+1}^{(N)};\mathbf{\eta}_{3,k,j}^{(N)},\mathbf{C}_{3,k,j}^{(N)})$ to generate the new particle $\mathbf{x}_{k+1,j}^{(N)}[n+1]$.
		
		\j \emph{Computation of the pseudo-measurements for} F$_1$: For $j=1$ to $N_{p}$: Compute $\mathbf{z}_{k,j}^{(L)}[n+1]$ (\ref{eq:z_L_evaluateda}), $\mathbf{\tilde{W}}_{4,k,j}^{(n)} $ (\ref{eq:W_pm_L_j1}) and $\mathbf{\tilde{w}}_{4,k,j}^{(n)} $ (\ref{eq:w_pm_L_j}), $\mathbf{\tilde{C}}_{4,k,j}^{(n)}=[\mathbf{\tilde{W}}_{4,k,j}^{(n)}]^{-1}$ and ${\tilde{\eta}}_{4,k,j}^{(n)} =\mathbf{\tilde{C}}_{4,k,j}^{(n)} \mathbf{\tilde{w}}_{4,k,j}^{(n)}$. Finally, compute $\mathbf{\eta }_{4,k}^{(n)} $ (\ref{eta_pm_l_k}) and $\mathbf{C}_{4,k}^{(n)} $ (\ref{C_pm_l_k}) (according to Eqs. (\ref{eta_pm_l_L_k})-(\ref{eta_pm_l_N_k}) and (\ref{C_pm_l_L_k_bis})-(\ref{C_pm_l_LN_k}), respectively).	
		}
		
		\k \emph{Compute forward predictions} (if $k<t$): For $j=1$ to $N_{p}$: set $\mathbf{x}_{k+1,j}^{(N)}=\mathbf{\bar{x}}_{k+1,j}^{(N)}[n_{i}]$ (these particles form the set $S_{k+1}[1]$). %
		Then, compute  $\mathbf{C}_{3,k}^{(n_{i}+1)}$ (\ref{W_fe2_EKF_k}) %
		and $\mathbf{\eta}_{3,k}^{(n_{i}+1)}$ (\ref{w_fe2_EKF_k}). Finally, compute $\mathbf{\eta }_{\mathrm{fp},k+1}$ (\ref{eta_fp_l+1_whole1}), $\mathbf{C}_{\mathrm{fp},k+1}$ (\ref{C_fp_l+1_whole}), $\mathbf{W}_{\mathrm{fp},k+1}=[\mathbf{C}_{\mathrm{fp},k+1}]^{-1}$ and $\mathbf{w}_{\mathrm{fp},k+1}=$ $\mathbf{W}_{\mathrm{fp},k+1} {\eta}_{\mathrm{fp},k+1}$.

		\caption{Dual Bayesian Filtering}
	\end{algorithm}

Following the same line of reasoning, a filtering method similar to DBF can
be developed for case \textbf{C.2}, i.e. for the second case considered in
the previous paragraph. Details are omitted for space limitations; however,
the relevant differences between this method (called \emph{simplified} DBF,
SDBF, in the following) and the DBF technique can be summarised as follows:

1) In phase I, $\mathbf{x}_{k}^{(N)}=\mathbf{\hat{x}}_{\mathrm{fp},k}^{(N)}$
is assumed in computing the \emph{first filtered pdf of} of $\mathbf{x}%
_{k}^{(L)}$, where $\mathbf{\hat{x}}_{\mathrm{fp},k}^{(N)}$ denotes the
prediction of $\mathbf{x}_{k}^{(N)}$ evaluated on the basis of the message $%
\vec{m}_{\mathrm{fp}}(\mathbf{x}_{k}^{(N)})$ (\ref{eq:eq:message_fp_PF})
provided by F$_{2}$.

2) In phase II, the message $\vec{m}_{4}^{(n)}(\mathbf{x}_{k})$ (\ref%
{message_pm_x_l_k}) is replaced by 
\begin{equation}
\vec{m}_{4}^{(n)}\left( \mathbf{x}_{k}^{(L)}\right) =\mathcal{\mathcal{N}}%
\left( \mathbf{x}_{k}^{(L)};\mathbf{\tilde{\eta}}_{4,k}^{(n)},\mathbf{\tilde{%
C}}_{4,k}^{(n)}\right) ,  \label{m_4_new}
\end{equation}%
since the pseudo-measurements computed in the F$_{2}{\rightarrow} $F$_{1}$
block refer to the linear state component only; here, $\mathbf{\tilde{\eta}}%
_{4,k}^{(n)}$ and $\mathbf{\tilde{C}}_{4,k}^{(n)}$ are given by Eqs. (\ref%
{eta_pm_l_L_k}) and (\ref{C_pm_l_L_k_bis}), respectively.

3) In phase III, $\mathbf{x}_{k}^{(N)}=\mathbf{\hat{x}}_{\mathrm{fe}%
,k}^{(N)} $ is assumed in computing the \emph{prediction} of $\mathbf{x}%
_{k+1}^{(L)}$, where $\mathbf{\hat{x}}_{\mathrm{fe},k}^{(N)}$ denotes the
estimate of $\mathbf{x}_{k}^{(N)}$ evaluated on the basis of the \emph{final}
filtered pdf computed by F$_{2}$.

\subsection{Computational complexity\label{comp_compl}}

The computational cost of the DBF and SDBF techniques has been carefully
assessed in terms of number of \emph{floating operations} (flops) to be
executed in each of their recursions. The general criteria adopted in
estimating the computational cost of an algorithm are the same as those
illustrated in \cite[App. A, p. 5420]{Hoteit_2016} and are not repeated here
for space limitations. A detailed analysis of the cost required by each task
accomplished by the DBF and the SDBF techniques is provided in Appendix \ref%
{app:CDBFA}. Our analysis leads to the conclusion that the computational
cost of the DBF and of the SDBF are approximately of order $\mathcal{O}%
(N_{DBF})$ and $\mathcal{O}(N_{SDBF})$, respectively, with 
\begin{eqnarray}
N_{DBF} &=&2PD^{2}+4P^{2}D+16D^{3}/3+14n_{i}D^{3}/3  \notag \\
&&+n_{i}\cdot N_{p}(2PD_{L}^{2}+2P^{2}D_{L}+2P^{3}/3  \notag \\
&&+6D_{L}^{3}+6D_{L}D_{N}^{2}+4D_{L}^{2}D_{N}+D_{N}^{3}/3)
\label{TF_comp_compl}
\end{eqnarray}%
and%
\begin{eqnarray}
N_{SDBF} &=&2PD_{L}^{2}+4P^{2}D_{L}+16D_{L}^{3}/3+14n_{i}D_{L}^{3}/3  \notag
\\
&&+n_{i}\cdot N_{p}(2PD_{L}^{2}+2P^{2}D_{L}+2P^{3}/3  \notag \\
&&+6D_{L}^{3}+6D_{L}D_{N}^{2}+4D_{L}^{2}D_{N}+D_{N}^{3}/3).
\label{STF_comp_compl}
\end{eqnarray}%
Each of the last two expressions has been derived as follows. First, the
costs of all the tasks identified in Appendix \ref{app:CDBFA} have been
summed (see Eqs. (\ref{C_F1_MU1})-(\ref{C_F2_TU})); then, the resulting
expression has been simplified, keeping only the dominant contributions due
to matrix inversions, matrix products and Cholesky decompositions and
discarding all the contributions that originate from the evaluation of the
matrices $\mathbf{A}_{k}^{(Z)}(\mathbf{x}_{k}^{(N)})$ (with $Z=L$ and $N$), $%
\mathbf{F}_{k}$, $\mathbf{H}_{k}$ and $\mathbf{B}_{k}$ and the functions $%
\mathbf{f}_{k}^{(Z)}(\mathbf{x}_{k}^{(N)})$ (with $Z=L$ and $N$), $\mathbf{f}%
_{k}(\mathbf{x}_{k})$ and $\mathbf{g}_{k}(\mathbf{x}_{k}^{(N)})$. Moreover,
the complexity of particle resampling has been ignored. A similar approach
has been followed for EKF, for RBPF and for the MPF technique described in 
\cite{Djuric_2013}; their complexities are approximately of order $\mathcal{O%
}(N_{EKF})$, $\mathcal{O}(N_{RBBF})$ and $\mathcal{O}(N_{MPF})$,
respectively, with 
\begin{equation}
N_{EKF}=2PD^{2}+2P^{2}D+2P^{3}/3+6D^{3},  \label{EKF_comp_compl}
\end{equation}%
\begin{eqnarray}
N_{RBPF} &=&N_{p}(4PD_{L}^{2}+6P^{2}D_{L}+2P^{3}/3+6D_{L}^{3}  \notag \\
&&+4D_{L}^{2}D_{N}+6D_{L}D_{N}^{2}+D_{N}^{3}/3)  \label{RBPF_comp_compl}
\end{eqnarray}%
and%
\begin{equation}
N_{MPF}=n(2\,M\,L\,d_{y}^{3}/3+M\,d_{x,i}^{3}/3);  \label{MPF_comp_compl}
\end{equation}%
note that the symbols appearing in the last formula are the same as those
defined in ref. \cite{Djuric_2013}. A detailed derivation of the eqs. (\ref%
{EKF_comp_compl})-(\ref{MPF_comp_compl}) is provided in the Appendices \ref%
{app:CEKF}-\ref{app:CMPF}.

It is important to keep in mind that a comparison among the computational
costs listed above does not fully account for the gap that can be observed
in the execution speed of the corresponding algorithms. In fact, distinct
filtering techniques may have substantially smaller memory requirements and,
as evidenced by our numerical results, this may influence their overall
execution speed. For instance, the DBF/SDBF techniques need to store the
state estimates and predictions generated by a single extended Kalman
filter, whereas RBPF needs to memorise those computed by a bank of $N_{p}$
Kalman filters running in parallel. Finally, it is worth stressing that $%
N_{DBF}$ (\ref{TF_comp_compl}) and $N_{SDBF}$ (\ref{STF_comp_compl}) exhibit
a linear dependence on the parameter $n_{i}$. Actually, in our computer
simulations, $n_{i}=1$ has been always selected, since marginal improvements
have been obtained by increasing $n_{i}$ beyond unity.

\section{Numerical Results\label{num_results}}

In this section we first compare, in terms of accuracy and execution time,
the DBF and SDBF techniques with an extended Kalman filter (corresponding to
F$_{1}$ of the DBF technique) and the RBPF technique (corresponding to the
combination of F$_{2}$ of the DBF technique with a bank of $N_{p}$ Kalman
filters) for a specific CLG SSM, denoted SSM \#1. This SSM is very similar
to the dynamic model described in \cite[Par. VII-A, p. 1531]{Vitetta_2019},
and refers to an agent moving on a plane and whose state is defined as $%
\mathbf{x}_{k}\triangleq \lbrack \mathbf{p}_{k}^{T},\mathbf{v}_{k}^{T}]^{T}$%
; here, $\mathbf{v}_{k}\triangleq \lbrack v_{x,k},v_{y,k}]^{T}$ and $\mathbf{%
p}_{k}\triangleq \lbrack p_{x,k},p_{y,k}]^{T}$ (corresponding to $\mathbf{x}%
_{k}^{(N)}$ and $\mathbf{x}_{k}^{(L)}$, respectively) represent the agent 
\emph{velocity} and its \emph{position}, respectively (their components are
expressed in m/s and in m, respectively). The dynamic models (see \cite[eqs.
(67)-(68), p. 1531]{Vitetta_2019}) 
\begin{equation}
\mathbf{v}_{k+1}=\rho \mathbf{v}_{k}+\left( 1-\rho \right) \mathbf{n}%
_{v,k}+\,\mathbf{a}\left( \mathbf{p}_{k},\mathbf{v}_{k}\right) T_{s},
\label{mod_1_v}
\end{equation}%
and%
\begin{equation}
\mathbf{p}_{k+1}=\mathbf{p}_{k}+\mathbf{v}_{k}T_{s}+\frac{1}{2}\mathbf{a}%
\left( \mathbf{p}_{k},\mathbf{v}_{k}\right) T_{s}^{2}+\mathbf{n}_{p,k}
\label{mod_1_p}
\end{equation}%
are adopted for the agent velocity and position, respectively; here, $\rho $
is a \emph{forgetting factor} ($0<\rho <1$), $T_{s}$ is the sampling
interval, $\{\mathbf{n}_{v,k}\}$ and $\{\mathbf{n}_{p,k}\}$ are mutually
independent \emph{additive white Gaussian noise} (AWGN) processes (whose
elements are characterized by the covariance matrices $\mathbf{I}_{2}$ and $%
\sigma _{p}^{2}\,\mathbf{I}_{2}$, respectively),%
\begin{equation}
\mathbf{a}\left( \mathbf{p}_{k},\mathbf{v}_{k}\right) =-(a_{0}/d_{0})\mathbf{%
p}_{k}-\tilde{a}_{0}f_{v}\left( \left\Vert \mathbf{v}_{k}\right\Vert \right) 
\mathbf{u}_{v,k}.  \label{acc_SSM1}
\end{equation}%
is the acceleration associated with position/velocity-dependent forces, $%
a_{0}$ and $\tilde{a}_{0}$ are scale factors (both expressed in m/s$^{2}$), $%
d_{0}$ is a \emph{reference distance, }$\mathbf{u}_{v,k}\triangleq \mathbf{v}%
_{k}/\left\Vert \mathbf{v}_{k}\right\Vert $ is the versor (i.e., the vector
of unit norm) associated with $\mathbf{v}_{k}$ and $f_{v}\left( x\right)
=(x/v_{0})^{3}$ is a continuous, differentiable and dimensionless function
(the parameter $v_{0}$ represents a \emph{reference velocity}). Moreover,
the measurement model 
\begin{equation}
\mathbf{y}_{k}=[\mathbf{p}_{k}^{T}\,\left\Vert \mathbf{v}_{k}\right\Vert
]^{T}+\mathbf{e}_{k},  \label{mod_1_y}
\end{equation}%
is adopted; here, $\{\mathbf{e}_{k}\}$ is an AWGN process, whose elements
are characterized by the covariance matrix $\mathbf{C}_{e}=$diag$(\sigma
_{e,p}^{2},\sigma _{e,p}^{2},\sigma _{e,v}^{2})$.

In our computer simulations, the estimation accuracy of the considered
filtering techniques for SSM\#1 has been assessed by evaluating two \emph{%
root mean square errors} (RMSEs), one for the linear state component, the
other for the nonlinear one, over an observation interval lasting $T=300$ $%
T_{s}$; these are denoted $RMSE_{L}($alg$)$ (m) and $RMSE_{N}($alg$)$ (m/s)
respectively, where `alg' denotes the algorithm these parameters refer to
(note also that $RMSE_{N}($DBF$)$ is computed on the basis of the estimate
of $\mathbf{v}_{k}$ generated by F$_{2}$, since this was found to be
slightly more accurate than that evaluated by F$_{1}$). Our assessment of 
\emph{computational requirements} is based, instead, on comparing $N_{DBF}$ (%
\ref{TF_comp_compl}), $N_{SDBF}$ (\ref{STF_comp_compl}), $N_{EKF}$ (\ref%
{EKF_comp_compl}) and $N_{RBPF}$ (\ref{RBPF_comp_compl}), and on assessing
the average \emph{execution time} required by each algorithm over the whole
observation interval. Moreover, the following values have been selected for
the parameters of SSM\#1: $\rho =0.99$, $T_{s}=0.1$ s, $\sigma _{p}$ $=0.01$
m, $\sigma _{e,p}=5\cdot 10^{-2}$ m, $\sigma _{e,v}=5\cdot 10^{-2}$ m/s, $%
a_{0}=1.5$ m/s$^{2}$, $d_{0}=0.5$ m, $\tilde{a}_{0}=0.05$ m/s$^{2}$ and $%
v_{0}=1$ m/s (the initial position $\mathbf{p}_{0}\triangleq \lbrack
p_{x,0},p_{y,0}]^{T}$ and the initial velocity $\mathbf{v}_{0}\triangleq
\lbrack v_{x,0},v_{y,0}]^{T}$ have been set to $[5$ m$,8$ m$]^{T}$ and $[4$
m/s$,$ $4$ m/s$]^{T}$, respectively). These values ensure that: a) the two
components of the position vector are represented by fast and damped
oscillations in the observation interval; b) the time variations of the
state vector can be accurately tracked by RBPF.

\begin{figure}[tbp]
\centering
\includegraphics[width=0.5\textwidth]{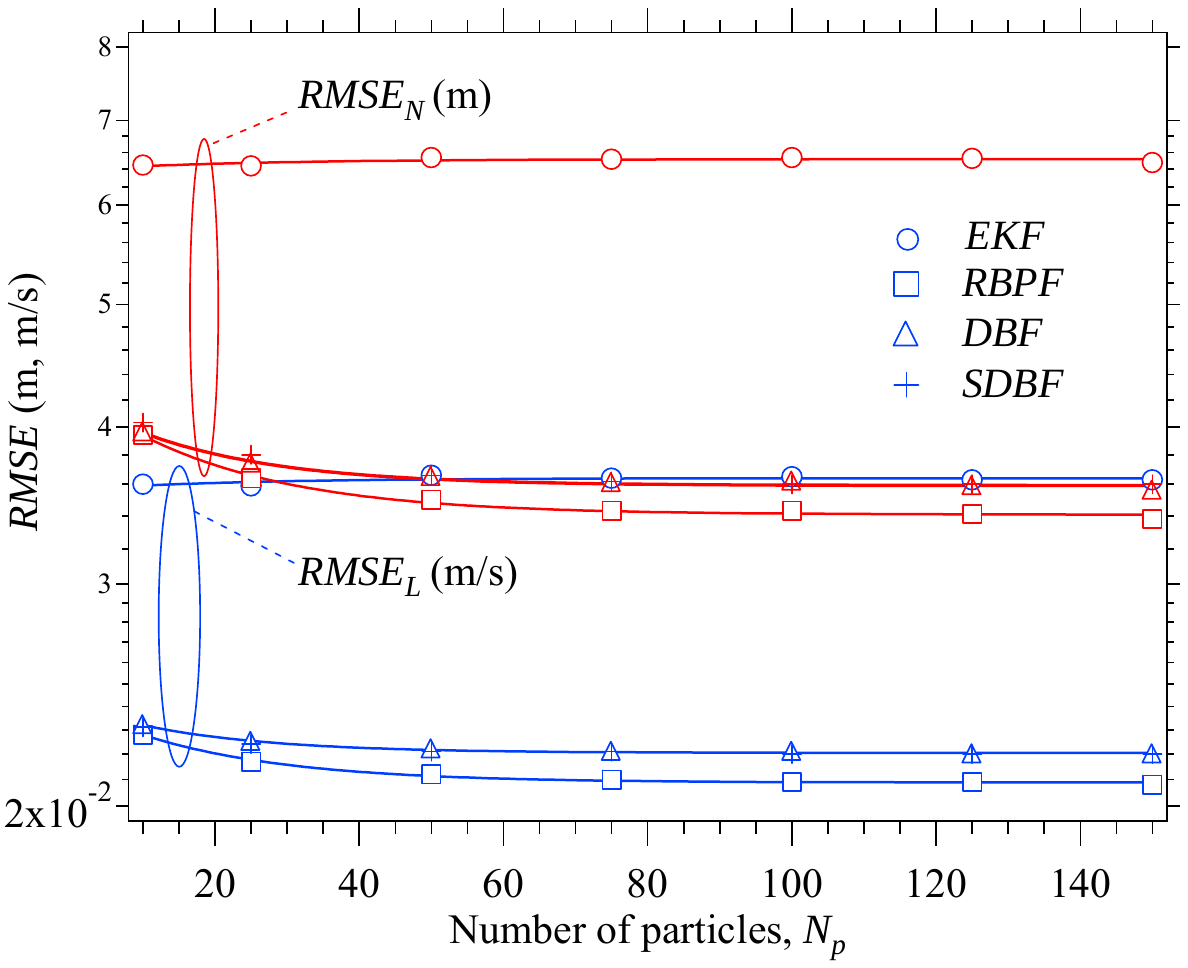}
\caption{RMSE performance versus $N_{p}$ for the linear component ($RMSE_{L}$%
; blue curves) and the nonlinear component ($RMSE_{N}$; red curves) of
system state (SSM\#1); EKF, RBPF, DBF and SDBF are considered.}
\label{Fig_1_sim}
\end{figure}

Some numerical results showing the dependence of $RMSE_{L}$ and $RMSE_{N}$
on the number of particles ($N_{p}$) for RBPF, EKF, DBF and SDBF are
illustrated in Fig. \ref{Fig_1_sim} (simulation results are indicated by
markers, whereas continuous lines are drawn to fit them, so facilitating the
interpretation of the available data); in this case, $n_{i}=1$ has been
selected for DBF/SDBF and the range $[10,150]$ has been considered for $%
N_{p} $. These results show that:

1) The EKF\ technique is appreciably outperformed by the other three
filtering algorithms in terms of both $RMSE_{L}$ and $RMSE_{N}$ for any
value of $N_{p}$; for instance, $RMSE_{L}($EKF$)$ ($RMSE_{N}($EKF$)$) is
about $1.65$\ ($1.80$) times larger than $RMSE_{L}($DBF$)$ ($RMSE_{N}($DBF$)$%
) for $N_{p}=100$.

2) DBF/SDBF perform slightly worse than RBPF for the same value of $N_{p}$
(for instance, $RMSE_{L}($DBF$)$ and $RMSE_{N}($DBF$)$ are about $5\%$
larger than the corresponding quantities evaluated for RBPF).

3) No real improvement in terms of $RMSE_{L}$ and $RMSE_{N}$ is found for $%
N_{p}\gtrsim 100$, if RBPF, DBF or SDBF are employed.

4) The SDBF performs very similarly as DBF; for this reason, in this
specific case, the presence of redundancy in the DBF does not allow to
achieve a better estimation accuracy.

\begin{figure}[tbp]
\centering
\includegraphics[width=0.55\textwidth]{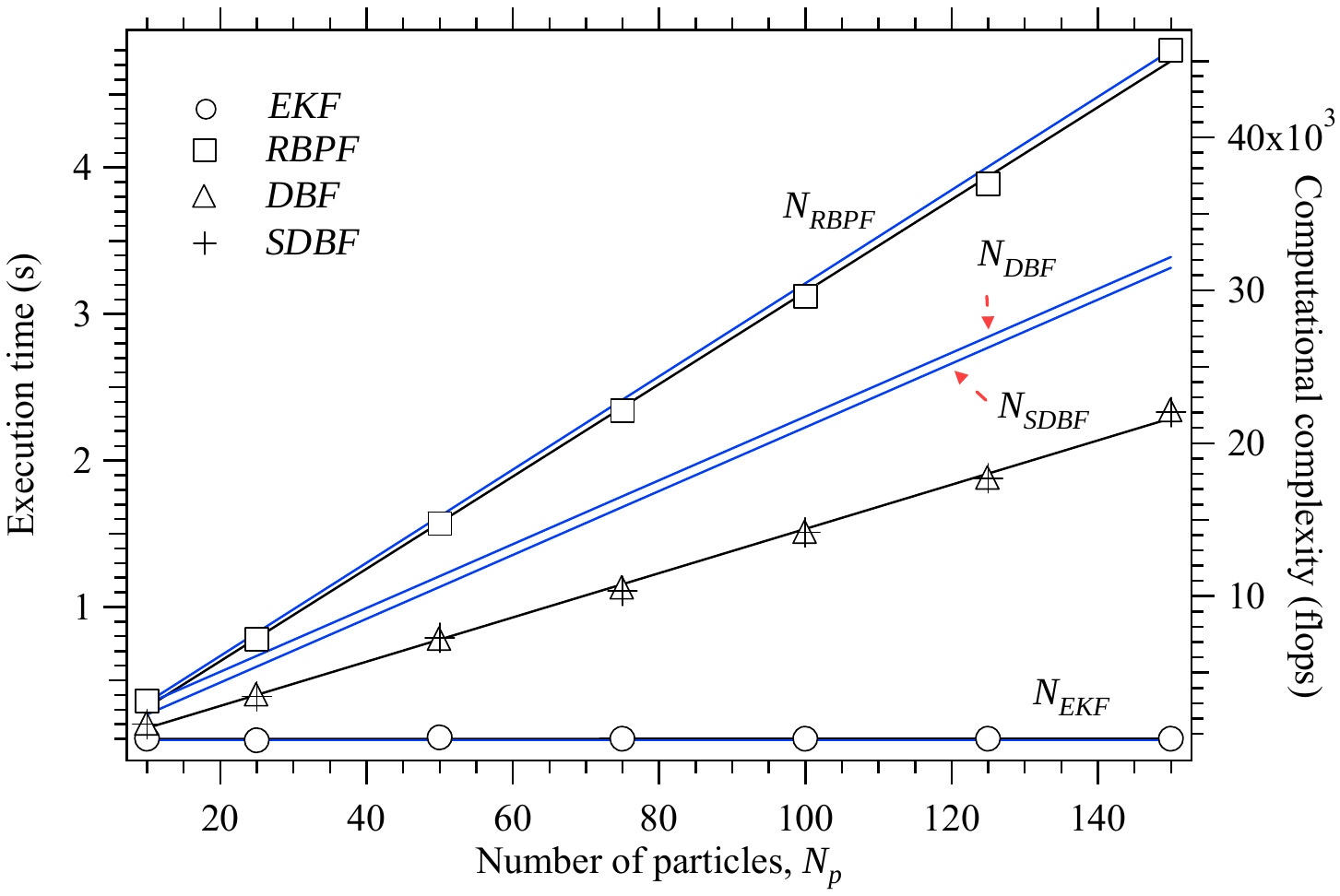}
\caption{ Execution time (black curves) and computational complexity (blue
curves) versus $N_{p}$ for EKF, RBPF, DBF and SDBF; SSM\#1 is considered. }
\label{Fig_2_sim}
\end{figure}

Despite their similar accuracies, RBPF, DBF and SDBF are characterized by
different computational complexities and execution times. This is evidenced
by the numerical results appearing in Fig. \ref{Fig_2_sim} and showing the
dependence of the execution time and the computational complexity on $N_{p}$
for the considered filtering algorithms. For instance, from these results it
is easily inferred that the DBF complexity is about $0.71$ times smaller
than that of RBPF for $N_{p}=100$; however, the gap in terms of execution
time is even larger mainly for the reasons illustrated at the end of
Paragraph \ref{comp_compl} (in particular, the execution time for the DBF is
approximately $0.61\ $times smaller than that required by RBPF). Moreover,
the results shown in Figs. \ref{Fig_1_sim}-\ref{Fig_2_sim} lead to the
conclusion that, in the considered scenario, DBF/SDBF achieve a better
accuracy-complexity tradeoff than RBPF.

The second SSM considered in this work has been inspired by refs. \cite%
{Closas_2012} and \cite{Djuric_2007}. In fact, it refers to a sensor network
employing $P$ sensors placed on the vertices of a square grid (partitioning
a square area whose side is equal to $l$ m) and receiving the reference
signals radiated, at the same power level and at the same frequency, by $N$
independent targets moving on a plane. Each target evolves according to the
motion model described by Eqs. (\ref{mod_1_v})-(\ref{mod_1_p}) with $\mathbf{%
a}\left( \mathbf{p}_{k},\mathbf{v}_{k}\right) =0$ for any $k$. In this case,
the considered SSM\ (denoted SSM\#2 in the following) refers to the whole
set of targets and its state vector $\mathbf{x}_{k}$ results from the
ordered concatenation of the vectors $\{\mathbf{x}_{k}^{(i)}$; $i=1$, $2$, $%
...$, $N\}$, where $\mathbf{x}_{k}^{(i)}\triangleq \lbrack (\mathbf{v}%
_{k}^{(i)})^{T},(\mathbf{p}_{k}^{(i)})^{T}]^{T}$, and \ $\mathbf{v}_{k}^{(i)}
$ and $\mathbf{p}_{k}^{(i)}$ represent the $i-$th target \emph{velocity} and
the \emph{position}, respectively. Moreover, the following additional
assumptions have been made about this SSM: 1) the process noises $\mathbf{n}%
_{p,k}^{(i)}$ and $\mathbf{n}_{v,k}^{(i)}$, affecting the $i-$th target
position and velocity, respectively, are given by $\mathbf{n}%
_{p,k}^{(i)}=(T_{s}^{2}/2)\,\mathbf{n}_{a,k}^{(i)}$ and $\mathbf{n}%
_{v,k}^{(i)}=T_{s}\mathbf{\,n}_{a,k}^{(i)}$, where $\{\mathbf{n}%
_{a,k}^{(i)}\}$ is two-dimensional AWGN, representing a random acceleration
and\ having covariance matrix $\sigma _{a}^{2}\,\mathbf{I}_{2}$ (with $i=1$, 
$2$, $...$, $N$); 2) the measurement acquired by the $q-$th sensor (with $q=1
$, $2$, $...$, $P$) in the $k$-th observation interval is given by 
\begin{equation}
y_{q,k}=10\,\mathrm{log}_{10}\left( {\Psi }\sum_{i=1}^{N}\frac{\,d_{0}^{2}}{%
\left\vert \left\vert \mathbf{s}_{q}-\mathbf{p}_{k}^{(i)}\right\vert
\right\vert ^{2}}\right) +e_{k},  \label{mod_2_y}
\end{equation}%
where the measurement noise $\{e_{k}\}$ is AWGN with variance $\sigma
_{e}^{2}$, ${\Psi }$ denotes the normalised power received by each sensor
from any target at a distance $d_{0}$ from the sensor itself and $\mathbf{s}%
_{q}$ is the position of the considered sensor; 3) the overall measurement
vector $\mathbf{y}_{k}$ results from the ordered concatenation of the
measurements $\{y_{q,k}$; $q=1$, $2$, $...$, $P\}$ and, consequently,
provides information about the position only; 4) the initial position $%
\mathbf{p}_{0}^{(i)}\triangleq \lbrack p_{x,0}^{(i)},p_{y,0}^{(i)}]^{T}$ and
the initial velocity $\mathbf{v}_{0}^{(i)}\triangleq \lbrack
v_{x,0}^{(i)},v_{y,0}^{(i)}]^{T}$ of the $i-$th target have been randomly
selected (with $i=1$, $2$, $...$, $N$). As far as the last point is
concerned, it is important to mention that, in our computer simulations,
distinct targets have been placed in different squares of the partioned area
in a random fashion; moreover, the initial velocity of each target has been
randomly selected within the interval $(v_{\mathrm{min}},v_{\mathrm{max}})$
in order to ensure that the trajectories of distinct targets do not cross
each other in the observation interval.

The following values have been selected for the parameters of SSM\#2: $P=25$%
, $l=10^{3}$ m, $T_{s}=1$ s, $\rho =1$, $\sigma _{a}^{2}=0.1$ m/s$^{2}$, $%
\sigma _{e}^{2}=-35$ dB, ${\Psi =1}$, $d_{0}=1$ m, $v_{\mathrm{min}}=0$ m/s
and $v_{\mathrm{min}}=0.1$ m/s. Moreover, a number $N$ of targets ranging
from $1$ to $5$ has been observed for $T=120$ $T_{s}$ s.

Our computer simulations for SSM\#2 have aimed at evaluating: a) the
accuracy achieved by different filtering algorithms in tracking the position
of $N$ targets; b) the probability that each filtering algorithm \emph{%
diverges} in the considered observation interval (this parameter is denoted $%
P_{FD}$ in the following). In practice, the accuracy achieved in position
tracking has been assessed by estimating the RMSE characterizing the whole
set $\{\mathbf{p}_{k}^{(i)}$; $i=1$, $2$, $...$, $N\}$ over each instant of
the considered observation interval; note that, if the $i-$th target is
considered, its position $\mathbf{p}_{k}^{(i)}$ represents the \emph{%
nonlinear} component of the associated substate $\mathbf{x}_{k}^{(i)}$,
because of the nonlinear dependence of $\mathbf{y}_{k}$ on it (see Eq. (\ref%
{mod_2_y})). On the other hand, the probability $P_{FD}$ has been assessed
by carefully identifying all the simulation runs in which the tracking of at
least one of the $N$ targets fails. Moreover, the tracking accuracy and the
probability of divergence have been evaluated for the following six
filtering techniques: 1) EKF; 2) RBPF; 3) the MPF technique developed in 
\cite{Djuric_2013} and based on the interconnection of $N$ identical
particle filters (one for each target); 4) DBF; 5) SDBF; 6) a novel
filtering algorithm based on the interconnection of $N_{F}=N+1$ filters and
dubbed MBF \emph{algorithm} (MBFA). The last algorithm involves the
interconnection of an extended Kalman filter with $N$ particle filters, each
representing the filtered/predicted pdfs of a two-dimensional vector through 
$\tilde{N}_{p}$ weighted particles. More specifically, the $i-$th particle
filter estimates the position $\mathbf{p}_{k}^{(i)}$ of the $i-$th target
(with $i=1$, $2$, $...$, $N$); consequently, the degree of redundancy of the
MBFA is $N_{d}=2N$, i.e. the same as DBF. The computation of the messages
passed in the $k$-th recursion of the MFBA is based on the same equations as
those derived for DBF; the only modifications are due to the fact that:

1) The measurement update accomplished by the $i-$th particle filter\ of the
MBFA requires integrating out the dependence of the measurement vector $%
\mathbf{y}_{k}$ on the $(N-1)$ positions $\{\mathbf{p}_{k}^{(j)}$; $j\neq
i\} $. This marginalization is accomplished by exploiting the pdfs of the
positions $\{\mathbf{p}_{k}^{(i)}$; $i\neq n\}$ \emph{predicted} by the
other $(N-1)$ particle filters. Moreover, the computation of particle
weights requires drawing $L$ particles from the predicted pdfs of the other
filters (see \cite[eq. (7), p. 354]{Djuric_2013}).

2) The computation of the pseudo-measurements for the extended Kalman filter
requires a particle representation for the whole vector $\mathbf{p}_{k}$,
that results from the ordered concatenation of the vectors $\{\mathbf{p}%
_{k}^{(i)}$; $i=1$, $2$, $...$, $N\}$. In the MBFA, the $j-$th particle for $%
\mathbf{p}_{k}$ is generated by: a) taking the $j-$th element of the
particle set made available, after resampling, by each of the $N$ particle
filters (with $j=1$, $2$, $...$, $\tilde{N}_{p}$); b) concatenating the $N$
particles obtained in this way.

In our computer simulations, $N_{p}=500$ has been selected for RBPF, DBF and
SDBF. Moreover, in the MPF technique and in the MBFA, each of $N$ particle
filters makes use of $\tilde{N}_{p}=\left\lfloor N_{p}/N\right\rfloor $
particles, where $N_{p}=500$. Note also\ that: a) the parameter $\tilde{N}%
_{p}$ corresponds to the parameter $M$ of \cite[Sec. III]{Djuric_2013},
since $J=1$ is set in MPF (where $J$ denotes the number of children
generated in the time update step); b) in our simulations, the ratio $L/%
\tilde{N}_{p}$ is always close to $1/3$ for both the MPF technique and the
MBFA. These choices ensure that all the algorithms involving PF have
comparable execution times; for instance, the execution time of RBPF, MPF
and DBF is approximately $21.4$ \%, $3.4$ \% and $0.9$ \% larger,
respectively, than that of the MBFA for $N=5$ targets. Despite this, these
techniques exhibit different behaviours. In fact, our computer simulations
have evidenced that, on the one hand, EKF\ and SDBF\ quickly diverge after
their initialization and, therefore, are useless in the considered scenario.
On the other hand, the RBPF, the MPF and the DBF techniques, and the MBFA
achieve similar accuracies in tracking conditions, but are characterized by
different probabilities of divergence. This is evidenced by Fig. \ref%
{Fig_3_sim}, that shows the dependence of the probability $P_{FD}$ on the
overall number of targets. From these results it is easily inferred that, as
the number of target increases, the RBPF and the MPF techniques are
substantially outperformed by the DBF technique and the MBFA. These results
lead to the conclusion that the property of redundancy can play a key role
in some applications, since it can substantially reduce the probability of
divergence of a filtering algorithm.

\begin{figure}[tbp]
\centering\includegraphics[width=0.55\textwidth]{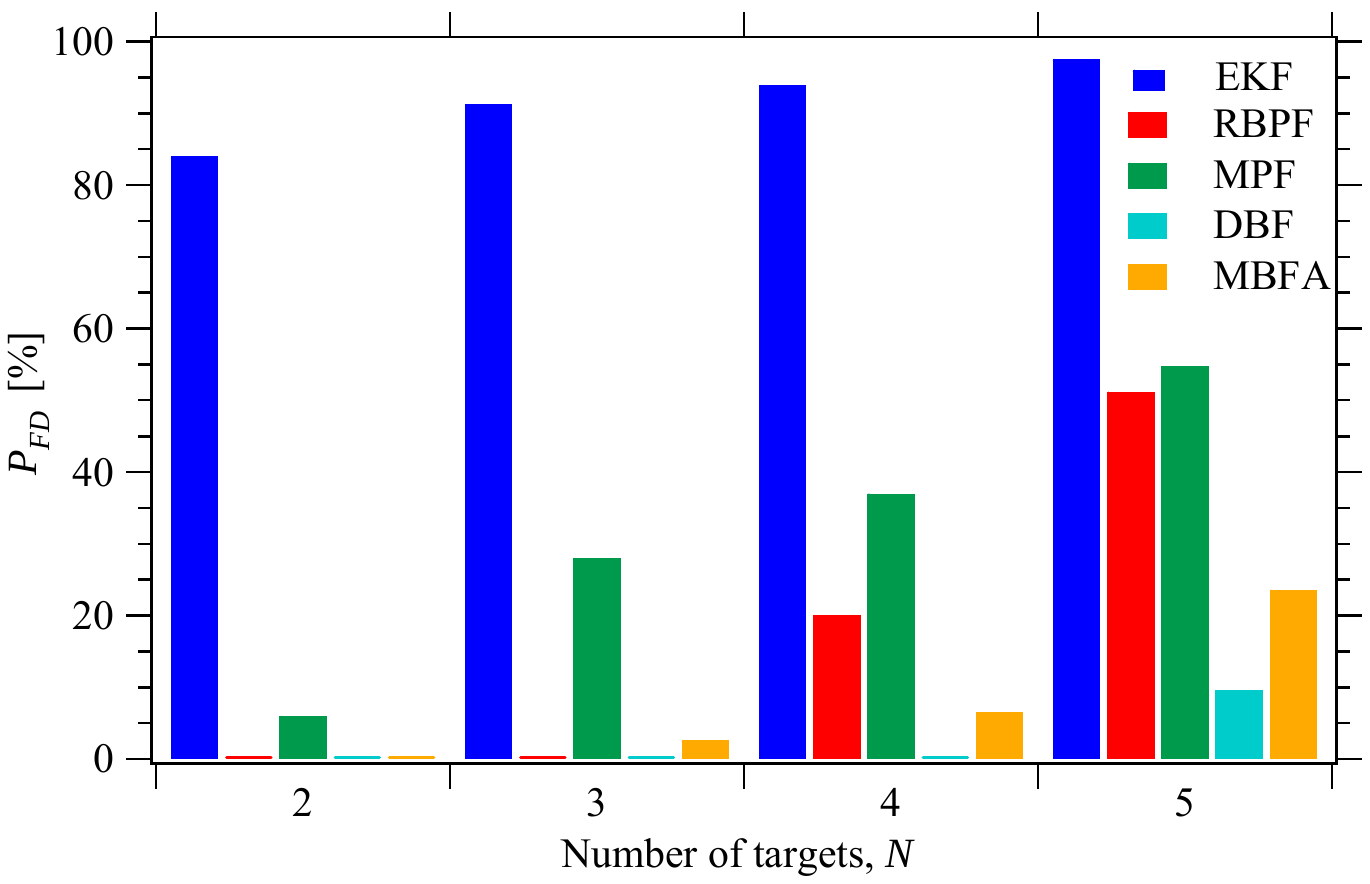}
\caption{Probability of divergence versus $N$ for the RBPF, the MPF and the
DBF techniques, and the MBFA.}
\label{Fig_3_sim}
\end{figure}

%\vspace{-0.3cm}

\section{Conclusions\label{sec:conc}}

In this manuscript, the problem of developing filtering algorithms that
involve multiple interconnected Bayesian filters running in parallel has
been investigated. The devised solution, called \emph{multiple Bayesian
filtering}, is based on the factor graph representation of Bayesian
filtering. The application of our graphical approach to a network consisting
of two Bayesian filters has been illustrated. Moreover, a specific instance
of the proposed approach has been analysed in detail for the case in which
the considered SSM is CLG, and the interconnected filters are an extended
Kalman filter and a particle filter. Simulation results for two specific
SSMs evidence that the devised filtering techniques perform closely to other
well known filtering methods, but are appreciably faster or offer a better
tracking capability.

\appendices

\section{\label{app:A}}

In this Appendix, the derivation of the expressions of various messages
evaluated in each of the three phases which the DBF technique consists of is
sketched.

\textbf{Phase I} - Message $\vec{m}_{1}\left( \mathbf{x}_{k}\right) $ (\ref%
{m_ms_EKF_bis}) conveys the pdf $\tilde{f}(\mathbf{y}_{k}|\mathbf{x}_{k})$ (%
\ref{Measurement_ekf_approx}); therefore, it can be expressed as $\vec{m}%
_{1}\left( \mathbf{x}_{k}\right) =\mathcal{N}\left( \mathbf{y}_{k};\mathbf{H}%
_{k}^{T}\,\mathbf{x}_{k}+\mathbf{v}_{k},\mathbf{C}_{e}\right) $. The last
formula can be easily put in the equivalent form (\ref{m_ms_EKF_bis}) (see 
\cite[Table 3, p. 1304, Eqs. (III.5) and (III.6)]{Loeliger_2007}). Then,
substituting eqs. (\ref{eq:eq:message_fp_l_EKF}) and (\ref{m_ms_EKF_bis}) in
the RHS of Eq. (\ref{m_fe1_EKF_bis}) and applying formula no. 2 of \cite[%
Table I]{Vitetta_2019} yields Eqs. (\ref{m_fe1_EKF_tris})--(\ref{w_fe1_EKF}).

\textbf{Phase II} - Step 1) The derivation of the formulas (\ref%
{m_fe2_EKF_bis}), (\ref{W_fe2_EKF_k}) and (\ref{w_fe2_EKF_k}) referring to
the message $\vec{m}_{3}^{(n)}\left( \mathbf{x}_{k}\right) $ can be
considered as a straightforward application of formula no. 2 of \cite[Table I%
]{Vitetta_2019}, since Eq. (\ref{m_fe2_EKF}) has exactly the same structure
as Eq. (\ref{CR_1}), and both $\vec{m}_{2}(\mathbf{x}_{k})$ and $\,\vec{m}%
_{4}^{(n-1)}(\mathbf{x}_{k})$ are Gaussian messages.$\,$

Step 2) - The expression (\ref{m_ms_j_k}) of the weight $w_{1,k,j}^{(n)}$
can be derived as follows. We first substitute Eq. (\ref{f_y_N})
(conditioned on $\mathbf{x}_{k}^{(N)}=\mathbf{x}_{k,j}^{(N)}[n]$) and Eq. (%
\ref{m_fe_L_EKF_2}) in the RHS\ of Eq. (\ref{m_1_x_N}); then, the resulting
integral is solved by applying formula no. 1 of \cite[Table II]{Vitetta_2019}%
.

Step 3) - The derivation of the expression (\ref{m_pm_x_N_l_ja}) for the
weight $w_{3,k,j}^{(n)}$ is similar to that illustrated for the particle
weights originating from the pseudo-measurements in \emph{dual }RBPF and can
be summarised as follows (additional mathematical details can be found in 
\cite[Sec. V, pp. 1528-1529]{Vitetta_2019}). Two different Gaussian
densities are derived for the random vector $\mathbf{z}_{k}^{(N)}$ (\ref%
{z_N_l}),\emph{\ conditioned on} $\mathbf{x}_{k}^{(N)}$. The expression of
the first density originates from the definition (\ref{z_N_l}) and from the
knowledge of the \emph{joint} pdf of $\mathbf{x}_{k}^{(L)}$ and $\mathbf{x}%
_{k+1}^{(L)}$; this joint density is obtained from: a) the statistical
information provided by the messages $\vec{m}_{2}(\mathbf{x}_{k}^{(L)})=%
\mathcal{N}(\mathbf{x}_{k}^{(L)};\mathbf{\tilde{\eta}}_{2,k},\mathbf{\tilde{C%
}}_{2,k})$ and $\vec{m}_{3}^{(n)}(\mathbf{x}_{k}^{(L)})=\mathcal{N}(\mathbf{x%
}_{k}^{(L)};\mathbf{\tilde{\eta}}_{3,k}^{(n)},\mathbf{\tilde{C}}%
_{3,k}^{(n)}) $, resulting from the marginalization of $\vec{m}_{2}(\mathbf{x%
}_{k})$ (\ref{m_fe1_EKF_tris}) and $\vec{m}_{3}^{(n)}\left( \mathbf{x}%
_{k}\right) $ (\ref{m_fe2_EKF_bis}), respectively, with respect to $\mathbf{x%
}_{k}^{(N)}$; b) the Markov model $f(\mathbf{x}_{k+1}^{(L)}|\mathbf{x}%
_{k}^{(N)},\mathbf{x}_{k}^{(L)})$ (\ref{f_x_L}). This leads to the pdf 
\begin{equation}
f_{1}^{(n)}\left( \mathbf{z}_{k}^{(N)}\left\vert \mathbf{x}_{k}^{(N)}\right.
\right) =\mathcal{N}\left( \mathbf{z}_{k}^{(N)};\mathbf{\check{\eta}}%
_{z,k}^{(n)}\left( \mathbf{x}_{k}^{(N)}\right) ,\mathbf{\check{C}}%
_{z,k}^{(n)}\left( \mathbf{x}_{k}^{(N)}\right) \right) ,  \label{f_z_N}
\end{equation}%
where 
\begin{equation}
\mathbf{\check{\eta}}_{z,k}^{(n)}\left( \mathbf{x}_{k}^{(N)}\right) =\mathbf{%
A}_{k}^{(L)}\left( \mathbf{x}_{k}^{(N)}\right) \left[ \mathbf{\tilde{\eta}}%
_{3,k}^{(n)}-\mathbf{\tilde{\eta}}_{2,k}\right] +\mathbf{f}_{k}^{(L)}\left( 
\mathbf{x}_{k}^{(N)}\right) ,  \label{eta_mess_z_N}
\end{equation}%
and%
\begin{eqnarray}
\mathbf{\check{C}}_{z,k}^{(n)}\left( \mathbf{x}_{k}^{(N)}\right) =\mathbf{C}%
_{w}^{(L)}+\mathbf{A}_{k}^{(L)}\left( \mathbf{x}_{k}^{(N)}\right) \left[ 
\mathbf{\tilde{C}}_{3,k}^{(n)}-\mathbf{\tilde{C}}_{2,k}\right] \left( 
\mathbf{A}_{k}^{(L)}\left( \mathbf{x}_{k}^{(N)}\right) \right) ^{T}.
\label{C_mess_Z_N_bis}
\end{eqnarray}%
The second pdf of $\mathbf{z}_{k}^{(N)}$, instead, results from the fact
that this vector $\mathbf{z}_{k}^{(N)}$ must equal the sum (\ref{z_N_l_bis}%
); consequently, it is given by%
\begin{equation}
f_{2}\left( \mathbf{z}_{k}^{(N)}\left\vert \mathbf{x}_{k}^{(N)}\right.
\right) =\mathcal{N}\left( \mathbf{z}_{k}^{(N)};\mathbf{f}_{k}^{(L)}\left( 
\mathbf{x}_{k}^{(N)}\right) ,\mathbf{C}_{w}^{(N)}\right) .  \label{fZ_N_bis}
\end{equation}%
Given the pdfs (\ref{f_z_N}) and (\ref{fZ_N_bis}), the message $\vec{m}%
_{3}^{(n)}(\mathbf{x}_{k}^{(N)})$ is expressed by their \emph{correlation},
i.e. it is computed as 
\begin{equation}
\vec{m}_{3}^{(n)}(\mathbf{x}_{k}^{(N)}) = \int f_{1}^{(n)}\left( \mathbf{z}%
_{k}^{(N)}\left\vert \mathbf{x}_{k}^{(N)}\right. \right) \cdot f_{2}\left( 
\mathbf{z}_{k}^{(N)}\left\vert \mathbf{x}_{k}^{(N)}\right. \right) d\mathbf{z%
}_{k}^{(N)}.  \label{m_pm_x_N_l_j_first}
\end{equation}%
Substituting (\ref{f_z_N}) and (\ref{fZ_N_bis}) in the RHS of the last
expression, setting $\mathbf{x}_{k}^{(N)}=\mathbf{x}_{k,j}^{(N)}[n]$ and
applying formula no. 4 of \cite[Table II]{Vitetta_2019} to the evaluation of
the resulting integral yields Eq. (\ref{m_pm_x_N_l_ja}); note that $\mathbf{%
\check{\eta}}_{z,k,j}^{(n)}$ (\ref{eta_pm_x_N_l_j}) and $\mathbf{\check{C}}%
_{z,k,j}^{(n)}$ (\ref{Cov_pm_x_N_l_j}) represent the values taken on by $%
\mathbf{\check{\eta}}_{z,k}^{(n)}(\mathbf{x}_{k}^{(N)})$ (\ref{eta_mess_z_N}%
) and $\mathbf{\check{C}}_{z,k}^{(n)}(\mathbf{x}_{k}^{(N)})$ (\ref%
{C_mess_Z_N_bis}), respectively, for $\mathbf{x}_{k}^{(N)}=\mathbf{x}%
_{k,j}^{(N)}[n]$.

Step 4) - Formula (\ref{m_fe_2_x_N_l_bis}), that refers to the message $\vec{%
m}_{4,j}^{(n)}(\mathbf{x}_{k}^{(N)})$, is obtained by substituting $\vec{m}%
_{2,j}(\mathbf{x}_{k}^{(N)})$ (\ref{eq:message_fe_l_j_PF_bis}) in the RHS of
Eq. (\ref{m_fe2_k_PF}) and observing that $w_{3,k,j}^{(n)}$ (\ref%
{m_pm_x_N_l_ja}) represents the value taken on by the message $\vec{m}%
_{3}^{(n)}(\mathbf{x}_{k}^{(N)})$ for $\mathbf{x}_{k}^{(N)}=\mathbf{x}%
_{k,j}^{(N)}[n]$.

Step 5) Eq. (\ref{eq:message_5_Nb}) is results from substituting Eqs. (\ref%
{m_fe_L_EKF_2}) and (\ref{m_pm_x_N_l_j_first}) in Eq. (\ref{eq:message_5_Na}%
) and, then, applying formula no. 1 of \cite[Table III]{Vitetta_2019} to
evaluate the resulting integral.

Step 6) - The message $\vec{m}_{4}^{(n)}\left( \mathbf{x}_{k}\right) $ (\ref%
{message_pm_x_l_k}) results from merging, in the F$_{2}{\rightarrow }$F$_{1}$
block, the statistical information about the \emph{nonlinear} state
component conveyed by the message $\vec{m}_{4}^{(n)}(\mathbf{x}_{k}^{(N)})$
(and, consequently, by its components $\{\vec{m}_{4,j}^{(n)}(\mathbf{x}%
_{k}^{(N)})\}$; see Eq. (\ref{m_fe_2_x_N_l_tris})) with those provided by
the pseudo-measurement $\mathbf{z}_{k}^{(L)}$ (\ref{eq:z_L_l}) about the 
\emph{linear} state component. The method employed for processing this
pseudo-measurement is the same as that developed for RBPF and can be
summarised as follows (additional mathematical details can be found in \cite[%
Sec. IV, p. 1527]{Vitetta_2019}):

a) The particles $\mathbf{x}_{k,j}^{(N)}[n+1]$ and $\mathbf{\bar{x}}%
_{k+1,j}^{(N)}[n+1]$, conveyed by the messages $\vec{m}_{4,j}^{(n)}(\mathbf{x%
}_{k}^{(N)})$ (\ref{m_fe_2_x_N_l_tris}) and $\vec{m}_{\mathrm{fp},j}^{(n)}(%
\mathbf{x}_{k+1}^{(N)})$ (\ref{eq:message_5_N_l+1}), respectively, are
employed to compute the $j-$th realization $\mathbf{z}_{k,j}^{(L)}\left[ n+1%
\right] $ (\ref{eq:z_L_evaluateda}) of the vector $\mathbf{z}_{k}^{(L)}$ (%
\ref{eq:z_L_l}) according to Eq. (\ref{eq:z_L_evaluateda}).

b) The pseudo-measurement $\mathbf{z}_{k,j}^{(L)}\left[ n+1\right] $ (\ref%
{eq:z_L_evaluateda}) is exploited to generate the (particle-dependent)
message%
\begin{equation}
\vec{m}_{4,j}^{(n)}\left( \mathbf{x}_{k}^{(L)}\right) =\mathcal{\mathcal{N}}%
\left( \mathbf{x}_{k}^{(L)};\mathbf{\tilde{\eta}}_{4,k,j}^{(n)},\mathbf{%
\tilde{C}}_{4,k,j}^{(n)}\right) ,  \label{eq:message_pm_L_j_tris}
\end{equation}%
that conveys pseudo-measurement information about $\mathbf{x}_{k}^{(L)}$;
the covariance matrix $\mathbf{\tilde{C}}_{4,k,j}^{(n)}$ and the mean vector 
$\mathbf{\tilde{\eta}}_{4,k,j}^{(n)}$ of this message are computed on the
basis of the precision matrix $\mathbf{\tilde{W}}_{4,k,j}^{(n)}$ (\ref%
{eq:W_pm_L_j1})\ and the transformed mean vector $\mathbf{\tilde{w}}%
_{4,k,j}^{(n)}$ (\ref{eq:w_pm_L_j}), respectively. Finally, the message $%
\vec{m}_{4}^{(n)}\left( \mathbf{x}_{k}\right) $ (\ref{message_pm_x_l_k})
results from merging the message $\vec{m}_{4}^{(n)}(\mathbf{x}_{k}^{(N)})$
(its $j-$th component is expressed by Eq. (\ref{m_fe_2_x_N_l_tris})) with
the pdfs $\{\vec{m}_{4,j}^{(n)}(\mathbf{x}_{k}^{(L)})\}$ (see Eq. (\ref%
{eq:message_pm_L_j_tris})); the adopted approach is based on the fact that:
a) as it can be easily inferred from our previous derivations, the Gaussian
message $\vec{m}_{4,j}^{(n)}(\mathbf{x}_{k}^{(L)})$ (\ref%
{eq:message_pm_L_j_tris}) is evaluated under the condition that $\mathbf{x}%
_{k}^{(N)}=\mathbf{x}_{k,j}^{(N)}[n+1]$; b) the messages $\vec{m}%
_{4,j}^{(n)}(\mathbf{x}_{k}^{(N)})$ and $\vec{m}_{4,j}^{(n)}(\mathbf{x}%
_{k}^{(L)})$ provide \emph{complementary} information, because they refer to
the two different components of the overall state $\mathbf{x}_{k}$.
Consequently, the statistical information conveyed by the sets $\{\vec{m}%
_{4,j}^{(n)}(\mathbf{x}_{k}^{(N)})\}$ and $\{\vec{m}_{4,j}^{(n)}(\mathbf{x}%
_{k}^{(L)})\}$ can be merged in the joint pdf%
\begin{equation}
f^{(k)}(\mathbf{x}_{k}^{(L)},\mathbf{x}_{k}^{(N)})\triangleq
w_{p}\sum\limits_{j=1}^{N_{p}}\vec{m}_{4,j}^{(n)}\left( \mathbf{x}%
_{k}^{(N)}\right) \vec{m}_{4,k}^{(n)}\left( \mathbf{x}_{k}^{(L)}\right) .
\label{joint_pdf}
\end{equation}%
referring to $\mathbf{x}_{k}$. Then, the message $\vec{m}_{4}^{(n)}\left( 
\mathbf{x}_{k}\right) $ (\ref{message_pm_x_l_k}) is computed by projecting
the pdf $f^{(k)}(\mathbf{x}_{k}^{(L)},\mathbf{x}_{k}^{(N)})$ (\ref{joint_pdf}%
) onto a single Gaussian pdf having the same \emph{mean }and \emph{covariance%
}.

\textbf{Phase III} - The message $\vec{m}_{\mathrm{fp}}\left( \mathbf{x}%
_{k+1}\right) $ (\ref{m_fp_l+1_whole_tris}) is computed as follows.
Substituting the expressions (\ref{Markov_ekf_approx}) of $\tilde{f}\left( 
\mathbf{x}_{k+1}\left\vert \mathbf{x}_{k}\right. \right) $ and (\ref%
{m_fe2_EKF_final}) of $\vec{m}_{3}^{(n_{i}+1)}\left( \mathbf{x}_{k}\right) $
in the RHS of Eq. (\ref{m_fp_l+1_whole}) and applying formula no. 1 of \cite[%
Table II]{Vitetta_2019} to the evaluation of the resulting integral produces
Eqs. (\ref{m_fp_l+1_whole_tris})--(\ref{C_fp_l+1_whole}).

\section{Computational complexity of the DBF and SDBF techniques \label%
{app:CDBFA}}

In this appendix, the computational complexity of the tasks accomplished in
a single recursion of the DBF technique is assessed in terms of flops.
Moreover, we comment on how the illustrated results can be also exploited to
assess the computational complexity of a single recursion of the SDBF
technique. In the following, $\mathcal{C}_{\mathbf{H}}$, $\mathcal{C}_{%
\mathbf{B}}$, $\mathcal{C}_{\mathbf{F}}$, $\mathcal{C}_{\mathbf{A}^{(L)}}$
and $\mathcal{C}_{\mathbf{A}^{(N)}}$, and $\mathcal{C}_{\mathbf{g}}$, $%
\mathcal{C}_{\mathbf{f}^{(L)}}$, $\mathcal{C}_{\mathbf{f}^{(N)}}$ and $%
\mathcal{C}_{\mathbf{f}_{k}}$denote the cost due to the evaluation of the
matrices $\mathbf{H}_{k}$, $\mathbf{B}_{k}$, $\mathbf{F}_{k}$, $\mathbf{A}%
_{k}^{(L)}(\mathbf{x}_{k}^{(N)})$ and $\mathbf{A}_{k}^{(N)}(\mathbf{x}%
_{k}^{(N)})$, and of the functions $\mathbf{g}_{k}(\mathbf{x}_{k}^{(N)})$, $%
\mathbf{f}_{k}^{(L)}(\mathbf{x}_{k}^{(N)})$, $\mathbf{f}_{k}^{(N)}(\mathbf{x}%
_{k}^{(N)})$ and $\mathbf{f}_{k}(\mathbf{x}_{k})$, respectively. Moreover,
similarly as \cite{Hoteit_2016}, it is assumed that the computation of the
inverse of any covariance matrix involves a Cholesky decomposition of the
matrix itself and the inversion of a lower or upper triangular matrix.

1. \textbf{Filter F$_{1}$, first measurement update}

The overall computational cost of this task is (see Eqs. (\ref{w_ms_EKF})-(%
\ref{W_ms_EKF}) and (\ref{W_fe1_EKFa})-(\ref{w_fe1_EKF})) 
\begin{equation}
\mathcal{C}_{MU1}^{(1)}=\mathcal{C}_{\mathbf{W}_{2,k}}+\mathcal{C}_{\mathbf{w%
}_{2,k}}+\mathcal{C}_{\mathbf{C}_{2,k}}+\mathcal{C}_{\eta _{2,k}}
\label{C_F1_MU1}
\end{equation}%
Moreover, we have that: 1) the cost $\mathcal{C}_{\mathbf{W}_{2,k}}$ is
equal to $\mathcal{C}_{\mathbf{H}}+2P{D}^{2}+2P^{2}D-PD$ flops; 2) the cost $%
\mathcal{C}_{\mathbf{w}_{2,k}}$ is equal to $\mathcal{C}_{\mathbf{B}}+%
\mathcal{C}_{\mathbf{g}}+2P^{2}D+5PD_{L}+3PD_{N}-P$ flops ($\mathbf{H}_{k}$
has been already computed at point 1); 3) the cost $\mathcal{C}_{\mathbf{C}%
_{2,k}}$ is equal to $2D^{3}/3+3D^{2}/2+5D/6$ flops; 4) the cost $\mathcal{C}%
_{\eta _{2,k}}$ is equal to $D(2D-1)$ flops. The expressions listed at
points 1)-4) can be exploited for the SDBF too; in the last case, however, $%
D_{N}=0$ and $D=D_{L}$ must be assumed.

2. \textbf{Filter F$_{2}$, second measurement update}

The overall computational cost of this task is (see Eqs. (\ref{W_fe2_EKF_k}%
)-(\ref{w_fe2_EKF_k})) 
\begin{equation}
\mathcal{C}_{MU2}^{(1)}=n_{i}\left( \mathcal{C}_{\mathbf{C}_{3,k}^{(n)}}+%
\mathcal{C}_{\eta _{3,k}^{(n)}}\right) ,  \label{C_F1_MU2}
\end{equation}%
where the costs $\mathcal{C}_{\mathbf{C}_{3,k}^{(n)}}$ and $\mathcal{C}%
_{\eta _{3,k}^{(n)}}$ are equal to $D^{2}(2D-1)$ flops and $4D^{2}-D$ flops,
respectively; if the SDBF is considered, we have that $D=D_{L}$ in the last
two expressions.

3. \textbf{Filter F$_{2}$, first measurement update}

The overall computational cost of this task is (see Eqs. (\ref{m_ms_j_k})-(%
\ref{eq:C_fe_PF})) 
\begin{equation}
\mathcal{C}_{MU1}^{(2)}=n_{i}\,N_{p}\left( \mathcal{C}_{\tilde{\eta}%
_{1,k,j}^{(n)}}+\mathcal{C}_{\mathbf{\tilde{C}}_{1,k,j}^{(n)}}+\mathcal{C}%
_{w_{1,k,j}^{(n)}}\right) .  \label{C_F2_MU1}
\end{equation}%
Moreover, we have that: 1) the cost $\mathcal{C}_{\tilde{\eta}%
_{1,k,j}^{(n)}} $ is equal to $\mathcal{C}_{\mathbf{B}}+\mathcal{C}_{\mathbf{%
g}}+2PD_{L}$ flops; 2) the cost $\mathcal{C}_{\mathbf{\tilde{C}}%
_{1,k,j}^{(n)}}$ is equal to $2PD_{L}^{2}+2P^{2}D_{L}-PD_{L}$ flops (the
cost for computing $\mathcal{C}_{\mathbf{B}}$ has been already accounted for
at point 1)); 3) the cost $\mathcal{C}_{w_{1,k,j}^{(n)}}$ is equal to $%
(4P^{3}+21P^{2}+17P+6)/6$ flops.

4. \textbf{Filter F$_{2}$, second measurement update}

The overall computational cost of this task is (see Eqs. (\ref{w_fe_2_x_N_l}%
) and (\ref{W_fe_2_x_N_l})) 
\begin{equation}
\mathcal{C}_{MU2}^{(2)}=\mathcal{C}_{w_{4,k,j}^{(n)}}+\mathcal{C}%
_{W_{4,k,j}^{(n)}}+n_{i}\,\mathcal{C}_{R}(N_{p}),  \label{C_F2_MU2}
\end{equation}%
where the costs $\mathcal{C}_{w_{4,k,j}^{(n)}}$ and $\mathcal{C}%
_{W_{4,k,j}^{(n)}}$ are equal to $n_{i}N_{p}$ flops and $2N_{p}-1$ flops,
respectively, and $\mathcal{C}_{R}(N_{p})$ denotes the total cost of the
resampling step (that involves a particle set of size $N_{p}$).

5. \textbf{Computation of the pseudo-measurements for filter F$_{2}$}

The overall computational cost of this task is (see Eqs. (\ref{m_pm_x_N_l_ja}%
)-(\ref{eta_pm_x_N_l_j})) 
\begin{align}
\mathcal{C}_{1\rightarrow 2}=& n_{i}\,N_{p}\left( \mathcal{C}_{\check{\eta}%
_{z,k}^{(n)}}+\mathcal{C}_{\mathbf{\check{C}}_{z,k}^{(n)}}+\mathcal{C}_{%
\mathbf{\check{W}}_{z,k}^{(n)}}+\mathcal{C}_{\mathbf{\check{w}}%
_{z,k}^{(n)}}+\right.  \notag \\
& \left. \mathcal{C}_{\mathbf{\check{W}}_{3,k,j}^{(n)}}+\mathcal{C}_{\mathbf{%
\check{w}}_{3,k,j}^{(n)}}+\mathcal{C}_{\mathbf{\check{C}}_{3,k,j}^{(n)}}+%
\mathcal{C}_{\mathbf{\eta }_{3,k,j}^{(n)}}+\mathcal{C}_{w_{3,k,j}^{(n)}}%
\right) .  \label{C_F2_PM}
\end{align}%
Moreover, we have that: 1) the cost $\mathcal{C}_{\check{\eta}_{z,k}^{(n)}}$
is equal to $\mathcal{C}_{\mathbf{A}^{(L)}}+\mathcal{C}_{\mathbf{f}%
^{(L)}}+2D_{L}^{2}+D_{L}$ flops; 2) the cost $\mathcal{C}_{\mathbf{\check{C}}%
_{z,k}^{(n)}}$ is equal $4D_{L}^{3}$ flops (since the cost for computing $%
\mathcal{C}_{\mathbf{A}^{(L)}}$ has been already accounted for at point 1);
3) the cost $\mathcal{C}_{\mathbf{\check{W}}_{z,k}^{(n)}}$ is equal to $%
2D_{L}^{3}/3+3D_{L}^{2}/2+5D_{L}/6$ flops; 4) the cost $\mathcal{C}_{\mathbf{%
\check{w}}_{z,k}^{(n)}}$ is equal to $D_{L}(2D_{L}-1)$ flops; 5) the cost $%
\mathcal{C}_{\mathbf{\check{W}}_{3,k,j}^{(n)}}$ is equal to $D_{L}^{2}$
flops; 6) the cost $\mathcal{C}_{\mathbf{\check{w}}_{3,k,j}^{(n)}}$ is equal
to $2D_{L}^{2}$ flops (the cost for computing $\mathcal{C}_{\mathbf{f}%
^{(L)}} $ has been already accounted for at point 1); 7) the cost $\mathcal{C%
}_{\mathbf{\check{C}}_{3,k,j}^{(n)}}$ is equal to $%
2D_{L}^{3}/3+3D_{L}^{2}/2+5D_{L}/6$ flops; 8) the cost $\mathcal{C}_{\mathbf{%
\eta }_{3,k,j}^{(n)}}$ is equal to $D_{L}(2D_{L}-1)$ flops; 9) the cost $%
\mathcal{C}_{w_{3,k,j}^{(n)}}$ is equal to $6D_{L}^{2}+3D_{L}+1$ flops (the
cost for computing $\mathcal{C}_{\mathbf{f}^{(L)}}$ has been already
accounted for at point 1)).

6. \textbf{Computation of the pseudo-measurements for filter F$_{1}$}

The overall computational cost of this task is (see Eqs. (\ref{eta_pm_l_L_k}%
)-(\ref{eq:z_L_evaluateda})) 
\begin{align}
\mathcal{C}_{2\rightarrow 1}=& n_{i}N_{p}\left( \mathcal{C}_{\mathbf{z}%
_{k,j}^{(L)}}+\mathcal{C}_{\mathbf{\tilde{W}}_{4,k,j}^{(n)}}+\mathcal{C}_{%
\mathbf{\tilde{w}}_{4,k,j}^{(n)}}+\mathcal{C}_{\mathbf{\tilde{C}}%
_{4,k,j}^{(n)}}+\mathcal{C}_{\tilde{\eta}_{4,k,j}^{(n)}}\right)  \notag \\
& +\mathcal{C}_{\mathbf{C}_{4,k}^{(n)}}+\mathcal{C}_{\eta _{4,k}^{(n)}}+%
\mathcal{C}_{\mathbf{W}_{k}^{(n)}}.  \label{C_F1_PM}
\end{align}%
Moreover, we have that: 1) the cost $\mathcal{C}_{\mathbf{z}_{k,j}^{(L)}}$
is equal to $D_{N}$ flops (the cost for computing $\mathcal{C}_{\mathbf{f}%
^{(N)}}$ has been already accounted for in the time update of filter F$_{2}$%
); 2) the cost $\mathcal{C}_{\mathbf{\tilde{W}}_{4,k,j}^{(n)}}$ is equal to $%
DD_{L}(2D_{N}-1)$ flops (the cost for computing $\mathcal{C}_{\mathbf{A}%
^{(N)}}$ has been already accounted for in the time update of filter F$_{2}$%
); 3) the cost $\mathcal{C}_{\mathbf{\tilde{w}}_{4,k,j}^{(n)}}$ is equal to $%
D_{L}(2D_{N}^{2}+D_{N}-1)$ flops (the cost for computing $\mathcal{C}_{%
\mathbf{A}^{(N)}}$ has been already accounted for in the time update of
filter F$_{2}$); 4) the cost $\mathcal{C}_{\mathbf{\tilde{C}}_{4,k,j}^{(n)}}$
is equal to $2D_{L}^{3}/3+3D_{L}^{2}/2+5D_{L}/6$ flops; 5) the cost $%
\mathcal{C}_{\tilde{\eta}_{4,k,j}^{(n)}}$ is equal to $D_{L}(2D_{L}-1)$
flops; 6) the cost $\mathcal{C}_{\mathbf{C}_{4,k}^{(n)}}$ is equal to $%
n_{i}(2N_{p}D_{L}^{2}+N_{p}D_{N}^{2}+2D_{L}^{2}+2D_{N}^{2}+N_{p}D_{L}D_{N}+2D_{L}D_{N}+3N_{p}) 
$ flops; 7) the cost $\mathcal{C}_{\eta _{4,k}^{(n)}}$ is equal to $%
n_{i}(D(N_{p}-1)+1)$ flops; 8) the cost $\mathcal{C}_{\mathbf{W}_{k}^{(n)}}$
is $n_{i}(16D^{3}+9D^{2}+5D)/6$ flops.\newline
If the SDBF is considered, the total costs 1)-5) remain unchanged, whereas $%
D_{N}=0$ and $D=D_{L}$ in the costs $\mathcal{C}_{\eta _{4,k}^{(n)}}$ and $%
\mathcal{C}_{\mathbf{W}_{k}^{(n)}}$ (see points 7) and 8)); moreover, the
cost $\mathcal{C}_{\mathbf{C}_{4,k}^{(n)}}$ becomes $%
n_{i}(2N_{p}D_{L}^{2}+2D_{L}^{2}+N_{p})$ flops (see point 6)).

7. \textbf{Filter F$_{1}$, time update}

The overall computational cost of this task is (see Eqs. (\ref%
{eta_fp_l+1_whole1})-(\ref{C_fp_l+1_whole})) 
\begin{equation}
\mathcal{C}_{TU}^{(1)}=\mathcal{C}_{\eta _{\mathrm{fp},k+1}}+\mathcal{C}_{%
\mathbf{C}_{\mathrm{fp},k+1}}+\mathcal{C}_{\mathbf{W}_{\mathrm{fp},k+1}}+%
\mathcal{C}_{\mathbf{w}_{\mathrm{fp},k+1}},  \label{C_F1_TU}
\end{equation}%
since $\mathcal{C}_{\mathbf{C}_{3,k}^{(n_{i}+1)}}$ and $\mathcal{C}_{\eta
_{3,k}^{(n_{i}+1)}}$ have been already computed in the previous time update
of filter F$_{2}$. Moreover, we have that:\newline
1) $\mathcal{C}_{\eta _{\mathrm{fp},k+1}}$ is equal to $\mathcal{C}_{\mathbf{%
f}_{k}}$ flops; 2) $\mathcal{C}_{\mathbf{C}_{\mathrm{fp},k+1}}$ is equal to $%
\mathcal{C}_{\mathbf{F}}+D^{2}(4D-1)$ flops; 3) $\mathcal{C}_{\mathbf{W}_{%
\mathrm{fp},k+1}}$ is equal to $2D^{3}/3+3D^{2}/2+5D/6$ flops; 4) $\mathcal{C%
}_{\mathbf{w}_{\mathrm{fp},k+1}}$ is equal to $D(2D-1)$ flops. If the SDBF
is considered, $D=D_{L}$ is set in the expressions of the costs listed at
points 1)-4).

8. \textbf{Filter F$_{2}$, time update}

The overall computational cost of this task is (see Eqs. (\ref{eq:eta_5_N})-(%
\ref{eq:message_5_N_l+1})) 
\begin{equation}
\mathcal{C}_{TU}^{(2)}=n_{i}\,N_{p}\left( \mathcal{C}_{\eta _{3,k,j}^{(N)}}+%
\mathcal{C}_{\mathbf{C}_{3,k,j}^{(N)}}+\mathcal{C}_{\mathbf{x}%
_{k+1,j}^{(N)}}\right) .  \label{C_F2_TU}
\end{equation}%
Moreover, we have that: 1) the cost $\mathcal{C}_{\eta _{3,k,j}^{(N)}}$ is
equal to $\mathcal{C}_{\mathbf{A}^{(N)}}+\mathcal{C}_{\mathbf{f}%
^{(N)}}+2D_{L}D_{N}$ flops; 2) the cost $\mathcal{C}_{\mathbf{C}%
_{3,k,j}^{(N)}}$ is equal to $D_{L}D_{N}(2D-1)$ flops ($\mathcal{C}_{\mathbf{%
A}^{(N)}}$ has been already accounted for at point 1)); 3) the cost $%
\mathcal{C}_{\mathbf{x}_{k+1,j}^{(N)}}$ is equal to $%
D_{N}^{3}/3+3D_{N}^{2}+5D_{N}/3$ flops.

\section{Computational complexity of the EKF technique\label{app:CEKF}}

In this appendix analysis of EKF complexity is illustrated; the notation is
the same as \cite[pp. 194-195]{Anderson_1979}. In the following, $\mathcal{C}%
_{\mathbf{H}}$ and $\mathcal{C}_{\mathbf{F}}$, $\mathcal{C}_{\mathbf{h}_{k}}$
and $\mathcal{C}_{\mathbf{f}_{k}}$ denote the cost due to the evaluation of
the matrices $\mathbf{H}_{k}$ and $\mathbf{F}_{k}$, and of the functions $%
\mathbf{h}_{k}(\mathbf{x}_{k})$ and $\mathbf{f}_{k}(\mathbf{x}_{k})$,
respectively. Moreover, similarly as \cite{Hoteit_2016}, it is assumed that
the computation of the inverse of any covariance matrix involves a Cholesky
decomposition of the matrix itself and the inversion of a lower or upper
triangular matrix.

1. \textbf{Measurement update}

The overall computational cost of this task is 
\begin{equation}
\mathcal{C}_{MU}=\mathcal{C}_{\mathbf{\Omega}_{k}}+\mathcal{C}_{\mathbf{L}%
_{k}}+\mathcal{C}_{\eta_{k|k}} + \mathcal{C}_{\mathbf{C}_{k|k}}.
\end{equation}
Moreover, we have that: 1) the cost $\mathcal{C}_{\mathbf{\Omega}_{k}}$ is
equal to $\mathcal{C}_{\mathbf{H}} + 2P^2D+2PD^2-PD $ flops; 2) $\mathcal{C}%
_{\mathbf{L}_{k}}$ is equal to $2P^3/3 + 3P^2/2 + 5P/6 + 2PD^2 + 2P^2D - 2PD 
$ flops; 3) $\mathcal{C}_{\eta_{k|k}}$ is equal to $\mathcal{C}_{\mathbf{h}%
_k} + 2PD+P $ flops; 4) $\mathcal{C}_{\mathbf{C}_{k|k}}$ is equal $2D^3
+2PD^2- D^2 $ flops.

2. \textbf{Time update}

The overall computational cost of this task is 
\begin{equation}
\mathcal{C}_{TU}=\mathcal{C}_{\eta_{k+1|k}} + \mathcal{C}_{\mathbf{C}%
_{k+1|k}},
\end{equation}
where the costs $\mathcal{C}_{\eta_{k+1|k}}$ and $\mathcal{C}_{\mathbf{C}%
_{k+1|k}}$ are equal to $\mathcal{C}_{\mathbf{f}_k} $ flops and $\mathcal{C}%
_{\mathbf{F}} + 4D^3 - D^2 $ flops, respectively.

\section{Computational complexity of the RBPF technique\label{app:CRBPF}}

In this appendix a detailed analysis of the RBPF complexity is provided; the
adopted notation is the same as \cite{Vitetta_2019}. In the following, $%
\mathcal{C}_{ \mathbf{B}}$, $\mathcal{C}_{\mathbf{A}^{(L)}}$ and $\mathcal{C}%
_{\mathbf{A}^{(N)}}$, and $\mathcal{C}_{\mathbf{g}}$, $\mathcal{C}_{\mathbf{f%
}^{(L)}}$ and $\mathcal{C}_{\mathbf{f}^{(N)}}$ denote the cost due to the
evaluation of the matrices $\mathbf{B}_{k}$, $\mathbf{A}_{k}^{(L)}(\mathbf{x}%
_{k}^{(N)})$ and $\mathbf{A}_{k}^{(N)}(\mathbf{x}_{k}^{(N)})$, and of the
functions $\mathbf{g}_{k}(\mathbf{x}_{k}^{(N)})$, $\mathbf{f}_{k}^{(L)}(%
\mathbf{x}_{k}^{(N)})$ and $\mathbf{f}_{k}^{(N)}(\mathbf{x}_{k}^{(N)})$,
respectively. Moreover, similarly as \cite{Hoteit_2016}, it is assumed that
the computation of the inverse of any covariance matrix involves a Cholesky
decomposition of the matrix itself and the inversion of a lower or upper
triangular matrix.

1. \textbf{Measurement update nonlinear part}

The overall computational cost of this task is 
\begin{equation}
\mathcal{C}_{MU}^{(N)}=N_{p}\left(\mathcal{C}_{{\eta}_{1,k,j}^{(N)}}+%
\mathcal{C}_{\mathbf{C}_{1,k,j}^{(N)}}+\mathcal{C}_{w_{fe,k,j}}\right)+%
\mathcal{C}_{W_{fe,k,j}} + \mathcal{C}_{R}(N_{p}).
\end{equation}
Moreover, we have that: 1) the cost $\mathcal{C}_{{\eta}_{1,k,j}^{(N)}}$ is
equal to $\mathcal{C}_{\mathbf{B}}+\mathcal{C}_{\mathbf{g}}+2PD_{L}$ flops;
2) $\mathcal{C}_{{C}_{1,k,j}^{(N)}}$ is equal to $2PD_{L}^2+2P^2D_{L}-PD_{L}$
flops ($\mathcal{C}_{\mathbf{B}}$ has been already accounted for at point
1)); 3) $\mathcal{C}_{w_{fe,k,j}}$ is equal to $(4P^3+21P^2+17P+6)/6$ flops;
4) $\mathcal{C}_{W_{fe,k,j}}$ is equal to $2N_{p}-1$ flops; 5) $\mathcal{C}%
_{R}(N_{p})$ denotes the total cost of the resampling step (that involves a
particle set of size $N_{p}$).

2. \textbf{First measurement update linear part}

The overall computational cost of this task is 
\begin{equation}
\mathcal{C}_{MU1}^{(L)}=N_{p}\left(\mathcal{C}_{\mathbf{w}_{1,k,j}^{(L)}}+%
\mathcal{C}_{\mathbf{W}_{1,k,j}^{(L)}}+\mathcal{C}_{\mathbf{C}%
_{2,k,j}^{(L)}}+\mathcal{C}_{{\eta}_{2,k,j}^{(L)}}\right).
\end{equation}
Moreover, we have that: 1) the cost $\mathcal{C}_{\mathbf{w}_{1,k,j}^{(L)}}$
is equal to $\mathcal{C}_{\mathbf{B}}+\mathcal{C}_{\mathbf{g}%
}+2P^2D_{L}+2PD_{L}-PD_{L}-D_{L}+P$ flops; 2) $\mathcal{C}_{\mathbf{W}%
_{1,k,j}^{(L)}}$ is equal to $2PD_{L}^2+2P^{2}D_{L}-D_{L}^2-PD_{L}$ flops;
3) $\mathcal{C}_{\mathbf{C}_{2,k,j}^{(L)}}$ is equal to $4D_{L}^3/3 +
4D_{L}^2+5D_{L}/3$ flops; 4) $\mathcal{C}_{{\eta}_{2,k,j}^{(L)}}$ is equal
to $D_{L}(4D_{L}-1)$ flops.

3. \textbf{Second measurement update linear part}

The overall computational cost of this task is 
\begin{equation}
\mathcal{C}_{MU2}^{(L)} = N_{p}\left(\mathcal{C}_{\mathbf{z}_{k,j}^{(L)}} + 
\mathcal{C}_{\mathbf{C}_{4,k,j}^{(L)}}+\mathcal{C}_{\eta_{4,k,j}^{(L)}}%
\right).
\end{equation}
Moreover, we have that: 1) the cost $\mathcal{C}_{\mathbf{z}_{k,j}^{(L)}}$
is equal to $\mathcal{C}_{\mathbf{f}^{(N)}} + D_{N} $ flops; 2) $\mathcal{C}%
_{\mathbf{C}_{4,k,j}^{(L)}}$ is equal to $\mathcal{C}_{\mathbf{A}^{(N)}} +
2D_{L}^3/3 + 2D_{L}^2D_{N}+2D_{L}D_{N}^2+3D_{L}^2/2 -D_{L}D_{N} +5D_{L}/6$
flops; 3) $\mathcal{C}_{\eta_{4,k,j}^{(L)}}$ is equal to $2D_{L}D_{N}^2 +
2D_{L}^2 +D_{L}D_{N} -2D_{L} + D_{N}$ flops.

4. \textbf{Time update nonlinear part}

The overall computational cost of this task is 
\begin{equation}
\mathcal{C}_{TU}^{(N)}=N_{p}\left(\mathcal{C}_{\eta_{3,k,j}^{(N)}}+\mathcal{C%
}_{\mathbf{C}_{3,k,j}^{(N)}}+\mathcal{C}_{\mathbf{x}_{\mathrm{fp}%
,k+1,j}^{(N)}}\right).
\end{equation}
Moreover, we have that: 1) the cost $\mathcal{C}_{\eta_{3,k,j}^{(N)}}$ is
equal to $\mathcal{C}_{\mathbf{A}^{(N)}} + \mathcal{C}_{\mathbf{f}^{(N)}} +
2D_{L}D_{N} $ flops; 2) $\mathcal{C}_{\mathbf{C}_{3,k,j}^{(N)}}$ is equal to 
$D_{L}D_{N}(2D-1) $ flops ($\mathcal{C}_{\mathbf{A}^{(N)}}$ has been already
accounted for at point 1)); 3) $\mathcal{C}_{\mathbf{x}_{\mathrm{fp}%
,k+1,j}^{(N)}}$ is equal to $D_{N}^3/3+3D_{N}^2+5D_{N}/3 $ flops.

5. \textbf{Time update linear part}

The overall computational cost of this task is 
\begin{equation}
\mathcal{C}_{TU}^{(L)}=N_{p}\left(\mathcal{C}_{\eta_{\mathrm{fp}%
,k+1,j}^{(L)}}+\mathcal{C}_{\mathbf{C}_{\mathrm{fp},k+1,j}^{(L)}}\right).
\end{equation}
where the costs $\mathcal{C}_{\eta_{\mathrm{fp},k+1}}$ and $\mathcal{C}_{%
\mathbf{C}_{\mathrm{fp},k+1,j}^{(L)}}$ are equal to $\mathcal{C}_{\mathbf{A}%
^{(L)}} +\mathcal{C}_{\mathbf{f}^{(L)}} + 2D_{L}^2$ flops and $%
D_{L}^2(4D_{L}-1) $ flops, respectively.

\section{Computational complexity of the MPF technique developed in ref. 
\protect\cite{Djuric_2013} \label{app:CMPF}}

In this appendix a detailed analysis of the MPF complexity is illustrated.
The notation is the same as \cite{Djuric_2013}. In the following, $\mathcal{C%
}_{f_x}$ and $\mathcal{C}_{f_y}$ denote the cost due to the evaluation of
the functions $f_{x}({x}_{t-1},{u}_t)$ and $f_{y}({x}_t,{v}_t)$,
respectively. Moreover, similarly as \cite{Hoteit_2016}, it is assumed that
the computation of the inverse of any covariance matrix involves a Cholesky
decomposition of the matrix itself and the inversion of a lower or upper
triangular matrix.

1. \textbf{Measurement update}

The overall computational cost of this task is 
\begin{equation}
\mathcal{C}_{MU}=n\left(\mathcal{C}_{x_{i,t-1}^{(m)}}+\mathcal{C}%
_{w_{i,t}^{(m)}}+\mathcal{C}_{x_{i,t}^{(m)}} + \mathcal{C}_{R}(M)\right).
\end{equation}
Moreover, we have that: 1) $\mathcal{C}_{x_{i,t-1}^{(m)}}$ is equal $M L
(N_{f} - 1) $ flops; 2) $\mathcal{C}_{w_{i,t}^{(m)}}$ is equal to $ML (6%
\mathcal{C}_{f_y} +4d_{y}^3+21d_{y}^2+17d_{y})/6+2ML+2M-1 $ flops; 3) $%
\mathcal{C}_{x_{i,t}^{(m)}}$ is equal to $d_{x,i} (2M-1) $ flops; 4) $%
\mathcal{C}_{R}(M)$ denotes the total cost of the resampling step (that
involves a particle set of size $M$).

2. \textbf{Time update}

The overall computational cost of this task is 
\begin{equation}
\mathcal{C}_{TU}=n\mathcal{C}_{x_{i,t}^{(m)}},
\end{equation}
where the cost $\mathcal{C}_{x_{i,t}^{(m)}}$ is equal to $M(3\mathcal{C}%
_{f_x} + d_{x,i}^3+9d_{x,i}^2+5d_{x,i})/3 $ flops.\newline

\end{document}